\documentclass[a4paper]{gtart}
\gtart

\usepackage{amsthm}
\usepackage{amsgen}
\usepackage{amsmath,amssymb}
\usepackage{curvesls}
\usepackage{epic}
\usepackage[all]{xy}
\usepackage[dvips]{graphicx}
\usepackage{psfrag}
\usepackage[T1]{fontenc}
\usepackage[sc]{mathpazo}
\usepackage{color}
\usepackage{float}
\newtheorem{thm}{Theorem}[section]
\newtheorem{lem}[thm]{Lemma}
\newtheorem{cor}[thm]{Corollary}
\newtheorem{prop}[thm]{Proposition}

\newtheorem{question}[thm]{Question}
\newtheorem{constr}[thm]{Construction}

\theoremstyle{definition}
\newtheorem{defn}[thm]{Definition}

\newcommand{\RProj}{{\mathbb RP}}
\newcommand{\CProj}{{\mathbb CP}}
\newcommand{\C}{{\mathbb C}}
\newcommand{\Real}{{\mathbb R}}

\newcommand{\Zed}{{\mathbb Z}}

\newcommand{\Rat}{{\mathbb Q}}
\newcommand{\splice}{{\bowtie}}

\newcommand{\mubar}{{\overline{\mu}}}

\newcommand{\Ext}{{\mathrm{Ext}}}
\newcommand{\Hom}{{\mathrm{Hom}}}
\newcommand{\rank}{{\mathrm{rank}}}
\newcommand{\image}{{\mathrm{img}}}
\newcommand{\Det}{{\mathrm{Det}}}
\newcommand{\Id}{{\mathrm{Id}}}
\newcommand{\trace}{{\mathrm{tr}}}
\newcommand{\linkno}{{\mathrm{lk}}}
\newcommand{\spin}{{\mathrm{Spin}}}
\newcommand{\spinc}{{{\mathrm{Spin}}^c}}
\newcommand{\signature}{{\mathrm{sig}}}
\newcommand{\KR}[2]{{\langle {\mathrm {#1}}: {#2} \rangle}}
\newcommand{\lau}{{\Lambda}}
\newcommand{\fof}{{\mathbb Q(\Lambda)}}

\newcommand{\ttc}{{ }}

\begin{document}

\title{Embeddings of $3$-manifolds in $S^4$ from the point of view of the $11$-tetrahedron census}

\authors{Ryan Budney \\ Benjamin A. Burton}

\addresses{Mathematics and Statistics, University of Victoria, Canada \\
School of Mathematics and Physics,  University of Queensland, Australia}
\emails{rybu@uvic.ca \\ bab@maths.uq.edu.au}

\begin{abstract} 
This is a collection of notes on embedding problems for 3-manifolds. The main question explored is 
``which 3-manifolds embed smoothly in $S^4$?'' The terrain of exploration is the 
Burton/Martelli/Matveev/Petronio census of triangulated prime closed 3-manifolds built from 
11 or less tetrahedra. There are 13766 manifolds in the census, of which 13400 are orientable. 
Of the 13400 orientable manifolds, only 149 of them have hyperbolic torsion linking forms and 
are thus candidates for embedability in $S^4$. The majority of this paper is devoted to the 
embedding problem for these 149 manifolds. At present \ref{emb_num} are known to embed in $S^4$. 
Among the remaining manifolds, embeddings into homotopy $4$-spheres are constructed for 
\ref{emb_homolsph_num}. \ref{nonemb_num} manifolds in the list are known to not embed in 
$S^4$. This leaves \ref{unknown_num} unresolved cases, of which only \ref{lastgeometric} 
are geometric manifolds i.e. having a trivial JSJ-decomposition. 
\end{abstract}


\maketitle

\section{Introduction}\label{INTRODUCTION}

Given a smooth manifold $M$, let $ed(M)$ denote the minimum of all integers $n$ such that $M$ 
admits a smooth embedding into $S^n$.  The purpose of these notes is to get a sense for how 
difficult it is to determine if $ed(M)=4$, when $M$ is compact, boundaryless $3$-dimensional manifold. 

Whitney proved that $ed(M) \leq 2n$ for all $n$-manifolds $M$ by a combination of a general 
position/transversality argument and a double point creation and destruction process now 
called The Whitney Trick.  A basic argument using characteristic classes shows that 
$ed(\RProj^{2^k}) = 2\cdot 2^k$ for all $k$, and so Whitney's result is generally the 
best one can expect for arbitrary $n$ (see for example Theorem 4.5 and Corollary 11.4 
in \cite{MS}).  For $3$-manifolds, C.T.C. Wall improved on Whitney's result, showing every
compact $3$-manifold embeds in $S^5$ \cite{Wall}.  Thus, for closed $3$-manifolds distinct 
from $S^3$, the embedding dimension can be one of two possible numbers $ed(M) \in \{4,5\}$.  

Recently Skopenkov has given a complete isotopy classification of embeddings of $3$-manifolds 
into $S^6$ \cite{Skopenkov}.  At the other extreme, the question of which $3$-manifolds (with 
boundary) embed in $S^3$ is quite a difficult problem \cite{CGLS, Agol, MJ} although there is 
much known -- for example, consider the case of $M$ compact, orientable with boundary a 
collection of tori. If $M$ embeds in $S^3$, then there is another embedding of $M$ in $S^3$ 
so that it is the complement of a link \cite{Schubert}.  By a Wirtinger presentation, $\pi_1 M$ 
is generated by the conjugates of $n$ curves on $\partial M$ corresponding to the meridians of 
the $n$-component link.  By the resolution of the Poincar\'e conjecture, the converse is 
true -- simply fill $M$ along the curves in $\partial M$ to get a homotopy $3$-sphere.  
Although this is an `answer' it is rather difficult to implement in a computationally-effective way \cite{MJ}. 

The hope of this paper is that the `intermediate' question of whether or not a $3$-manifold 
embeds in $S^4$ is tractable. This is also problem 3.20 on Kirby's problem list.  The point of 
view of this paper is that there is no better way to discover than to get one's hands dirty. 
The census of prime $3$-manifolds which can be triangulated by $11$ or less tetrahedra \cite{BBurton} 
is chosen as a `generic supply' of test cases.  Of course, there is good reason to think this 
problem could be very difficult.  There are several significant, closely-related outstanding 
problems such as the Sch\"onflies problem, and the smooth Poincar\'e Conjecture in dimension 
$4$ which indicate possible pitfalls.  Sometimes in this paper embeddings of $3$-manifolds 
are constructed into homotopy $4$-spheres. Likely all the homotopy $4$-sphere constructed are 
the standard $S^4$ but we have not always determined this.   There are perhaps simpler 
obstacles to overcome -- at present in the literature there are no known examples of $3$-manifolds 
that embed smoothly in a homology $4$-sphere but do not embed in $S^4$. It is rather remarkable 
that all the obstructions used in this paper are obstructions to embedding into homology 
$4$-spheres, and at least so far they have largely sufficed to determine which $3$-manifolds 
embed in $S^4$. 

In Section \ref{Known results} a brief survey is given of known obstructions to a $3$-manifold 
embedding in $S^4$. Many useful techniques to construct embeddings in $S^4$ are also listed.

We apply the results from Section \ref{Known results} to the census of $3$-manifolds in Sections 
\ref{embeddable_man}, \ref{emb_htpy} and \ref{nonembeddable_man}.  To keep the paper a reasonable 
length, only the manifolds which pass the torsion linking form test (Theorem \ref{Hant}) are 
listed in these sections.

Section \ref{embeddable_man} describes embeddings in $S^4$ for the manifolds in the 
census which are known to embed in $S^4$.  

Section \ref{emb_htpy} describes embeddings into homotopy $4$-spheres of the manifolds in the
census that are known to embed in homotopy $4$-spheres -- these homotopy $4$-spheres are 
likely to be diffeomorphic to $S^4$ but this has not been determined.

Section \ref{nonembeddable_man} provides obstructions for the manifolds in the census which are 
known not to embed in $S^4$ (or any homology $4$-sphere). Manifolds which fail the torsion linking
form test (Theorem \ref{Hant}) are not listed as these are too numerous. Manifolds that fail
the torsion linking form test are available via the software Regina, see Section \ref{hyperbolicity}.

Section \ref{unknown_embeddable_man} lists the manifolds for which it is not yet known if 
they embed in $S^4$ or homology $4$-spheres.  Moreover, a list of computed obstructions is
provided.

Section \ref{hyperbolicity} provides sketches of some techniques used to compute various invariants
of the manifolds from the census.  If the reader ever gets lost in the notation used in the
tables, usually this section or Section \ref{Known results} is the appropriate place to look for clarification. 

Section \ref{qsec} contains various observations and comments on the data.  

Many of the obstructions and constructions present in this paper were described to the first author 
by Danny Ruberman.  Thanks to Brendan Owens and Sa\v{s}o Strle who kindly let us use their software 
to compute the $d$-invariant of Seifert fibred rational homology spheres.  Thanks also to Jonathan 
Hillman, Ahmed Issa, Gregor Masbaum, Peter Landweber, Lee Rudolph, Ronald Fintushel, Ronald Stern, 
Ian Agol, Scott Carter, Nathan Dunfield, Jeff Weeks, Peter Teichner, Tom Goodwillie and Mike Freedman 
for their suggestions and/or encouragement (whether they remember it or not).  A paper such as this 
requires immense amounts of time for thousands of hand and computer-aided computations.   The first author 
would especially like to thank the Max Planck Institute for Mathematics (Bonn) and the Institut des 
Hautes \'Etudes Scientifiques for giving him the freedom to initiate this open-ended project.  
Thanks also to IPMU (Tokyo) and KITP (Santa Barbara) for hosting the first author.

\section{Obstructions and embedding constructions}\label{Known results}

There are only a few completely general obstructions to a closed $3$-manifold embedding in $S^4$.  
The first is of course orientability, coming from the generalized Jordan Curve Theorem. There 
are no other tangent-bundle derived obstructions since the tangent bundle of an orientable 3-manifold 
is trivial (Stiefel's Theorem) \cite{Kirby}. A powerful and easy-to-compute obstruction comes 
from the torsion linking form of a $3$-manifold.  

\begin{defn}\label{pdrev} (Torsion Linking Form) In a compact, boundaryless oriented $n$-manifold $M$ there is
a canonical, natural isomorphism (Poincar\'e duality) 
$$H_i(M,\Zed) \simeq H^{n-i}(M,\Zed) \hskip 3mm \forall i \in \{0,1,\cdots,n\}.$$
This is a natural short exact sequence (the homology-cohomology Universal Coefficient Theorem)
$$0 \to \Ext_\Zed(H_{i-1}(M,\Zed),\Zed) \to H^i(M,\Zed) \to \Hom(H_i(M,\Zed),\Zed) \to 0 \hskip 3mm \forall i \in \{0,1,\cdots,n\} $$
and a canonical isomorphism
$$\Ext_\Zed(H_i(M,\Zed),\Zed) \simeq \Hom_\Zed(\tau H_i(M,\Zed),\Rat/\Zed) \hskip 3mm \forall i \in \{0,1,\cdots,n\}$$
where $\tau H_i(M,\Zed)$ is the subgroup of torsion elements of $H_i(M,\Zed)$. This gives
two duality pairings on the homology of $M$, the `intersection product'
and the `torsion linking form' respectively:
$$fH_i(M,\Zed) \otimes fH_{n-i}(M,\Zed) \to \Zed \hskip 10mm \tau H_i(M,\Zed) \otimes 
 \tau H_{n-i-1}(M,\Zed) \to \Rat/\Zed$$
where $fH_i(M,\Zed) = H_i(M,\Zed) / \tau H_i(M,\Zed)$ is the `free part' of $H_i(M,\Zed)$.
\end{defn}

See the discussion preceding Figure \ref{dualbits} in Section \ref{hyperbolicity} for details on 
how one computes the torsion linking form for a triangulated manifold, in practice.  In short, 
given $[x] \in \tau H_i(M)$ and $[y] \in \tau H_{n-i-1}(M)$, let $ky = \partial Y$, i.e. assume
$[y]$ is $k$-torsion, then the linking form is defined as $\langle [x], [y] \rangle = \frac{x \pitchfork Y}{k}$
where $x \pitchfork Y$ indicates the transverse signed intersection number. 

\begin{thm}\cite{Han, KK}\label{Hant} (Hantzsche Test)
If $M$ is a compact, boundaryless, connected, oriented $3$-manifold which embeds in
a homology $S^4$ then there is a splitting $\tau H_1(M,\Zed) = A \oplus B$, 
inducing a splitting 
$$\Hom_\Zed(\tau H_1(M,\Zed),\Rat/\Zed) \simeq \Hom_\Zed(A,\Rat/\Zed) \times \Hom_\Zed(B,\Rat/\Zed)$$
which is reversed by Poincar\'e duality, in the sense that the P.D. isomorphism
$$\tau H_1(M,\Zed) \to \Hom_\Zed(\tau H_1(M,\Zed),\Rat/\Zed)$$
restricts to isomorphisms $A \to \Hom_\Zed(B,\Rat/\Zed)$ and 
$B \to \Hom_\Zed(A,\Rat/\Zed)$.  This uses the convention that $\Hom_\Zed(A,\Rat/\Zed)$ 
is the submodule of $\Hom_\Zed(\tau H_1(M,\Zed), \Rat/\Zed)$ which is zero on $B$, similarly
$\Hom_\Zed(B,\Rat/\Zed)$ is the submodule which is zero on $A$.
\begin{proof}
$M$ separates the homology $4$-sphere into two manifolds, call them $V_1$ and
$V_2$, $S^4 = V_1 \cup_{M} V_2$.   Let $A=\tau H_1(V_1,\Zed)$ and $B=\tau H_1(V_2,\Zed)$,
when we have the isomorphism $A \oplus B \simeq H_1(M,\Zed)$ by the Mayer-Vietoris sequence for 
$S^4 = V_1 \cup_{M} V_2$. Let $\{i,j\} = \{1,2\}$. If a homology class in $\tau H_1(M)$ comes
from $\tau H_1(V_i)$ then it must bound a $2$-cycle in $V_j$. Thus the torsion linking form is 
zero on $\tau H_1(V_i) \otimes \tau H_1(V_i)$, giving the result.  
\end{proof}
\end{thm}

An immediate corollary of Theorem \ref{Hant} is that the only
lens space that admits a smooth embedding into $S^4$ is $S^3$. Our convention
is that a lens space is a manifold quotient of $S^3$ by a group of isometries, 
so we exclude $S^1 \times S^2$.
Kawauchi and Kojima call torsion linking forms which have such 
a splitting `hyperbolic' \cite{KK}.  Kawauchi and Kojima's test for hyperbolicity of the
torsion linking form has been implemented by the author in the freely-available 
open-source software package `Regina' \cite{BBurton}.  

%
%

As stated in the abstract, 
there are only $149$ manifolds in the census with hyperbolic torsion linking
forms, and they are listed in Sections \ref{embeddable_man}, \ref{emb_htpy}, 
\ref{nonembeddable_man} and \ref{unknown_embeddable_man}.  
Since the hyperbolicity computation plays a significant role in this
paper, a sketch of the algorithm is given in Section \ref{hyperbolicity}.

In general, if a $3$-manifold $M$ embeds in a homology $4$-sphere $\Sigma^4$, 
$V_1 \cup_M V_2 = \Sigma^4$. The argument in the proof of Theorem \ref{Hant} gives us 
(for $\{i,j\} = \{1,2\}$):

\begin{align*}
H_1 V_i &\simeq fH_1 V_i \oplus \Hom_\Zed(\tau H_1 V_j, \Rat/\Zed) \\ 
H_2 V_i &\simeq \Hom_\Zed (fH_1 V_j, \Zed)  \\
H_3 V_i &\simeq *.
\end{align*}

If $M$ is a rational homology sphere, the manifolds $V_1$ and $V_2$ are rational homology balls.  
If $M$ is a rational homology $S^1 \times S^2$, one of $V_1, V_2$ is a rational homology
$S^1 \times D^3$, and the other a rational homology $D^2 \times S^2$.  
If $H_1 M \simeq \Zed^2 \oplus \tau H_1 M$ then there are two possibilities: 
in the first case, one would be a rational genus two $1$-handlebody 
$S^1 \times D^3 \#_\partial S^1 \times D^3$ and the other a rational genus two $2$-handlebody 
$(D^2 \times S^2) \#_{\partial} (D^2 \times S^2)$, 
in the second case both manifolds would be rational $(S^1 \times D^3) \#_\partial (D^2 \times S^2)$. 

By and large, these complications do not
come up much in the census as the majority (13173 of 13766) are rational homology spheres. 
There are only 201 rational homology $S^1 \times S^2$ manifolds in the census. There are 25 manifolds
in the census that have $fH_1(M) \simeq \Zed^2$, and there is only one manifold in the census with
$fH_1(M) \simeq \Zed^3$, the manifold $S^1 \times S^1 \times S^1$.  There are no manifolds in the
census with $\rank(H_1 M)>3$. Thus, intersection forms such as $H_2 M \otimes H_2 M \to H_1 M$
gives no useful obstruction to census $3$-manifolds embedding in $S^4$.  

Kawauchi developed an obstruction to a rational homology $S^1 \times S^2$ bounding a rational homology
$S^1 \times D^3$, which will be described below.

\begin{defn}(Alexander Polynomial)
If $h : H_1(M,\Zed) \to \Zed$ is an epimorphism, let $M_h \to M$ denote the normal
abelian covering space corresponding to $h$, and let $h$ play a double-role as 
the corresponding generator of the group of covering transformations.
Consider $H_1(M_h,\Rat)$ to be a module over $\lau \equiv \Rat[\Zed]$ (the group ring of 
the integers $\Zed$ with coefficients in the rationals $\Rat$), 
where the action of $\Zed$ on $M_h$ is generated 
by the covering transformation $h$. Notice that $\Rat[\Zed]$ is isomorphic to a 
Laurent polynomial ring $\Rat[h^{\pm 1}]$, which is a principal ideal domain.  
By the classification of finitely-generated modules over PIDs, 
$H_1(M_h,\Rat) \simeq \lau^k \oplus (\oplus_{p \in P} \lau/p)$ for various non-zero 
polynomials $P$.  The order ideal of the $\lau$-torsion submodule of
$H_1(M_h,\Rat)$ is called the Alexander polynomial of $h$, and will be denoted 
$\Delta(h) = \prod_{p \in P} p \in \Rat[h^{\pm 1}]$.  Since 
it is representing an ideal, it is only well-defined up to multiplication by a unit.
We will use the notation $\fof$ for the field of fractions of $\lau$. 
\end{defn}

We use the symbol $\equiv$ to denote either a definition or a canonical identification, while $\simeq$ 
denotes abstract isomorphism. 

When $M$ is compact, orientable and boundaryless, Poincar\'e duality (of the Blanchfield
variety -- see for example \cite{LinkInv})
and basic linear algebra provides isomorphisms 
$$\tau_\lau H_1(M_h,\Rat) \simeq \tau_\lau \overline{H^2(M_h,\Rat)} \simeq
  \overline{\Ext_{\lau}(H_1(M_h,\Rat),\lau)}$$
$$ \hskip 14mm \simeq \overline{\Hom_{\lau}(\tau_\lau H_1(M_h,\Rat),\fof/\lau)} $$
where cohomology is `cohomology with compact support.'
The inclusion $\lau \subset \fof$ is the submodule consisting of elements whose 
denominator is $1$. Given a $\lau$-module $A$, $\overline{A}$ indicates the conjugate
$\lau$-module -- as a $\Rat$-vector space it is identical to $A$, but the action of
$\Zed$ on $A$ is the inverse action.  This statement is the $\lau$ analogue of the isomorphisms in
Definition \ref{pdrev}. Since $\lau$ is a PID, 
$\tau_\lau H_1(M_h,\Rat)$ has a diagonal presentation matrix, thus there
is a (not natural) isomorphism between
$\tau_\lau H_1(M_h,\Rat)$ and $\Hom_{\lau}(\tau_\lau H_1(M_h,\Rat),\fof/\lau)$.
Thus, the Alexander polynomial is symmetric $\Delta(h) = \Delta(h^{-1})$.

Notice if $h : H_1(M,\Zed) \to \Zed$ is an epimorphism, and if $M$ embeds in
a homology $S^4$, then one can write $S^4$ as a union $V_1 \cup_M V_2$ and so the 
homomorphism $h$ factors as a composite
$$\xymatrix{ H_1(M,\Zed) \ar[r]^-h \ar[d] & \Zed \\
             H_1(V_i,\Zed) \ar[ur] & }$$
for some $i \in \{1,2\}$.

\begin{thm}\label{KawCond}\cite{Kawauchi} (Kawauchi Test)
If $M$ is a rational homology $S^1 \times S^2$ with $h : H_1(M,\Zed) \to \Zed$ onto,
and if $M$ admits an embedding into a homology $S^4$ then 
$\Delta(h) = f(h)f(h^{-1})$ for some Laurent polynomial $f(h) \in \Rat[h^{\pm 1}]$.
\begin{proof}
Let $V$ be the rational homology $S^1 \times D^3$ bounding $M$, as above. 
Consider the Poincar\'e Duality long exact sequence of the pair
$(V_h,M_h)$:
$$\xymatrix{ \ar[r]^-{j_*} & H_2(V_h,M_h,\Rat) \ar[r]^-\partial \ar[d]^-{PD} & H_1 (M_h,\Rat) \ar[r]^-{i_*} \ar[d]^-{PD} & 
                   H_1 (V_h,\Rat) \ar[d]^-{PD} \ar[r]^-{j_*} &  \\
            \ar[r]^-{j^*} & \overline{H^2 (V_h,\Rat)} \ar[r]^-{i^*} & \overline{H^2 (M_h,\Rat)} \ar[r]^-{\delta} & 
                   \overline{H^3 (V_h,M_h,\Rat)} \ar[r]^-{j^*} &  
}$$
The next step is to show all six $\lau$-modules in the above exact ladder are $\lau$-torsion. 
First consider $H_2(V_h,M_h,\Rat)$. By the Poincar\'e Duality isomorphism
$H_2(V_h,M_h,\Rat) \simeq \overline{H^2(V_h,\Rat)}$. The Universal Coefficient Theorem reduces
this to showing that $H_2(V_h,\Rat)$ is a $\lau$-torsion module.  Consider the long exact sequence
$$ \xymatrix{ \ar[r] & H_2(V_h,\Rat) \ar[r]^{(t-1)} & H_2(V_h,\Rat) \ar[r]^{p_*} & H_2(V,\Rat) \ar[r]^{\partial} & H_1(V_h,\Rat) \ar[r] & }$$  
Where `$(t-1)$' indicates multiplication by $(t-1)$, and $p : V_h \to V$ is the covering projection. 
$H_2(V,\Rat)=0$ therefore multiplication by $(t-1)$ is onto $H_2(V_h,\Rat)$, thus $H_2(V_h,\Rat)$ 
is $\lau$-torsion.  Similarly, $H_1(V_h,\Rat)$ is $\lau$-torsion.  
$H_1(M_h,\Rat)$ is an extension of a quotient of $H_2(V_h,M_h,\Rat)$, and
a submodule of $H_1(V_h,\Rat)$, so it is also torsion.

Poincar\'e Duality combined with the Universal Coefficient Theorem gives us isomorphisms
of the three short exact sequences:

$$\xymatrix{ 0 \ar[r] & \image(\partial) \ar[r] \ar[d]^{PD} & H_1(M_h,\Rat) \ar[r] \ar[d]^{PD} & 
                               \image(i_*) \ar[r] \ar[d]^{PD} & 0  \\
             0 \ar[r] & \overline{\image(i^*)} \ar[r] \ar[d]^{UCT} & \overline{H^2 (M_h,\Rat)} \ar[r] \ar[d]^{UCT} &
                               \overline{\image(\delta)} \ar[r] \ar[d]^{UCT} & 0 \\
             0 \ar[r] & \overline{\image(\Ext(i_*))} \ar[r] & \overline{\Ext(H_1(M_h,\Rat),\lau)} \ar[r] & 
                              \overline{\image(\Ext(\partial))} \ar[r] & 0 }$$
where
$\Ext(i_*) : \Ext(H_1(V_h,\Rat),\lau) \to \Ext(H_1(M_h,\Rat),\lau)$ and 
$\Ext(\partial) : \Ext(H_1(M_h,\Rat),\lau) \to \Ext(H_2(V_h,M_h,\Rat),\lau)$ are the 
$\Ext(\cdot,\lau)$-duals to $i_*$ and $\partial$ respectively.

The remainder follows from the following well-known lemma.
\end{proof}
\end{thm}

\begin{lem}
\begin{itemize}
\item Given a short exact sequence of finitely generated torsion $\lau$-modules
$0 \to A \to B \to C \to 0$, the order ideal of $B$ is the product of the order ideals of 
$A$ and $C$ respectively.  
\item If $f : A \to B$ is a homomorphism of finitely generated torsion $\lau$ modules then
$\image(f)$ and $\image(\Ext(f))$ have the same order ideals, where
$\Ext(f) : \Ext(B,\lau) \to \Ext(A,\lau)$ is induced from $f$.
\end{itemize}
\end{lem}

Theorem \ref{KawCond} can be generalized to an obstruction for a $3$-manifold $M$ to bound a $4$-manifold
$W$ provided $H_1(M,\Rat) \to H_1(W,\Rat)$ is onto with kernel of dimension at most $1$, 
$\rank (H_1 M) > 0$, and $H_2 W=0$ \cite{Kawauchi}. 
Unfortunately, this is not quite an obstruction to a $3$-manifold embedding in a homology $4$-sphere $\Sigma$,
provided $\rank (H_1 M) > 1$.  Take for example a $3$-manifold $M$ with $\rank (H_1 M)=2$.  If 
$\Sigma = V_1 \sqcup_M V_2$, this obstruction could be used to argue that neither $V_1$ nor $V_2$ are
rational homology $(S^1\times D^3) \#_\partial (S^1 \times D^3)$, but it can't be used to rule out the
possibility that $V_1$ and $V_2$ are rational homology $(S^1\times D^3) \#_\partial (S^2\times D^2)$'s.

Let $M$ be a rational homology $S^1 \times S^2$. As with knots, the Alexander polynomial can be 
defined integrally in terms of $H_1(M_h,\Zed)$, giving an integral normalization of 
$\Delta(h) \in \Zed[\Zed]$ (the group-ring of the integers with coefficients in $\Zed$).  One definition is
to let $\tau_\Zed H_1 (M_h,\Zed)$ denote the $\Zed$-torsion submodule of $H_1 (M_h,\Zed)$, and
$\tau_{\Zed[\Zed]} H_1 (M_h,\Zed)$ denote the $\Zed[\Zed]$-torsion submodule of $H_1 (M_h,\Zed)$. 
Let $f_{\Zed} H_1 (M_h,\Zed)$ be the maximal free quotient $\Zed$-module of $\tau_{\Zed[\Zed]}H_1(M_h,\Zed)$ i.e. 
$f_{\Zed} H_1(M_h,\Zed) = \tau_{\Zed[\Zed]} H_1 (M_h,\Zed) / \tau_{\Zed} H_1 (M_h,\Zed)$.
Define the Alexander polynomial of $h$ to be $\Delta(h) = \Det(hI - h_*)$, where $h_*$ is the automorphism of 
$f_{\Zed} H_1 (M_h,\Zed)$ (thought of as a finitely-generated free $\Zed$-module),
$I$ the identity automorphism, and $h$ is a variable effectively making the expression $hI - h_*$ a matrix
with entries in $\Zed[h^{\pm 1}] \equiv \Zed[\Zed]$. 

The group $H_1 (M_h,\Zed)$ is $\Zed$-torsion free.  This follows from Poincar\'e duality, which
when followed by Universal Coefficients gives the Farber-Levine isomorphism \cite{LinkInv}
$$\tau_\Zed H_1(M_h,\Zed)
\simeq \Hom_\Zed(\tau_\Zed H_0(M_h,\Zed), \Rat/\Zed) = 0.$$ Consider the homology long exact sequence induced from the short exact sequence
$$\xymatrix{ 0 \ar[r] & C_* (M_h,\Zed) \ar[r]^{t-1} & C_* (M_h,\Zed) \ar[r] & C_*(M,\Zed) \ar[r] & 0}.$$
This allows us to compute $\Delta(h=1)$ (the Alexander polynomial evaluated at $h=1$) 
as $\Delta(h=1) = \pm |\tau_\Zed H_1(M,\Zed)|$.  
This condition together with the symmetry of the Alexander polynomial provide redundancies
that are helpful  when doing hand computations of the Alexander polynomial. 

There are further obstructions to a rational homology $S^1 \times S^2$ embedding
in a homology $S^4$, called signature invariants. As we have seen above there is a 
canonical isomorphism of $\lau$-modules
$$ H_1 (M_h,\Rat) \simeq \overline{\Hom_{\lau}(H_1 (M_h,\Rat), \fof/\lau ) }$$
which we think of as a sesquilinear duality pairing
$$ \langle \cdot, \cdot \rangle : H_1(M_h,\Rat) \times H_1(M_h,\Rat) \to \fof/\lau. $$
Here is how one computes the pairing.  Let $[v], [w] \in H_1(M_h,\Rat)$ be homology classes, with
$v,w \in C_1(M_h,\Rat)$ the corresponding cycle representatives.  Since they are $\lau$-torsion
classes, let $A_v, A_w \in \lau$ be non-zero such that $A_v v = \partial S_v$ and $A_w w = \partial S_w$. Then
 
$$ \langle [v],[w] \rangle = \frac{1}{A_v} \sum_{i \in \Zed} \left( S_v \pitchfork h^i w \right) h^i = 
  \frac{1}{\overline{A_w}} \sum_{i \in \Zed} \left( v \pitchfork h^i S_w \right) h^i \in \Rat[h^\pm] \equiv \lau.$$ 

$\overline{A_w}$ indicates we are taking the conjugate polynomial (conjugation is
the automorphism of $\lau \equiv \Rat[h^{\pm}]$ induced by the non-trivial automorphism
of $\Zed$, or equivalently by the operation on polynomials $h \longmapsto h^{-1}$). 
The symbol $\pitchfork$ indicates we are taking the oriented intersection number -- i.e. one first perturbs
the chains to be transverse and then takes the signed intersection number. 
That the pairing $\langle \cdot,\cdot\rangle$ is sesquilinear means that it is $\Rat$-linear in both
variables and $h\langle x,y\rangle = \langle hx,y\rangle = \langle x,h^{-1}y\rangle$ for all
$x,y \in H_1 (M_h,\Rat)$.  Moreover, $\langle x,y\rangle = \overline{\langle y,x\rangle}$ for all $x,y \in H_1 (M_h,\Rat)$, 
where the conjugation is the involution of $\fof/\lau$ induced by conjugation on
$\lau$. 
From the pairing $\langle \cdot,\cdot \rangle$ we construct an anti-symmetric pairing
$\left[ \cdot,\cdot \right] : H_1(M_h,\Rat) \times H_1(M_h,\Rat) \to \Rat$ by composing with 
the `Trotter trace' function $\trace : \fof/\lau \to \Rat$, i.e. 
$\left[x,y\right] = \trace(\langle x,y\rangle)$. See page 182 of \cite{Trotter} for details on the 
trace function. In brief:
\begin{itemize}
\item[(a)] $\trace$ is a $\Rat$-linear function such that $\trace(\overline{x}) = -\trace(x)$ for all $x \in \fof/\lau$.
\item[(b)] Given $p,q \in \lau$ where $q$ is not a unit
nor divisible by $1-h$, and assuming the lowest (resp. highest) degree non-zero coefficient of 
$p$ has degree $\geq$ (resp. $\leq$) the lowest (resp. highest) degree non-zero coefficient of
$q$ (say, via the division algorithm), $\trace(p/q)$ is defined to be the derivative evaluated at $1$, 
$\trace(p/q) = (p/q)'(1)$.  
\item[(c)] If $q$ is a unit or is a power of $1-h$, let $\trace(p/q) = 0$.
\item[(d)] $\trace$ is defined on $\fof/\lau$ by extending the definitions (b) and (c) linearly.
\item[(e)] An essential property of the Trotter trace is that provided we're in case (b)
and that the highest-order non-zero term of $p$ is strictly smaller than the highest-order non-zero term for $q$, then $\trace((h-1)p/q) = (p/q)(1)$. 
\end{itemize}
From this it follows that composition with the Trotter trace gives an isomorphism
$$ \Hom_{\lau} (H_1(M_h,\Rat), \fof/\lau) \to \Hom_{\Rat} (H_1(M_h,\Rat), \Rat).$$
Thus, the pairing $\left[\cdot,\cdot\right]$ is non-degenerate, anti-symmetric and multiplication by $h$ is
an isometry $\left[ hx,hy\right] = \left[ x,y\right]$.   We construct a symmetric bilinear form
$H_1 (M_h,\Rat) \times H_1 (M_h,\Rat) \to \Rat$ via the formula
$\left\{ x,y \right\} = [x,ty] + [y,tx]$.  Notice that this symmetric form can potentially be degenerate: $[x,ty] + [y,tx] = 0$ if and only if $[x,(t^2-1)y] = 0$.  Assume $x \neq 0$ is fixed.  Since multiplication by $t-1$ is an isomorphism on $H_1 (M_h,\Rat)$, $[x,(t^2 -1)y] = 0$ for all $y \in H_1 (M_h,\Rat)$ if and only if $[x,(t+1)y] = 0$ for all $y$.  Therefore if we restrict $\left\{ \cdot, \cdot \right\}$ to the maximal $\lau$-submodule of $H_1(M_h,\Rat)$ on which multiplication by $t+1$ is an isomorphism, we get a non-degenerate symmetric form. Let $\sigma_h \in \Zed$ be the signature of this form.  Let $p$ be any prime symmetric factor of $\Delta(h)$.  By further restricting the above symmetric form to the submodule killed by a power of $p$, we get further signature invariants $\sigma_{p,h} \in \Zed$, called Milnor signature invariants.  These are closely related to Tristram-Levine invariants \cite{LivingstonSurv, LinkInv}. The relations among these signature invariants appears in slightly different form in \cite{K, Kawauchi, Erle}.

\begin{thm}\label{TLthm} (Signature Test)
If $M$ is a rational homology $S^1 \times S^2$ and if $M$ embeds in a homology $S^4$, then
all the above signature invariants are zero.  
\begin{proof}
The proof of Theorem \ref{KawCond} gives a commuting ladder
$$\xymatrix{ 0 \ar[r] & \image(\partial) \ar[r]^-{\partial} \ar[d] & H_1(M_h,\Rat) \ar[r]^-{i_*} \ar[d] & 
                               \image(i_*) \ar[r] \ar[d] & 0  \\
             0 \ar[r] & \overline{\image((i_*)^*)} \ar[r]^-{(i_*)^*} & \overline{\Hom_{\lau}(H_1(M_h,\Rat),\fof/\lau)} \ar[r]^-{\partial^*} & 
                              \overline{\image(\partial^*)} \ar[r] & 0 }$$
where the upper stars indicate $\Hom_{\lau}(\cdot,\fof/\lau)$-duals. 
Thus, the domain of the form $\left\{ \cdot,\cdot \right\}$ splits into two subspaces of
equal dimension, and the form is zero on one of them.  For a non-degenerate form this can happen
if and only if the signature is zero.
\end{proof}
\end{thm}

There are a few obstructions related to particular families of manifolds. 
For the geometric 3-manifolds among the geometries: $S^3$, $S^2$-fibre,
$\mathbb E^3$, $Sol$ and $Nil$, Crisp and Hillman \cite{CrispH} 
computed precisely which of these manifolds embed in $S^4$. 
They do this by a combination of the above obstructions together
with a new obstruction derived as a generalization of 
the Massey-Whitney Theorem on the normal Euler class of 
$2$-manifolds in homology $4$-spheres. 

Let $E$ be the total space of a $D^2$-bundle $p : E \to \Sigma$ over a closed 
surface $\Sigma$.  Let $q : \partial E \to \Sigma$ be the corresponding $S^1$-bundle. 
The Whitney class $W_2(q) \in H^2(\Sigma, \mathcal B) \simeq \Zed$ 
is {\it the} obstruction to the existence of an everywhere non-zero section of the
bundle $p : E \to \Sigma$. $W_2(q)$ is an element of the 2nd cohomology group of $\Sigma$ with 
coefficients in the bundle of groups $\mathcal B = \{(s,\pi_1 q^{-1} (s)) : s \in \Sigma\}$.

\begin{thm}\label{chthm}(Whitney-Massey-Crisp-Hillman) \cite{Whitney, Massey, CrispH}
The total space of a disc bundle $p : E \to \Sigma$ embeds in $S^4$ (equivalently, a homology $S^4$) if and only if 
\begin{itemize}
\item $W_2 = 0$, provided $\Sigma$ is orientable
\item $W_2 \in \{2\chi -4, 2\chi, 2\chi + 4, \cdots, 4-2\chi\}$ if $\Sigma$ is non-orientable, 
where $\chi$ is the Euler characteristic of $\Sigma$.
\end{itemize}
A circle bundle over a surface embeds in $S^4$ (equivalently a homology $4$-sphere) if and only if
\begin{itemize}
\item $W_2 \in \{-1,0,1\}$ provided $\Sigma$ is orientable
\item $W_2 \in \{2\chi -4, 2\chi, 2\chi + 4, \cdots, 4-2\chi\}$ provided $\Sigma$ is non-orientable.
\end{itemize}
Here $\chi \equiv \chi(\Sigma)$ is the Euler characteristic of the surface $\Sigma$, and $W_2$
is the Whitney class of the associated disc bundle (i.e. the obstruction to a section of the
circle bundle). 
\end{thm}
 
Circle bundles over surfaces are Seifert fibred manifolds with no 
singular fibres.  I will use the notation of Regina \cite{BBurton} which is
consistent with Orlik's unnormalized Seifert notation \cite{Orlik}.
A circle bundle over a surface $\Sigma$ with Euler number $W_2 = k$  is
denoted SFS$\left[ \Sigma : k \right]$.
Whitney constructed all the above embeddings of $D^2$-bundles over surfaces
\cite{Whitney} and conjectured it was the complete list of $D^2$-bundles
that embed in $S^4$. Massey went on to prove his conjecture \cite{Massey}.
Crisp and Hillman proved the extension for $S^1$-bundles over surfaces \cite{CrispH}.  

The proof that $W_2=0$ when $\Sigma$ is orientable follows from the observation that 
$W_2$ is the self-intersection number of $\Sigma$ in $S^4$, and that 
$\Sigma$ can be isotoped off itself in $S^4$. When $\Sigma$ is non-orientable, the same 
observation tells us that $W_2$ is even. To get the restriction 
$W_2 \in \{2\chi -4, 2\chi, 2\chi + 4, \cdots, 4-2\chi\}$ Massey employed the 
$G$-signature Theorem to show that $W_2$ is the signature of a certain form. Precisely, 
let $X$ be the
$\Zed_2$-branched cover of $S^4$ branched over $\Sigma$ corresponding to the non-trivial
element of $H_1(S^4 \setminus \Sigma,\Zed) \simeq \Zed_2$.  The $G$-signature Theorem
states that the Euler class of $\Sigma$ in $X$ is the signature of the form
$\langle x, T_* y\rangle$ on $H_2(X,\Rat)$ where $T : X \to X$ is the covering
transformation and $\langle \cdot,\cdot \rangle$ is the intersection product.
The result follows from the computations $H_2(X,\Rat) \simeq \Rat^{2-\chi(\Sigma)}$,
 $T_* = -\Id_{H_2(X,\Rat)}$ and that the Euler class of $\Sigma$ in $S^4$ is twice
that of $\Sigma$ in $X$. 

In the $S^1$-bundle case, with $\Sigma$ orientable the torsion linking
form is the appropriate embedding obstruction.  When $\Sigma$ is non-orientable, and
$M$ is an $S^1$-bundle over $\Sigma$, the torsion linking form test tells us that $W_2$ 
must be even.  Crisp and Hillman generalized \cite{CrispH} the above argument of Massey's.  
Since $W_2$ is even, $H_1(M,\Zed) \simeq \Zed^{g-1} \oplus \Zed_2 \oplus \Zed_2$ with
one of the $\Zed_2$ factors being generated by the fundamental class of the fibre. 
So if $M$ embeds in $S^4$, we have $S^4 = V_1 \cup_M V_2$ and so $H_1 V_1 \oplus H_1 V_2 \simeq
H_1 M$, and so the $\Zed_2$-summand corresponding to the fibre inclusion belongs
(WLOG) to $H_1 V_1$.  Let $W$ be $V_1$ union the $D_2$-bundle over $\Sigma$ with Euler class
$W_2$.  Let $W'$ be the $\Zed_2$-branched cover of $W$ branched over $\Sigma$, and apply
the $G$-signature Theorem as in the previous case.

Crisp and Hillman make similar but increasingly complex applications of the 
$\Zed_2$-signature Theorem as formulated in \cite{JO} to get further obstructions to the 
embedding of Seifert-fibred and $Sol$ manifolds. The idea being to use the homology of
$M$ to construct $2$-sheeted covering spaces $\tilde M$ of $M$, and to attach to it the associated
covers of $V_1$ or $V_2$, or some associated $\Zed_2$-space whose boundary is $\tilde M$ and for
which the fixed point set is understood. See Proposition 1.2 and Theorem 1.4 of
\cite{CrispH}.

A link $L \subset S^3$ is said to be {\it slice} if there is a manifold $D \subset D^4$ such that
$\partial D = L$ and $D$ is diffeomorphic to a disjoint union of discs $D^2$. $D$ is called
{\it slice discs} for $L$.

\begin{constr}\label{surgical} ({\it 0-Surgical Embeddings}):
Let $M$ be a $0$-surgery along a link $L \subset S^3$ where $L$ is the union of two
links $L = L_1 \cup L_2$ such that $L_i$ is smoothly slice for $i \in \{1,2\}$. Then
$M$ admits a smooth embedding into $S^4$.

\begin{figure}[H]
{\psfrag{d1}[tl][tl][0.8][0]{$D_1$}
\psfrag{d2}[tl][tl][0.8][0]{$D_2$}
\psfrag{s3}[tl][tl][1][0]{$S^3$}
\psfrag{s4}[tl][tl][1][0]{$S^4$}
\psfrag{M}[tl][tl][1][0]{$M$}
$$\includegraphics[width=8cm]{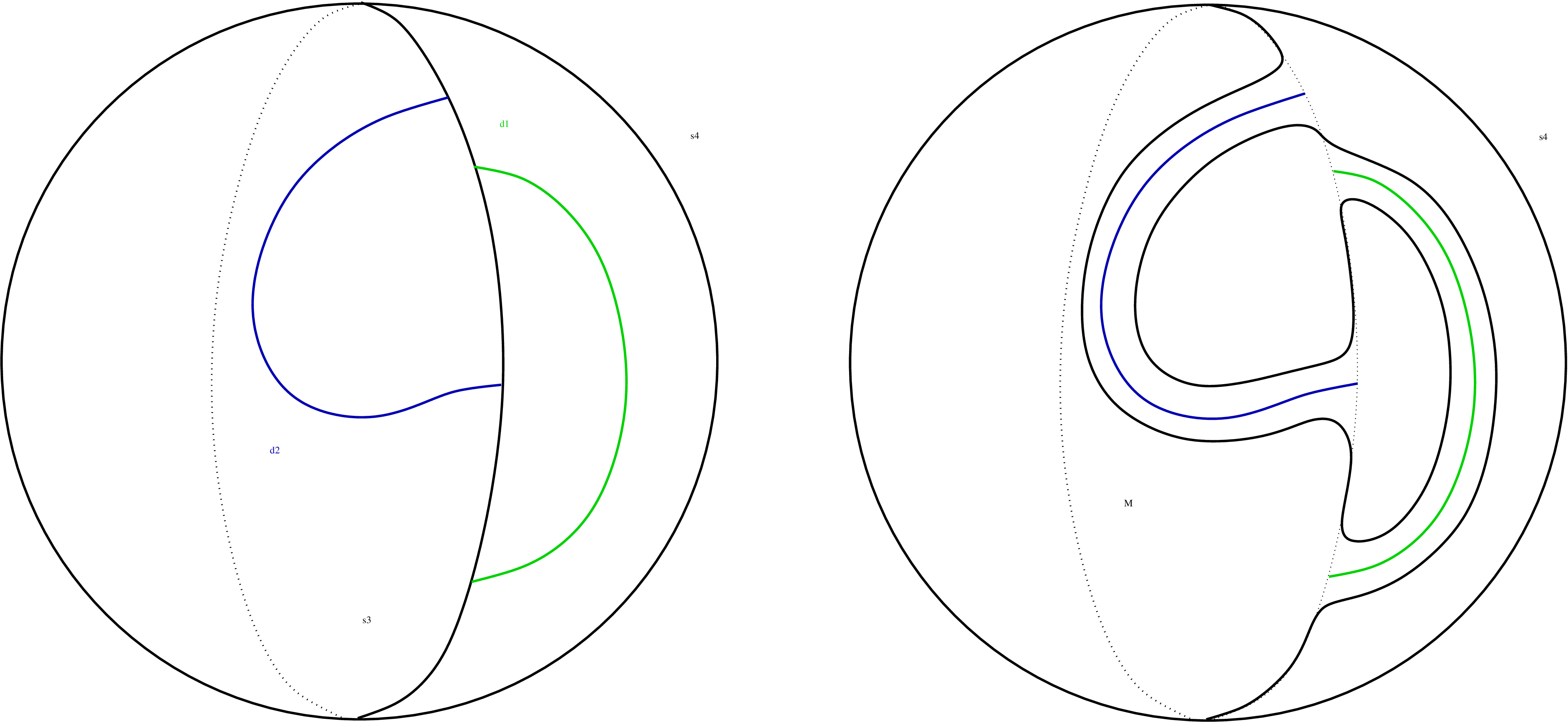}$$ }
\caption{\label{0-surg-emb} A $0$-surgical embedding.}
\end{figure}

\begin{proof}
The idea of the proof is to consider $S^4$ as the union of two $4$-balls, separated via a
great $3$-sphere.  Let $D_1$ be a collection of slice discs in the first hemi-sphere whose
boundary is $L_1$, and let $D_2$ be a collection of slice discs in the second hemi-sphere
whose boundary is $L_2$.  Then $M$ can be obtained by an embedded surgery on the
great $3$-sphere along the discs $D_1 \cup D_2$, see Figure \ref{0-surg-emb}.
\end{proof}
\end{constr}

Some examples of links which are the disjoint union of two slice links are: 
the Hopf link ($M = S^3$), 
the Whitehead link ($M=S^1 \ltimes_{\begin{pmatrix}  1 & 1 \\ 0 & 1\end{pmatrix}} (S^1\times S^1)$), 
and the Borromean rings ($M=S^1\times S^1 \times S^1$) \cite{CrispH}. 

\begin{figure}[H]
{\psfrag{d4}[tl][tl][1][0]{$D^4$}
\psfrag{d+p}[tl][tl][0.8][0]{$D_1'$}
\psfrag{d+}[tl][tl][0.8][0]{$D_1$}
\psfrag{d-}[tl][tl][0.8][0]{$D_{-1}$}
\psfrag{d-p}[tl][tl][0.8][0]{$D_{-1}'$}
$$\includegraphics[width=6cm]{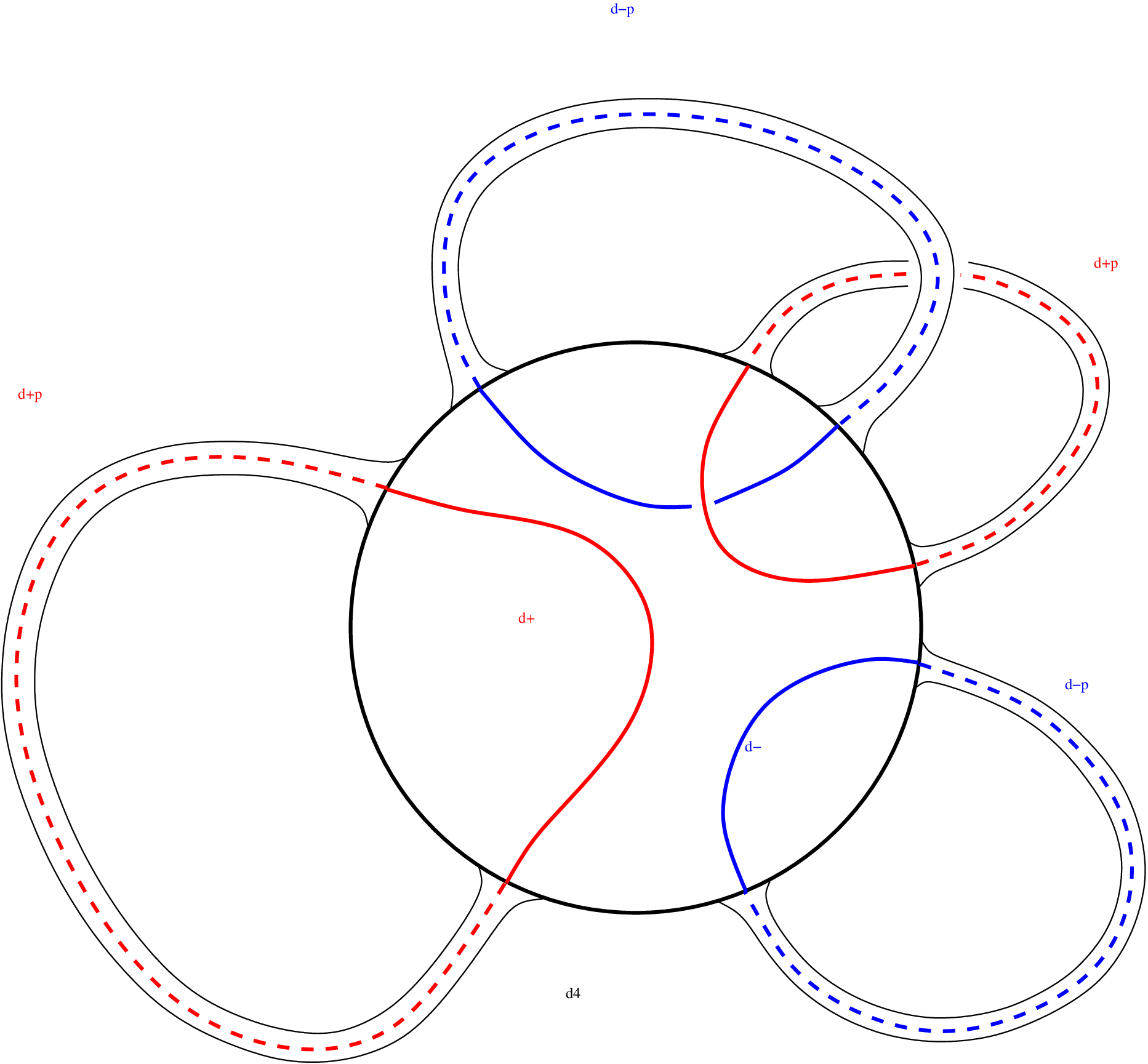}$$ }
\caption{\label{fig-surg} A $1$-surgical embedding.}
\end{figure}

\begin{figure}[H]
$$\includegraphics[width=7cm]{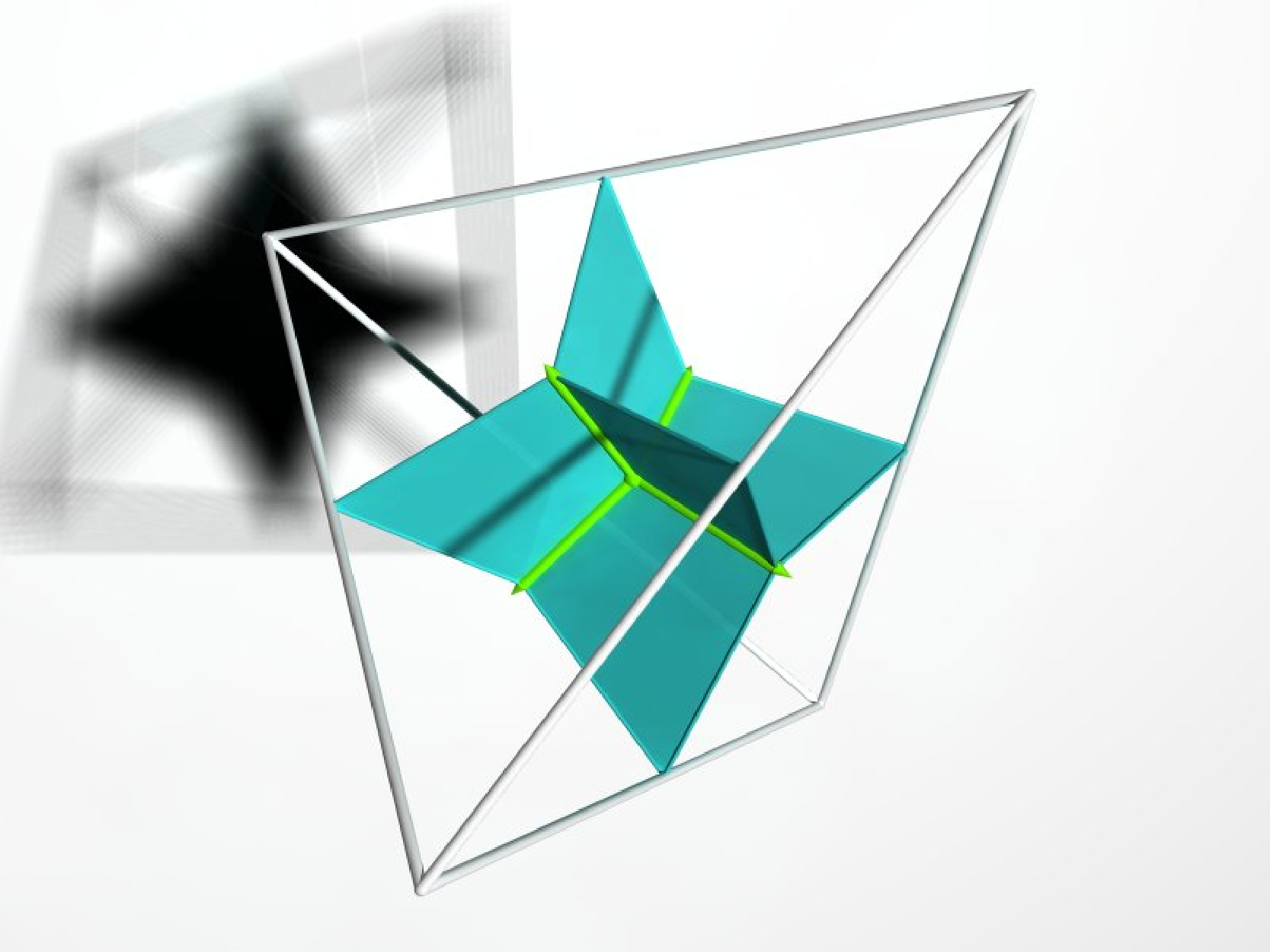}$$
\caption{\label{dualbits}Dual polyhedral bits inside a tetrahedron $\Delta_3$}
\end{figure}

\begin{constr}\label{tenh}({\it 1-Surgical Embedding}): 
If $M$ is a surgery on a smooth slice link such that the surgery coefficients all
belong to the set $\{1,-1\}$, then $M$ admits a smooth embedding into 
a homotopy $4$-sphere. 
\begin{proof}
Write the link $L \subset S^3$ as the union of two disjoint sublinks $L= L_{-1} \cup L_1$ 
where the surgery coefficents for the $L_i$ components are $i$ for $i \in \{-1,1\}$. 
Let $D \subset D^4$ be the slice discs for $L$, $D = D_{-1} \cup D_1$ with $\partial D_i = L_i$ 
for $i = \{-1,1\}$.  Attach $2$-handles to $D^4$ along the components of $L_i$ with framing numbers 
$i$ appropriately for $i \in \{-1,1\}$. Let $D_i'$ be the cores of the attaching handles,
thus $D_i' \cup D_i$ is a union of disjointly embedded $2$-spheres in $N$ whose normal bundles
have Euler number $i$ for $i \in \{-1,1\}$. Recall that $\CProj^2$ has this decomposition: it is a 
$D^2$-bundle over $S^2$ ($\CProj^1$) with Euler number $1$, capped-off with a $4$-handle.   
Thus we can replace a tubular neighbourhood of $D_i' \cup D_i$ with a union of $4$-handles
for $i \in \{-1,1\}$, giving a manifold $N'$ with $N = N' \# \CProj^2 \# \cdots \# \CProj^2 
\# -\CProj^2 \# \cdots \# -\CProj^2$.  This is more commonly known as a `blow-down' operation. 
Thus, $N'$ is contractible, and $\partial N' = M$ so
the double of $N'$ is a homotopy $S^4$ containing $M$. See Figure \ref{fig-surg}. 
\end{proof}
\end{constr}

Given $L \subset S^3$ a slice link with slice discs $D \subset D^4$, we say the slice
discs $D$ are in {\it ribbon position} if the function $d : D^4 \to \Real$ given by
$d(v) = |v|^2$ when restricted to $D$ ($f_{|D} : D \to \Real$) is a Morse function 
having no local maxima.  If a link $L$ has slice discs $D$ that can be put into ribbon position,
$L$ is called a {\it ribbon link}.  Whether or not every slice knot is ribbon is a long-standing
open problem in knot theory, due to Ralph Fox, and is called the {\it slice-ribbon problem}. 

\begin{prop}\label{slice-ribbon-cor}
$L$ is a ribbon link, then the manifold $N'$ in the proof of Construction \ref{tenh} 
admits a handle decomposition with a single $0$-handle, followed by only 
$1$-handle attachments and $2$-handle attachments, i.e. there are no $3$ 
or $4$-handles.
\begin{proof}
Let $A$ be the complement of an open tubular neighbourhood of $D$ in $D^4$. The 
distance function $d$ restricts to a Morse function (in the stratified sense)
$d_{|A} : A \to \Real$ with one local minima, and critical points of index 
$(+1,-2)$ on $\partial A$ corresponding to the critical points of index $(+1,-1)$
of $f_{|D} : D\to \Real$, and critical points of index $(+2,-1)$ corresponding
to critical points of index $(+2,0)$ for $f_{|D} : D \to \Real$. So $A$ consists
of a $4$-ball with $1$-handles and $2$-handles attached. Let $B$ be $N'$ with an 
open tubular neighbourhood of the spheres $\{D_i \cup D_i' : \forall i\}$ removed.
$B$ is $A$ with generalized handles (in the sense of Bott \cite{Bott}) attached.  The
generalized handles correspond to the spheres $D_i\cup D_i'$ for each $i$, and are
trivial $I$-bundles over their core $S^1 \times D^2$. For each $i$ one can
think of this generalized handle as a $2$-handle followed by a $3$-handle attachment.
We construct $N$ by attaching $4$-handles to $B$, one for each $i$.  
The $4$-handles cancel the above
$3$-handle attachments since they satisfy the conditions of Smale's Handle Cancellation
Lemma (see for example \cite{Kosinski} VI.7.4) -- i.e. the attaching sphere of the 
$4$-handle intersects the belt sphere of the $3$-handle transversely in a single point
(the belt sphere consists of two points one in $M$ and one not in $M$). 
Thus $N'$ has a handle decomposition with one $0$-handle, and only $1$ and $2$-handles
attached. 
\end{proof}
\end{prop}

Since $N'$ is contractible, the presentation of $\pi_1 N'$ coming from Proposition 
 \ref{slice-ribbon-cor} must be a presentation of the trivial group, moreover the number
of generators and relators is equal, this is called a ``balanced presentation.''  If the Andrews-Curtis
Conjecture were true \cite{GS} we could cancel the $1$ and $2$-handles of $N'$ using handle
slides, so $N'$ would be diffeomorphic to the standard $D^4$ and $M$ would embed in $S^4$. 
The upshot of this observation is that if we use ribbon links in Proposition \ref{tenh}
and verify the presentation of $\pi_1 N'$ can be trivialized by Andrews-Curtis moves, then
we have verified that the manifold $M$ embeds in $S^4$.  The presentation of $\pi_1 N'$
has the form 
$$\pi_1 N' = \langle g_1, \cdots, g_k : r_1, \cdots, r_j, R_1, \cdots, R_l \rangle$$
where the generators $g_i$ correspond to the local minima of $d$ on the slice discs, 
the relators $r_i$ correspond to the saddle points of $d$ on the slice discs, and
the relators $R_i$ correspond to the framing curves of the link $L$ -- so $k = j + l$. These
presentations are readily computed from a ribbon diagram for $L$.

Constructions \ref{tenh} and \ref{surgical} have relatively simple
implementations. For example, given a hyperbolic manifold which satisfies
Theorem \ref{Hant}, using SnapPea one can drill out selections of curves from $M$ 
then look for the resulting manifold in previously-enumerated tables of hyperbolic link complements. Frequently this technique finds useful surgery presentations. See the beginning of \S \ref{hyperbolicity} for 
details.

Notice that it is relatively easy to construct embeddings of many $3$-manifolds
in homology spheres, for example: A homology $3$-sphere embeds in a homology $4$-sphere if and
only if it is the boundary of a homology $4$-ball. The boundary of any homology $4$-ball
is a homology $3$-sphere, thus constructing embeddings of homology $3$-spheres in homology
$4$-spheres is essentially the same problem as constructing homology $4$-balls.   If 
$M$ is a homology $3$-sphere then $M \# (-M)$ embeds in a homology $4$-sphere -- simply drill 
out a tubular neighbourhood of $\{*\} \times I$ from $M \times I$ to construct a homology
$4$-ball bounding $M \# (-M)$.  If $B$ is an open $3$-ball in a homology $3$-sphere $M$, the manifold 
$(M \setminus B) \times S^1 \cup S^2 \times D^2$ is another homology $4$-sphere
containing $M \# (-M)$. 

If $M$ has non-trivial homology the situation is a little more subtle.  Consider when a $4$-manifold 
$W = V_1 \cup_M V_2$ is a homology $4$-sphere.  By a simple Mayer-Vietoris argument, this happens if and 
only if the manifolds $V_1$ and $V_2$ are orientable and the maps $H_1 M \to H_1 V_1 \oplus H_1 V_2$
and $H_2 M \to H_2 V_1 \oplus H_2 V_2$ are isomorphisms. By considering the long exact sequences
of the pairs $(V_i,M)$ for $i \in \{1,2\}$ and a Poincar\'e Duality argument, this is equivalent to
the statement that the horizontal maps in the commutative diagrams below ($\forall \ \{i,j\} = \{1,2\})$ 
are isomorphisms. 

$$\xymatrix{ H_2(V_j,M,\Zed) \ar[dr]^-{\partial} \ar[rr]^-{\simeq} & & H_1 (V_i,\Zed) & \\
		& H_1 (M,\Zed) \ar[ur]^-{i_*} &}$$

Thus, the problem of constructing an embedding of an arbitrary $3$-manifold into a homology $S^4$
can be thought of as a type of `simultaneous cobordism' problem. 

\begin{constr}\label{homosphemb}
Let $M$ be the result of a surgery along a link $L = L_1 \sqcup L_2$.  Assume that
$L_1$ is smooth slice and that the surgery coefficients for $L_1$ are all zero.  Further, 
assume that the matrix of linking numbers $\linkno_{i,j}$ where $i$ indexes the component of $L_1$
and $j$ indexes the components of $L_2$ is square and invertible, then $M$ is the
boundary of a homology $4$-ball.  If we weaken this last condition to  the matrix $\linkno_{i,j}$
is square with non-zero determinant, then $M$ is a rational homology sphere bounding a
rational homology ball. 
\end{constr}

In the case of manifolds that fibre over $S^1$ there is a spinning construction that produces
many embeddings.

\begin{constr}\label{fibre-const} Let $M$ be a closed orientable $3$-manifold which fibres over $S^1$.
Let $W$ be the fibre of the locally trivial fibre bundle $W \to M \to S^1$
and let $f : W \to W$ be the monodromy, i.e. $M = \Real \times_\Zed W$ where
$\Zed$ acts on $\Real$ by translation, and the action on $W$ is generated
by $f$. If $W$ admits an embedding into $S^3$ such that $f$ extends to 
an orientation-preserving diffeomorphism of $S^3$, then $M$ embeds smoothly
in $S^4$.  
\begin{proof}
The diffeomorphism $f : (S^3,W) \to (S^3,W)$ is isotopic to the identity
when considered as a diffeomorphism of $S^3$ \cite{Cerf}.  Let $F : [0,1] \times S^3 \to S^3$
be such an isotopy: $F(0,x) = x$ and $F(1,x)=f(x)$ for all $x \in S^3$.  
Let $B$ be an open $3$-ball which is disjoint from $W$ and fixed pointwise by $F$. 
Let $B'$ be the closure of the complement of $B$ in $S^3$, thus $f$ can be assumed to
be of the form $f : (B',W) \to (B',W)$, and $f$ restricts to the identity on
$\partial B'$.  Let $\hat F : I \times B' \to B'$ be the corresponding isotopy. 
Consider $S^4$ to be the union $S^4 = (D^3 \times S^1) \cup (S^2 \times D^2)$, 
where we identify $B'$ with $D^3$, then 
$\{ (\hat F(x,t), e^{2\pi i t}) : x \in W, t \in [0,1] \} \subset D^3 \times S^1$
is the embedding of $M$ in $S^4$. 
\end{proof}
\end{constr}

It seems appropriate to call such embeddings `deform-spun' due to the analogy with
Litherland's spinning construction for knots \cite{Litherland}.
It has been known since the work of Crisp and Hillman \cite{CrispH} that not all manifolds that fibre over
$S^1$ which embed in $S^4$ admit deform-spun embeddings. At present the only examples 
of this type that are known are $0$-surgeries on fibred smooth slice knots (see \S \ref{embeddable_man} item 
\ref{nondseg} for an example). 

Embeddings for some special families of $3$-manifolds have been worked out in the literature. 
A class that has received particular attention are the Seifert-fibred homology spheres.

\begin{thm}(Casson, Harer \cite{CH})\label{CHthm} The Brieskorn homology spheres 
$\Sigma(p,q,r)$ smoothly embed in $S^4$ provided $(p,q,r)$ is of the type:
\begin{enumerate}
\item $(p,pa +1, pa +2)$ or $(p,pa-2,pa-1)$ for $p$ odd.
\item $(p,pa-1,pa+1)$ for $p$ even and $a$ odd.
\item $(2,3,13)$ or $(2,5,7)$
\item $(2,5,9)$ or $(3,4,7)$
\end{enumerate}
\begin{proof}
Casson and Harer prove that these Brieskorn spheres $\Sigma$ bound contractible $4$-manifolds $M$
where $M$ has a handle decomposition with a single $0$, $1$ and $2$-handle, and no $3$ or $4$-handles. 
Thus the corresponding handle decomposition for $M \times I$ can be trivialized via handle-slides, 
making $M$ a smooth submanifold of $\partial (M\times I) \simeq S^4$.
\end{proof}
\end{thm}

The statement of Theorem \ref{CHthm} uses the numbering convention of \cite{CH} 
together with the observation that Casson and Harer's families (3) and (4) are finite. Other
useful related references are \cite{AkKi}, \cite{FinSter}.

\begin{thm}(Stern) \cite{Stern}
The Brieskorn spheres $\Sigma(p,q,r)$ bound contractible $4$-manifolds 
provided $(p,q,r)$ is of the form below. Thus, these Brieskorn homology spheres
embed in homotopy $4$-spheres.
\begin{itemize}
\item $(p, pa \pm 1, 2p(pa \pm 1) + pa \mp 1)$ for $p$ even and $a$ odd.
\item $(p, pa \pm 1, 2p(pa \pm 1) + pa \pm 2)$ for $p$ odd 
\item $(p, pa \pm 2, 2p(pa \pm 2) + pa \pm 1)$ for $p$ odd
\end{itemize}
\end{thm}

Stern's contractible $4$-manifolds are constructed from a $4$-ball by attaching
two $1$-handles and then two $2$-handles.

There is one further construction of embeddings of $3$-manifolds in $S^4$ due to Zeeman
and Litherland.  Let $K$ be a ``long knot'' i.e. an embedding $K : D^1 \to D^3$ which 
agrees with the standard inclusion $t \longmapsto (t,0,0)$ on $\{\pm 1\} = \partial D^1$.  
Let $f$ be a diffeomorphism of $D^3$ which fixes pointwise $\partial D^2$ and $\image(K)$.  
By Cerf's Theorem \cite{Cerf}, there is a smooth $1$-parameter family 
$F : D^3 \times [0,1] \to D^3$ such that $F(x,t) = x$ for all $t \in [0,1]$ and $x \in \partial D^3$, with 
$F(x,0) = x$ for all $x \in D^3$ and $F(x,1) = f(x)$ for all $x \in D^3$. 
$F(K(x),t)$ is an isotopy which starts and ends at $K$. Conversely, by the 
Isotopy Extension Theorem, an isotopy that returns $K$ to itself gives a diffeomorphism of the 
pair $(D^3,K)$. These two processes are mutually inverse in the sense that there is an 
isomorphism of the fundamental group of the `space of maps' of type $K$, and
the mapping class group of the pair $(D^3,K)$ (see for example \cite{Budney} 
for details).  Consider $S^4$ to be the union $(D^3 \times S^1) \cup (S^2 \times D^2)$,
then the deform spun knot corresponding to $f$ is the embedding
$$S^2 \equiv (D^1 \times S^1) \cup (S^0 \times D^2) \to (D^3 \times S^1) \cup (S^2 \times D^2) \equiv S^4$$
given by 
$$D^1 \times S^1 \ni (x, e^{2\pi i \theta} ) \longmapsto (F(K(x), \theta), e^{2 \pi i \theta}) \in D^3 \times S^1$$
$$S^0 \times D^2 \ni (a,b) \longmapsto ((a,0), b) \in S^2 \times D^2$$

\begin{thm}\label{DSthm}\cite{Litherland}
Let $M : (D^3,K) \to (D^3,K)$ denote the diffeomorphism induced from rotating $K$ by
$2\pi$ around the axis $[-1,1]\times \{0\}^2 \subset D^3$, a `meridional Dehn twist'. 
If $f : (D^3,K) \to (D^3,K)$ preserves a Seifert surface for $K$, then the complement of
the deform-spun knot associated to $M^n\circ f$ fibres over $S^1$, provided $n \neq 0$.  
\end{thm}

Zeeman proved Theorem \ref{DSthm} in the case that $f$ was the identity automorphism of
$D^3$.  He also went on to show that the fibre is the $n$-fold cyclic branch cover of $D^3$ 
branched over $K$.  So for example, if $n = \pm 1$ and
$f = \Id$, the associated deform-spun knot is trivial, as it bounds a disc. 
Litherland identified the fibre in the more general case. Let $\Sigma$ be
the preserved Seifert surface. This means that $\Sigma$ is an oriented surface in $D^3$ whose 
boundary consists of $K$ union a smooth arc in $\partial D^3$ connecting the endpoints of $K$ 
and that $f(\Sigma)=\Sigma$.  Let $C_K$ denote $D^3$ remove an open tubular neighbourhood of 
$K$, and let $X$ denote $C_K$ remove an open tubular neighbourhood of $C_K \cap \Sigma$.
Denote the two components of the boundary of the tubular neighbournood of $C_K \cap \Sigma$ in $C_K$
by $\Sigma_1$ and $\Sigma_2$  respectively   (thought of as the boundary of $\Sigma \times [1,2]$).
Litherland shows that the Seifert surface for the deform-spun knot is diffeomorphic to the
space $X \times \{1,2,\cdots,n\} / \sim$ where the equivalence relation
is defined by $( (s,1),i ) \sim ((s,2), i+1)$ for $i \in \{1,2,\cdots,n-1\}$ and
$( (s,2), n ) \sim ( (f(s),1), 1)$, where $(s,i) \in \Sigma_i$.  
If one goes on to write $f_{|\Sigma}$ as a product of Dehn twists, this allows the further description
of the Seifert surface as a surgery on a link in a cyclic branch cover of $(D^3,K)$. 

Theorem \ref{DSthm} gives us a rich source of $3$-manifold embeddings in $S^4$, for example, 
the lens spaces $L_{p,q}$ for $p$ odd are $2$-sheeted branched cover over $S^3$
with branch point set the corresponding $2$-bridge knot, thus punctured lens spaces with odd 
order fundamental group embed in $S^4$. Thus the connect sum $L_{p,q} \# -L_{p,q}$ embed smoothly in $S^4$. 
Similarly, a punctured Poincar\'e Dodecahedral Space embeds in $S^4$ by using the $5$-fold branch
cover of $(D^3,K)$ where $K$ is the trefoil. 

If $M_1$ and $M_2$ are lens spaces such that $M_1 \# M_2$ embeds in $S^4$, it follows from Theorem \ref{Hant} that $\pi_1 M_1 \simeq \pi_1 M_2$, and from the torsion linking form that the order of $\pi_1 M_i$ must be odd \cite{KK}. Historically the first proof of this is due to Epstein \cite{Ep}, who used different techniques. This led to one of the more interesting conjectures about $3$-manifolds embedding in $S^4$, due to Gilmer and Livingston \cite{GilLiv} concerning when a connect-sum of two lens spaces embeds in $S^4$.   The Gilmer-Livingston conjecture was solved by Fintushel and Stern \cite{FinSter}, and recently generalized by Andrew Donald \cite{ADon} to the case of an arbitrary connect-sum of lens spaces.

\begin{thm}\label{fsthm}\cite{FinSter, KK, GilLiv, Rub, ADon}
A manifold $M$ that is a connect sum of finitely many lens spaces smoothly embeds in $S^4$ if and only if $M$ is a connect-sum of finitely many manifolds of the form $L_{p,q} \# L_{p,-q}$ where $p$ is odd.  Stated another way, $M$ must be a balanced connect sum of lens spaces and their orientation-reverse, where the lens spaces are required to have fundamental groups of odd order. 
\end{thm}

So for example $L_{p,1}$ admits an orientation-reversing diffeomorphism, but it does not embed in $S^4$ provided $p \geq 2$.  But a connect sum of $k$ copies of $L_{p,1}$ embeds in $S^4$ if and only if $k$ is even and $p\geq 2$ is odd. 
Fintushel and Stern's result is the case of the above theorem where there is precisely two prime summands.  Donald's result is a generalization of \cite{Lisca} where Lisca determines when an arbitrary connect-sum of lens spaces bounds a rational homology ball. 

Donald makes use of a mixed branch cover / slice disc embedding construction. 

\begin{defn} A link $L \subset S^3$ is {\it doubly-slice} if there is an unknotted $2$-sphere $M \subset S^4$ such that $M$ intersects $S^3 \times \{0\}$ transversely, and $L \times \{0\} = M \cap (S^3 \times \{0\})$.  
\end{defn}

\begin{constr}\label{doubleslice}\cite{ADon} If a $3$-manifold $M$ is a finite cyclic branched cover of $(S^3,L)$ with $L$ doubly-slice, then $M$ embeds smoothly in $S^4$. 
\end{constr}

The proof amounts to observing that any finite cyclic branched cover of $S^4$ branched over an unknot is diffeomorphic to $S^4$.  Construction \ref{doubleslice} would not be useful if there wasn't a large class of doubly-slice links.  Donald does so in his Proposition 2.6, constructing a link $L_{a,n} \subset S^3$ such that the double branch cover of $(S^3,L_{a,n})$ is a Seifert fibre space of type $SFS[S^2 : \frac{1}{a}, \frac{1}{a}, -\frac{1}{a}, -\frac{n}{na+1} ]$. 


One can algorithmically construct all $3$-manifolds that embed smoothly in $S^4$.  The algorithm goes like
this: Start with any triangulation of $S^4$.  Enumerate the vertex-normal $3$-manifolds
in that triangulation.  In particular, find all vertex-normal solutions to the gluing
equations, and triangulate them.  Barycentrically subdivide the triangulation of $S^4$ and repeat.  
All $3$-manifolds that embed in $S^4$ eventually appear as vertex-normal solutions in 
{\it any} sufficiently-fine triangulation of $S^4$. This is a consequence of Whitehead's 
proof that smooth manifolds admit triangulations.   This procedure is implemented in Regina (as of version 5.0) and 
was used to construct several embedding examples in Section \ref{emb_htpy}.  This technique also recovers most of
the embeddings in Section \ref{embeddable_man}.  The downside to this technique is it's computationally extremely expensive.  The upside is it finds embeddings of manifolds that have not been found via any other technique. 

There are several obstructions to embedding rational homology spheres in $S^4$ which
utilize $\spin$-structures and $\spinc$-structures.  We summarise the useful properties
of these invariants, but first a quick review of orientation, $\spin$ and $\spinc$ structures on manifolds.
Helpful references for this material are \cite{KM, MS, Mil, Larry}.

The group $\spin(n)$ is the connected 2-sheeted cover of the Lie group $SO_n$, together with the
Lie group structure making $\spin(n) \to SO_n$ a homomorphism of Lie groups. Provided $n \geq 3$, $\spin(n)$ is the universal cover of $SO_n$.  The group $\spinc(n)$ is the twisted-product $\spin(n) \times_{\Zed_2} \spin(2)$ where $\Zed_2$ acts diagonally as the covering transformation of both factors. Thus, there are Lie group submersions:

$$ \Zed_2 \to \spin(n) \to SO_n \hskip 15mm \spin(2) \to \spinc(n) \to SO_n$$

Notice that there is a canonical isomorphism of Lie groups $U_2 \simeq \spinc(3)$, 
since $SU_2 \subset U_2$ is naturally isomorphic to $\spin(3) = S^3$, and the
diagonal matrices in $U_2$ are naturally isomorphic to $\spin(2)$, moreover, $SU_2$ intersects
the diagonal matrices at precisely $\pm 1$. More generally, $U_n \simeq SU_n \times_{\Zed_n} U_1$.

Given an $n$-manifold $N$ let $TN$ denote the tangent bundle of $N$, this the union
of all the tangent spaces to $N$. $TN$ is a vector bundle over $N$. 
The space of all bases to the tangent spaces of $N$ is called the principal $GL_n$-bundle 
associated to $N$, and will be denoted $GL_n(TN)$. $GL_n(TN)$ is a fibre bundle over $N$
with fibre the Lie group $GL_n$, thus there are fibrations:

$$ GL_n \to GL_n(TN) \to N \hskip 15mm GL_n(TN) \to N \to BGL_n$$

The map $N \to BGL_n$ is called the classifying map for the bundles
$TN \to N$ and $GL_n(TN) \to N$ respectively.  Since the inclusion 
$O_n \to GL_n$ is a homotopy-equivalence, a choice of a Riemannian metric on $N$
allows us to replace $GL_n$ by $O_n$ in the discussion above. 

An orientation of $N$ is a homotopy class of lifts of the classifying
map $N \to BO_n$ to $BSO_n$. For an oriented manifold $N$, a $\spinc(n)$-structure on 
$N$ is a homotopy class of lifts of maps $N \to BSO_n$ to maps $N \to B\spinc(n)$.  Similarly, a
$\spin(n)$-structure is a homotopy class of lifts of $N \to BSO_n$ to
$N \to B\spin(n)$.  Essentially by definition, two
$\spin(n)$-structures on $N$ differ by an element of $[N,B\Zed_2] \equiv H^1(N,\Zed_2)$.
Similarly, two $\spinc(n)$-structures on $N$ differ by an element of
$[N,B\spin(2)] \equiv H^2(N,\Zed)$.

Every orientable $3$-manifold has a trivial tangent bundle \cite{Kirby}, so it has both a 
$\spin(3)$ and a $\spinc(3)$-structure.  In general, a manifold $N$ has a $\spin(n)$ structure if
and only if it is orientable and the 2nd Stiefel-Whitney class is zero, $w_2(N)=0$.  
Equivalently, if its tangent bundle trivializes
over the $2$-skeleton of $N$ -- moreover, the $\spin(n)$-structure is taken to be a homotopy class
of such a trivialization, once restricted to the $1$-skeleton. 
 $N$ has a $\spinc(n)$-structure if and only if $w_2(N)$ is
the reduction of an integral cohomology class.  Equivalently, this is if and only if
a direct sum with a complex line bundle admits a $\spin$-structure.  Another equivalent
definition is that (if $N$ has odd dimension, stabilize by adding a trivial $1$-dimensional
vector bundle) a $\spinc$-structure is a homotopy class of almost complex structures over the 
$2$-skeleton such that a representative almost complex structure extends over the $3$-skeleton. 

\begin{thm}\label{spintests}\cite{Kirby, OzSz} If $M$ is a $\spin$ $3$-manifold there exists an invariant,
called the Rochlin invariant, taking values in $\Rat / 2\Zed$. The Rochlin invariant of $M$
is $\mu(M) = \frac{\signature(W)}{8} \in \Rat/2\Zed$ where $W$ is a $\spin$-manifold such that $\partial W = M$. 
$\signature(W)$ is the signature of the intersection form on $fH_2(W,\Zed) = H_2(W,\Zed)/\tau H_2(W,\Zed)$.
When $M$ is a homology sphere $\signature(W)$ is divisible by $8$, so $\mu(M) \in \Zed_2$.

The Rochlin invariant has an integral lift for homology spheres, called the 
$\overline{\mu}$-invariant \cite{Siebenmann}. 
$\overline{\mu}$ is a homology cobordism invariant for 
Seifert fibred homology spheres (see \cite{Sav} Corollary 7.34).

If $M$ is a rational homology $3$-sphere with a $\spinc$-structure, there is an invariant called 
the Ozsv\'ath-Szab\'o $d$-invariant or `correction term,' taking values in $\Rat$. It is a rational
homology $\spinc$-cobordism invariant and additive under connect-sum. 
\end{thm}

The above theorems explain why we're interested in $\spin$ and $\spinc$ structures -- the extra structure
given to the tangent bundle allows for more delicate constructions.  For our purposes, a $\spin$ structure is the
most sensitive tangent bundle structure we'll ever need. This is because a connected $4$-manifold 
which bounds a non-empty $3$-manifold has a trivial tangent bundle if and only if it has admits a $\spin$
structure -- to see this, notice such $4$-manifolds have the homotopy-type of a $3$-complex.  The
tangent bundle of a $4$-manifold with a $\spin$-structure trivializes over the $2$-skeleton, and the
obstruction to extending over the $3$-skeleton (and thus the entire manifold) lives in a $3$-dimensional
twisted cohomology group with coefficients $\pi_2 SO_4 = \pi_2 \spin(4) = \pi_2 (S^3 \times S^3) = 0$.  

\begin{defn}
Given a rational homology sphere $M$, let $\vec \mu (M)$ be the function whose domain is the $\spin$-structures on 
$M$ and whose values are the Rochlin invariants of $M$ with the associated $\spin$-structure.  Similarly, 
let $\vec d(M)$ be the function whose domain is the $\spinc$ structures on $M$ and whose values are 
the associated $d$-invariants. 
\end{defn}

\begin{cor}\label{vectest} ($\vec \mu$ and $\vec d$ tests)
Given a rational homology sphere $M$ which admits a smooth embedding into
a homology $4$-sphere, $|H_1(M,\Zed)|=k^2$ for some $k$.  Moreover, there are $2k-1$ zeros in
$\vec d(M)$.  Similarly, $|H_1(M,\Zed_2)| = l^2$ for some $l$, and there are 
$2l-1$ zeros in $\vec \mu(M)$.
\begin{proof}
Assume $M$ embeds in $S^4$, then $M$ separates $S^4$ into two rational homology
balls $V_1$ and $V_2$.  Since $V_1 \subset S^4$, $V_1$ has a trivial tangent bundle. 
If we fix a trivialization of $TV_1$, the $\spinc$ structures on $V_1$ correspond
to elements of $[V_1,B\spin(2)] = H^2(V_1,\Zed)$.  

Consider the problem of determining the $\spinc$-structures on $M$ which restrict from 
$\spinc$-structures on $V_1$.  If we use the trivialization of $TM$ coming from considering 
$M = \partial V_1$, this then amounts to determining the image of the restriction map 
$[V_1,B\spin(2)] \to [M,B\spin(2)]$ which by the Brown Representation Theorem is equivalent to 
the image of the map $H^2(V_1,\Zed) \to H^2(M,\Zed)$.  Via Poincar\'e duality this map is 
equivalent to $H_1(V_2,\Zed) \to H_1(M,\Zed)$ whose image is the kernel of the map
$H_1(M,\Zed) \to H_1(V_1,\Zed)$.  In other words, we have the hyperbolic splitting
$H_1(M,\Zed) \simeq H_1(V_1,\Zed) \oplus H_1(V_2,\Zed)$, and the $\spinc$-structures on
$M$ that extend to $V_1$ correspond to the subgroup $H_1(V_2,\Zed)$. Similarly, the $\spinc$-structures
on $M$ which extend to $V_2$ correspond to the subgroup $H_1(V_1,\Zed)$. 

Consider the $\vec \mu (M)$-test. We are considering the image of the map
$[V_1,B\Zed_2] \to [M,B\Zed_2]$, which is equivalent to the map 
$H^1(V_1,\Zed_2) \to H^1(M,\Zed_2)$.  The result is analogous, except here
we use the splitting $H^1(M,\Zed_2) \simeq H^1(V_1,\Zed_2) \oplus H^1(V_2,\Zed_2)$.
\end{proof}
\end{cor}

Corollary \ref{vectest} has a stronger statement, as the zeros in $\vec d$ and $\vec \mu$ have
the shape of an `affine X' in directions specified by the hyperbolic splitting of the
torsion linking form.  

Perhaps the simplest way to compute the Rochlin vector $\vec \mu(M)$ follows this procedure:
\begin{itemize}
\item Find a surgery presentation for $M$. For hyperbolic $3$-manifolds see 
\S \ref{hyperbolicity}. Graph manifolds in essence have canonical surgery presentations
given by their definition, this is also sketched in \S \ref{hyperbolicity}. 
\item Using inverse `slam-dunk' moves (see Figure 5.30 of \cite{GS}), find an
integral surgery presentation for $M$.
\item Enumerate the $\spin$-structures on $M$ via characteristic sublinks (see Proposition 5.7.11 of
\cite{GS}).
\item Use the Kaplan algorithm to find a $\spin$ 4-manifold bounding the $\spin$ 3-manifold specified
by a characteristic sublink (Theorem 5.7.14 of \cite{GS}).
\item From the surgery presentation, the signature is readily computed via basic linear algebra.
\end{itemize}

The reader will notice that the only obstructions to a $3$-manifold embedding in $S^4$
that we have mentioned are obstructions to embedding in homology $4$-spheres.  
Theorems \ref{fsthm} and \ref{chthm} completely describe, for a very limited class of 
$3$-manifolds, precisely which manifolds from that class admit embeddings in $S^4$.  Namely, for
connect-sums of two lens spaces, and for circle bundles over surfaces there is the curious phenomenon that 
these $3$-manifolds embed in $S^4$ if and only if they embed in a homology $4$-sphere.  

Recently, Issa and McCoy \cite{Issa} have made progress applying Donaldson's Theorem to
obstructing embeddings of $3$-manifolds into $S^4$.  Specifically, we are
referring to the theorem that states that compact, oriented, smooth $4$-manifolds have
diagonalizable intersection forms, and when the form is definite the diagonalization can be
performed over the integers. The basic idea of the argument is that if one has a $3$-manifold embedding
$M \to S^4$, one replaces the manifold $V_1 \cup_M V_2 = S^4$ by $X_M \cup_M V_2$, where
$X_M$ is an inspired choice.  Issa and McCoy's techniques work when one can find a $4$-manifold where
$X_M$ has a definite intersection pairing.  They then study the
induced map $H_2 X_M \to H_2 X_M \cup_M V_2$.  Donaldson's theorem characterises the 
geometry of the target, thus if one knows enough about $X_M$ one can obstruct such maps.
They take this argument quite far, using both sides ($V_1$ and $V_2$) of the splitting
to generate obstructions. 

As a warning to the reader, this paper is not exhaustive in its usage of known obstructions
to $3$-manifolds embedding in $S^4$.  Known obstructions to $3$-manifolds embedding in
homology spheres that have not been employed (yet) include: the Casson-Gordon invariants and 
their relatives \cite{FinSter}, and the $w$-invariant \cite{Sav}.  

\section{Manifolds from the census which embed smoothly in $S^4$}\label{embeddable_man}

In the list below, an attempt was made to give all the manifolds 
a more-or-less standard name.  The Seifert-fibred data is all un-normalized. 
This means (among other things) that if you sum up all the fibre-data numbers, you get the
Euler characteristic of the Seifert bundle over the base orbifold, see Orlik for
details \cite{Orlik}. 

\vskip 5mm
\centerline{$\star$ Spherical manifolds $\star$}
\vskip 5mm
\begin{enumerate}
\item $S^3$. $S^3$ is the equator in $S^4$. \ttc
\item $SFS \left[ S^2 : \frac{1}{2}, \frac{1}{2}, -\frac{1}{2}\right] = S^3/Q_8 = 
SFS \left[ \RProj^2 : 2 \right]$. \ttc $H_1 = \Zed_2^2$.
$Q_8$ is the quaternion group of order $8$, i.e. $Q_8 = \{ \pm 1, \pm i, \pm j, \pm k \}$. 
$S^3/Q_8$ appears as the boundary of a tubular neighbourhood of any 
embedding $\RProj^2 \to S^4$ \cite{Han}.  A standard embedding of $\RProj^2$ in 
$\Real^4$ is given by $(x,y,z) \longmapsto (xy,xz,y^2-z^2,2yz)$ where we think of
$S^2 \subset \Real^3$ as the universal cover of $\RProj^2$.
\vskip 5mm
\centerline{$\star$ The $\Real \times S^2$ manifold $\star$}
\vskip 5mm
\item $S^1 \times S^2$. $H_1 = \Zed$. Trivial deform-spun embedding (Construction \ref{fibre-const}), 
also $0$-surgery on unknot (Construction \ref{surgical}). 
\ttc
\vskip 5mm
\centerline{$\star$ Nil manifolds $\star$}
\vskip 5mm

\item \label{nondseg} SFS$\left[ S^1 \times S^1 : 1 \right] = (S^1\times S^1) \rtimes_{\begin{pmatrix}  1 & 1 \\ 0 & 1\end{pmatrix}} S^1$.  $H_1 = \Zed^2$. One obtains this manifold
as a zero surgery on the link $\KR{R}{5^2_1}$ \cite{CrispH}. 
\ttc
\item SFS$\left[ S^1 \ltimes S^1 : 4 \right] = (S^1 \times S^1) \rtimes_{\begin{pmatrix}  -1 & 4 \\ 0 & -1\end{pmatrix}} S^1$. 
$H_1 = \Zed \oplus \Zed_2^2$. This manifold is obtained by zero surgery
on the link $\KR{R}{9^3_{19}}$. Alternatively, it is the unit normal bundle to an embedding
of the Klein bottle in $S^4$ \cite{CrispH}. \ttc 
\item SFS$\left[ S^2 ; \frac{1}{3}, \frac{1}{3}, -\frac{1}{3} \right]$. $H_1 = \Zed_3^2$. This manifold is
obtained as the $0$-surgery on the $(2,6)$-torus link which is a disjoint 
union of two unknots \cite{CrispH} (Construction \ref{surgical}). \ttc 
\item SFS$\left[ \RProj^2 ; \frac{1}{2}, \frac{3}{2} \right]$. $H_1 = \Zed_4^2$. This manifold is obtained 
as zero surgery on the link $\KR{R}{8^2_2}$ \cite{CrispH} (Construction \ref{surgical}). \ttc 
\vskip 5mm
\centerline{$\star$ Euclidean manifolds $\star$}
\vskip 5mm

\item $S^1\times S^1 \times S^1$.  $H_1 =\Zed^3$. Trivial deform-spun embedding (Construction
\ref{fibre-const}), also $0$-surgery on Borromean rings (Construction \ref{surgical}). \ttc
\item $(S^1 \times S^1) \times_{\Zed_2} SO_2$ where $\Zed_2 \subset SO_2$
acts on $S^1 \times S^1$ by $\pi$-rotation on the square torus, so it admits
a deform-spun embedding. This manifold is also 
$SFS \left[ (S^1 \rtimes S^1) : 0 \right]$, so it is the boundary of a tubular
neighbourhood of an embedding of the Klein bottle in $S^4$. $H_1 = \Zed \oplus \Zed_2^2$. 
\vskip 5mm
\centerline{$\star$ Sol manifolds $\star$}
\vskip 5mm

Crisp and Hillman \cite{CrispH} determined the Sol manifolds that embed in $S^4$.  In particular, 
they showed that none of the Sol manifolds which fibre over $S^1$ embed in $S^4$, and of the remaining
Sol manifolds, only three of them embed.  Consider the Klein bottle to be 
$S^1 \times_{\Zed_2} S^1$ where $\Zed_2 = \{\pm 1\}$ acts by
$-1.(z_1,z_2)=(\overline{z_1},-z_2)$. Given a matrix $A = \begin{pmatrix} a & b \\ c & d \end{pmatrix}$
we can describe a $Sol$-manifold as the union of two orientable $I$-bundles over $S^1 \times_{\Zed_2} S^1$.  Precisely, if we consider $S^1 \times S^1$ to be the boundary of this $I$-bundle, the gluing map
$A_* : S^1 \times S^1 \to S^1 \times S^1$ is given by $A_*(z_1,z_2) = (z_1^az_2^b, z_1^cz_2^d)$.
Alternatively, these manifolds can be described as the union of two manifolds of the
form $SFS\left[D^2, \frac{1}{2},\frac{1}{2}\right]$. Identify the boundary with $S^1\times S^1$ where the first
coordinate indicates the fibre direction and the 2nd coordinate the `base' direction, thus such manifolds
are specified by a corresponding gluing matrix $B$, which in the notation of Regina would be
$B= \begin{pmatrix} d-b & b \\ d+c-b-a & b+1 \end{pmatrix}$.

\item SFS $\left[D: \frac{1}{2}, \frac{1}{2}\right]$ U/m SFS $\left[D: \frac{1}{2}, \frac{1}{2}\right]$ 
$m = \begin{pmatrix}  -1 & 3 \\ 0 & 1\end{pmatrix}$ $H_1 = \Zed_4^2$ \ttc
embeds \cite{CrispH} Crisp-Hillman notation: $\begin{pmatrix} 2 & 3 \\ 1 & 2\end{pmatrix}$. 
$0$-surgery on link $\KR{R}{9^2_{53}}$ (Construction \ref{surgical}).

\item SFS $\left[D: \frac{1}{2}, \frac{1}{2}\right]$ U/m 
SFS $\left[D: \frac{1}{2}, \frac{1}{2}\right]$ 
$m = \begin{pmatrix}  -3 & 5 \\ -2 & 3\end{pmatrix}$ $H_1 = \Zed_4^2$ \ttc
embeds \cite{CrispH} Crisp-Hillman notation: $\begin{pmatrix} 2 & -5 \\ 1 & -2\end{pmatrix}$.
$0$-surgery on $\KR{R}{9^2_{61}}$ (Construction \ref{surgical}).

\item SFS $\left[D: \frac{1}{2}, \frac{1}{2}\right]$ U/m 
SFS $\left[D: \frac{1}{2}, \frac{1}{2}\right]$ 
$m = \begin{pmatrix}  -7 & 9 \\ -4 & 5\end{pmatrix}$  $H_1 = \Zed_4^2$ 
embeds \cite{CrispH} Crisp-Hillman notation: $\begin{pmatrix} 2 & -9 \\ 1 & -4\end{pmatrix}$. 
$0$-surgery on $2$-component link (Construction \ref{surgical})
$$\includegraphics[width=2cm]{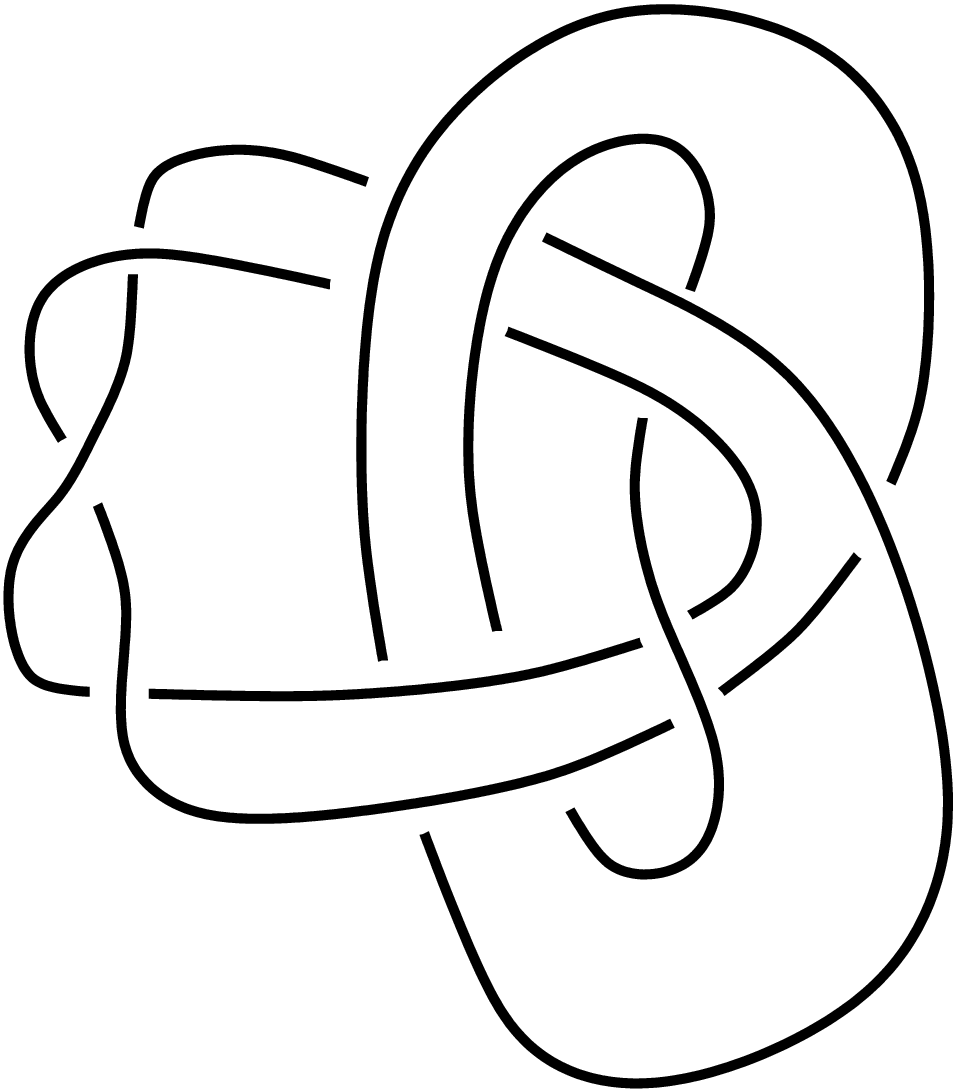}$$

\vskip 5mm
\centerline{$\star$ $SL_2 \Real$ (Brieskorn) homology spheres $\star$}
\vskip 5mm

\item SFS $\left[S^2: \frac{1}{3}, \frac{1}{4}, -\frac{3}{5}\right] = \Sigma(3,4,5)$ case (1) of Theorem \ref{CHthm}. \ttc
\item SFS $\left[S^2: \frac{1}{2}, \frac{1}{5}, -\frac{5}{7}\right] = \Sigma(2,5,7)$ case (2) of Theorem \ref{CHthm}. \ttc
\item SFS $\left[S^2: \frac{1}{3}, \frac{2}{7}, -\frac{5}{8}\right] = \Sigma(3,7,8)$ case (1) of Theorem \ref{CHthm}.
\item SFS $\left[S^2: \frac{1}{2}, \frac{2}{9}, -\frac{8}{11}\right] = \Sigma(2,9,11)$ case (2) of Theorem \ref{CHthm}.
\item SFS $\left[S^2: \frac{1}{2}, \frac{1}{3}, -\frac{11}{13}\right] = \Sigma(2,3,13)$ \ttc
case (3) of Theorem \ref{CHthm}. 
\vskip 5mm
\centerline{$\star$ $SL_2\Real$ rational homology spheres $\star$}
\vskip 5mm

\item SFS $\left[\RProj^2/n2: \frac{1}{3}, \frac{5}{3}\right]$ $H_1 = \Zed_6^2$.\ttc
Proposition 1.2 from Crisp-Hillman \cite{CrispH}.

\item SFS $\left[\RProj^2/n2: \frac{1}{4}, \frac{7}{4}\right]$ $H_1 = \Zed_8^2$. 
Proposition 1.2 from Crisp-Hillman \cite{CrispH}.

\item SFS $\left[\RProj^2/n2: \frac{2}{5}, \frac{8}{5}\right]$ $H_1 = \Zed_{10}^2$. 
Proposition 1.2 from Crisp-Hillman \cite{CrispH}. 

\item SFS $\left[S^2: \frac{1}{2}, \frac{1}{2}, \frac{1}{2}, \frac{1}{2}, -\frac{3}{2} \right]$ $H_1 = \Zed_2^4$. \ttc  
To construct an embedding of this manifold into $S^4$ notice that this manifold is obtained by 
surgery on a regular fibre in the manifold
$$SFS \left[S^2: \frac{1}{2}, -\frac{1}{2}, \frac{1}{2}, -\frac{1}{2}\right]$$ 
which embeds as the unit normal bundle to the `standard' embedding of the Klein bottle in $S^4$ ($W_2=0$).
The surgery curve bounds the disc pictured below -- thus the surgery can be realized as an embedded
surgery.
{\psfrag{C1}[tl][tl][0.8][0]{surgery}
\psfrag{C2}[tl][tl][0.8][0]{disc}
$$\includegraphics[width=6cm]{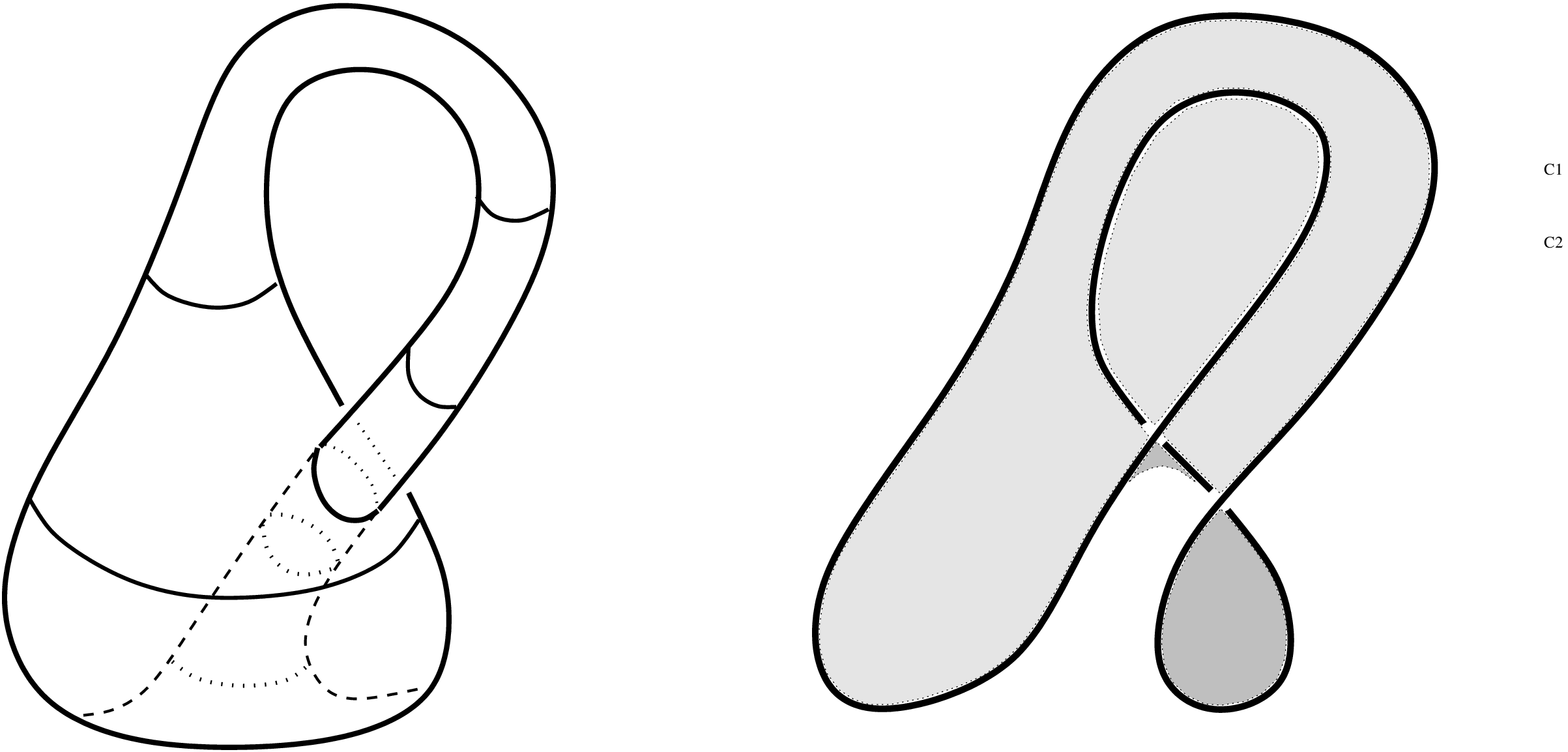}$$
\centerline{Constructing embedding of 
$SFS\left[S^2: \frac{1}{2}, \frac{1}{2}, \frac{1}{2}, \frac{1}{2}, -\frac{3}{2} \right]$}}

\item SFS $\left[S^2: \frac{1}{4}, \frac{1}{4}, -\frac{1}{4}\right]$ $H_1 = \Zed_4^2$. \ttc
$0$-surgery on the $(2,8)$-torus link, see Figure A4 of Crisp and Hillman \cite{CrispH}. 

\item SFS $\left[S^2: \frac{1}{4}, \frac{1}{4}, -\frac{7}{12}\right]$ $H_1 = \Zed_4^2$. 
$\vec d = \begin{pmatrix} 
 0 & 0 & 0 & 0 \\ 
 0 & \frac{1}{2} & -1 & -\frac{1}{2} \\
 0 & -1 & 0 & -1 \\
 0 & -\frac{1}{2} & -1 & \frac{1}{2}
\end{pmatrix}$. 
Characteristic links 

$(\{c\}, \{c,e,f\}, \{d,e\}, \{d,f\})$, 
$\vec \mu = (0,1,0,0)$. 
Surgery diagram $\includegraphics[height=1.5cm]{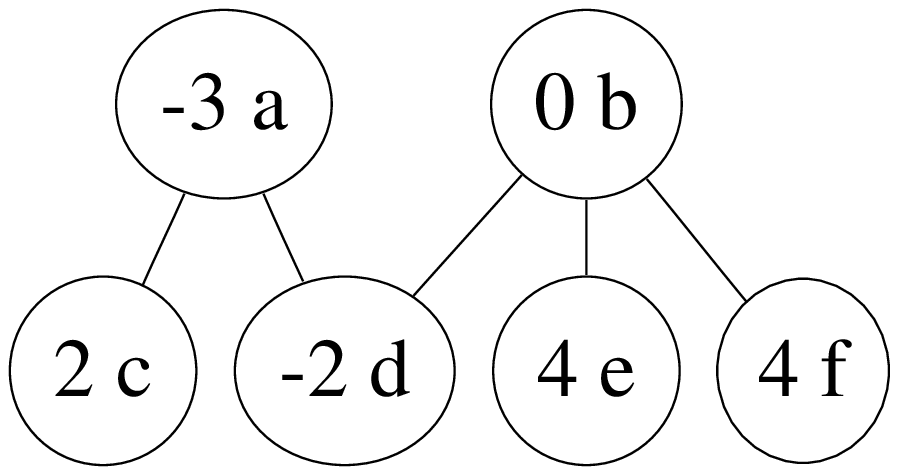}$.
A.Donald constructs an embedding in Example 2.14 \cite{ADon}.

\item SFS $\left[S^2: \frac{1}{2}, \frac{1}{2}, \frac{1}{2}, -\frac{5}{3}\right]$ $H_1 = \Zed_2^2$. \ttc
Characteristic links: $(\{b\}, \{b,c,d\}, \{b,c,e\}, \{b,d,e\} )$, $\vec \mu = (0,0,0, \frac{1}{2})$,
$\vec d = \begin{pmatrix} 
 1 & 0 \\ 0 & 0 \end{pmatrix}$,
surgery diagram: \includegraphics[height=1.5cm]{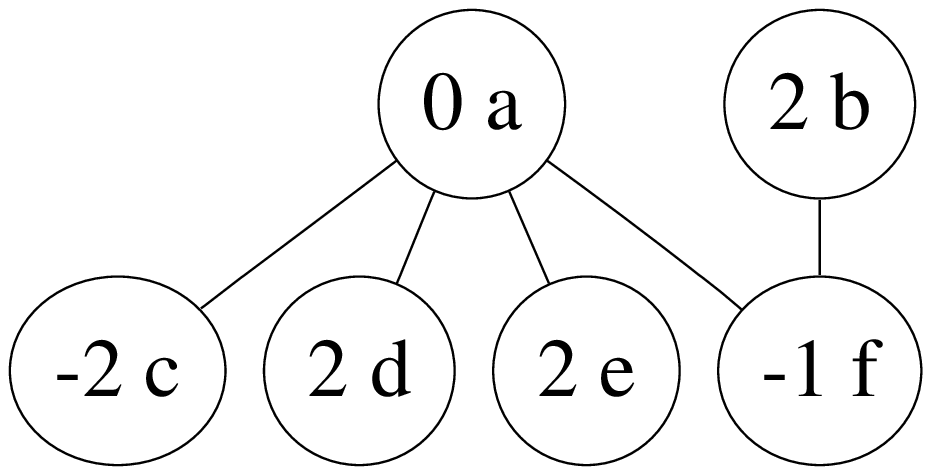}. 
Embeds via Construction \ref{doubleslice}, manifold is the double branch cover of $(S^3,L_{2,2})$.   


\item SFS $\left[S^2: \frac{1}{2}, \frac{1}{2}, \frac{1}{2}, -\frac{8}{5}\right]$ $H_1 = \Zed_2^2$. \ttc
$\vec d = \begin{pmatrix} 
 1 & 0 \\ 
 0 & 0 
\end{pmatrix}$. 
Characteristic links 

$(\{c,f\}, \{d,f\}, \{e,f\}, \{c,d,e,f\})$,
$\vec \mu = (\frac{1}{2}, 0, 0, 0)$.
Surgery diagram $\includegraphics[height=1.5cm]{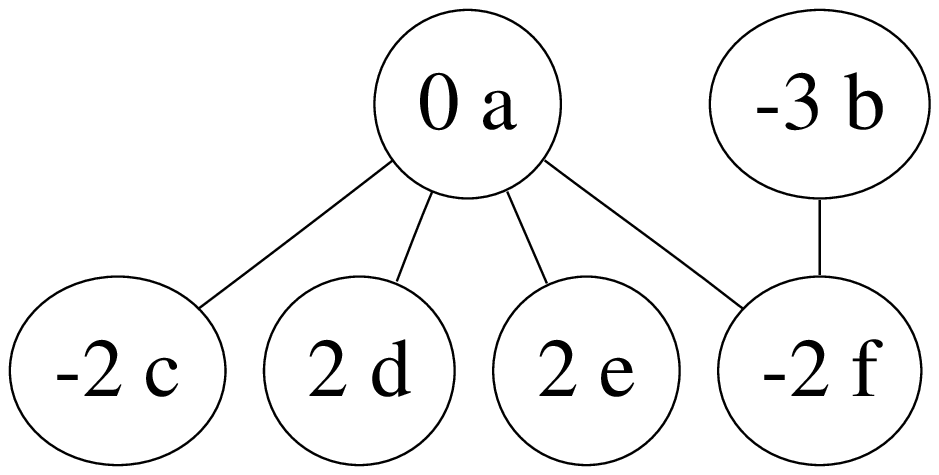}$
Embeds via Construction \ref{doubleslice}, manifold is the double branch cover of $(S^3,L_{2,-3})$.   

\item SFS $\left[S^2: \frac{1}{2}, \frac{1}{2}, \frac{1}{2}, -\frac{11}{7}\right]$ $H_1 = \Zed_2^2$. \ttc
$\vec d = \begin{pmatrix} 
 1 & 0 \\ 
 0 & 0 
\end{pmatrix}$. 
Characteristic links $(\phi, \{c,d\}, \{d,e\}, \{c,e\})$,
$\vec \mu = (0,0,-\frac{1}{2}, 0)$.
Surgery diagram $\includegraphics[height=1.5cm]{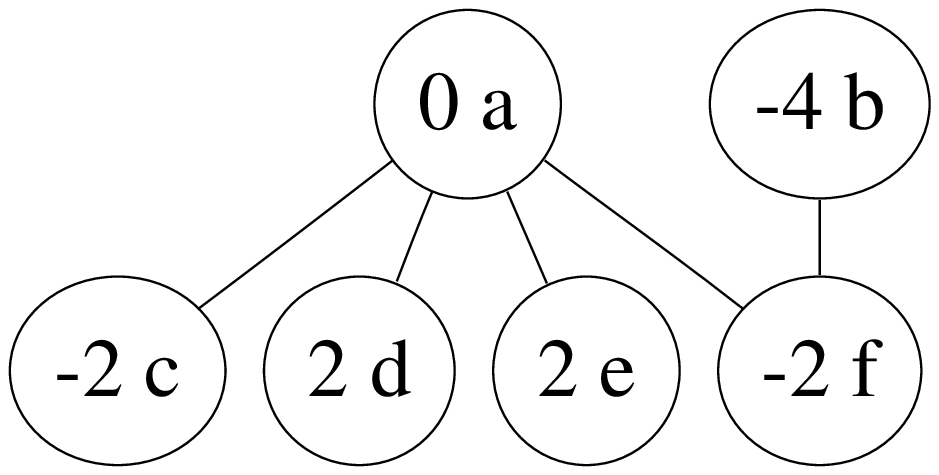}$
Embeds via Construction \ref{doubleslice}, manifold is the double branch cover of $(S^3,L_{2,-4})$.   

\item SFS $\left[S^2: \frac{1}{2}, \frac{1}{2}, \frac{1}{2}, -\frac{14}{9} \right]$ $H_1 = \Zed_2^2$. 
$\vec d = \begin{pmatrix} 
 1 & 0 \\ 
 0 & 0 
\end{pmatrix}$. 
Characteristic links 

$(\{c,f\}, \{d,f\}, \{e,f\}, \{c,d,e,f\})$,
$\vec \mu = (\frac{1}{2}, 0, 0, 0)$.
Surgery diagram $\includegraphics[height=1.5cm]{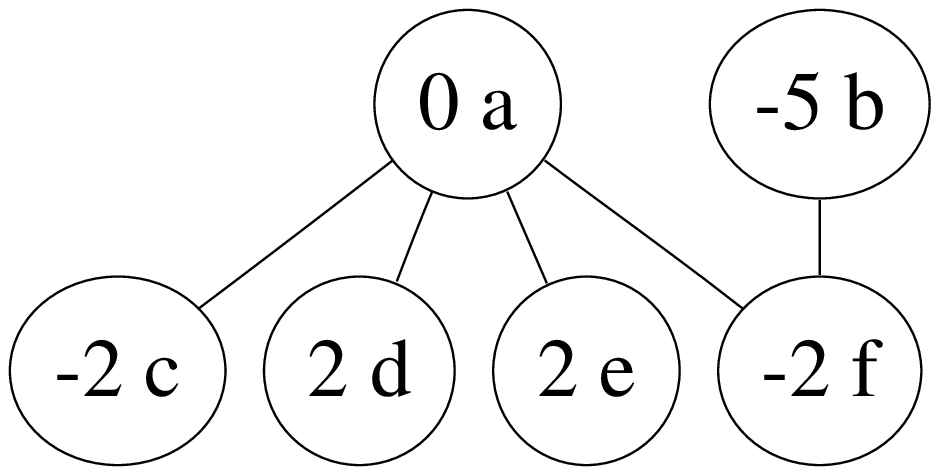}$
Embeds via Construction \ref{doubleslice}, manifold is the double branch cover of $(S^3,L_{2,-5})$.   

\item SFS $\left[S^2: \frac{1}{2}, \frac{1}{2}, \frac{1}{2}, -\frac{17}{11}\right]$ $H_1 = \Zed_2^2$.  
$\vec d = \begin{pmatrix} 
 1 & 0 \\ 
 0 & 0 
\end{pmatrix}$. 
Characteristic links $(\phi, \{c,d\}, \{d,e\}, \{c,e\})$,
$\vec \mu = (0,0,-\frac{1}{2}, 0)$.
Surgery diagram $\includegraphics[height=1.5cm]{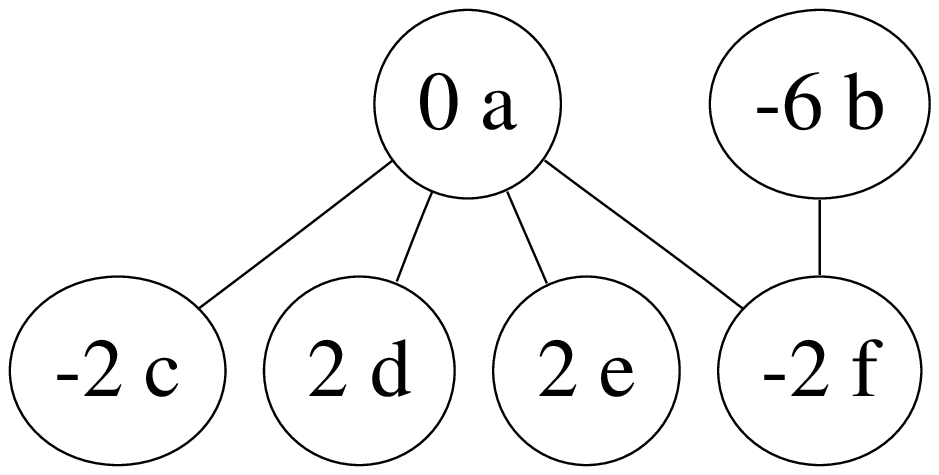}$
Embeds via Construction \ref{doubleslice}, manifold is the double branch cover of $(S^3,L_{2,-6})$.   

\item SFS $\left[S^2: \frac{1}{2}, \frac{1}{3}, \frac{1}{3}, -\frac{4}{3} \right]$ $H_1 = \Zed_3^2$. \ttc
$\vec d = \begin{pmatrix} 
 0 & 0 & 0 \\ 
 0 & \frac{2}{3} & \frac{4}{3} \\
 0 & \frac{4}{3} & \frac{2}{3}
\end{pmatrix}$. $\mu = 0$. 
Surgery diagram $\includegraphics[height=1.5cm]{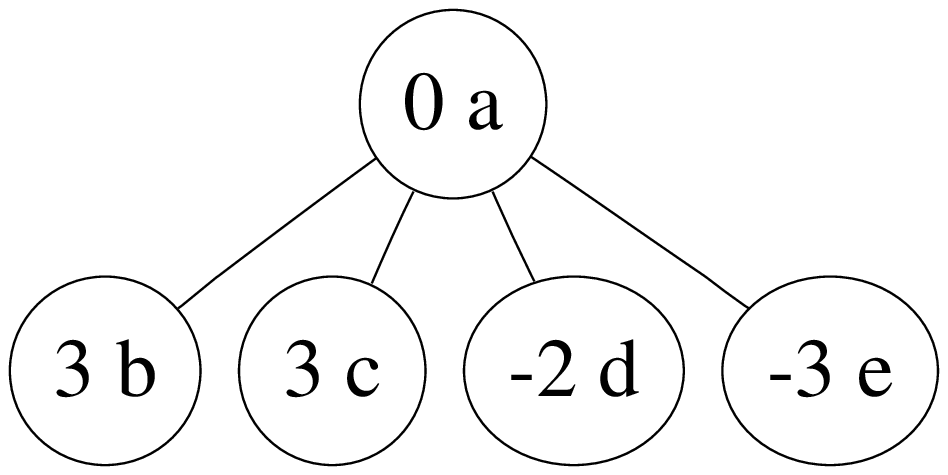}$
Embeds via Construction \ref{doubleslice}, manifold is the double branch cover of $(S^3,L_{3,-1})$.   

\item SFS $\left[S^2: \frac{1}{3}, \frac{1}{3}, \frac{2}{3}, -\frac{5}{4}\right]$ $H_1 = \Zed_3^2$. 
$\vec d = \begin{pmatrix} 
 0 & 0 & 0 \\ 
 0 & \frac{4}{3} & \frac{2}{3} \\
 0 & \frac{2}{3} & \frac{4}{3}
\end{pmatrix}$. $\mu = 0$. 
Surgery diagram $\includegraphics[height=1.5cm]{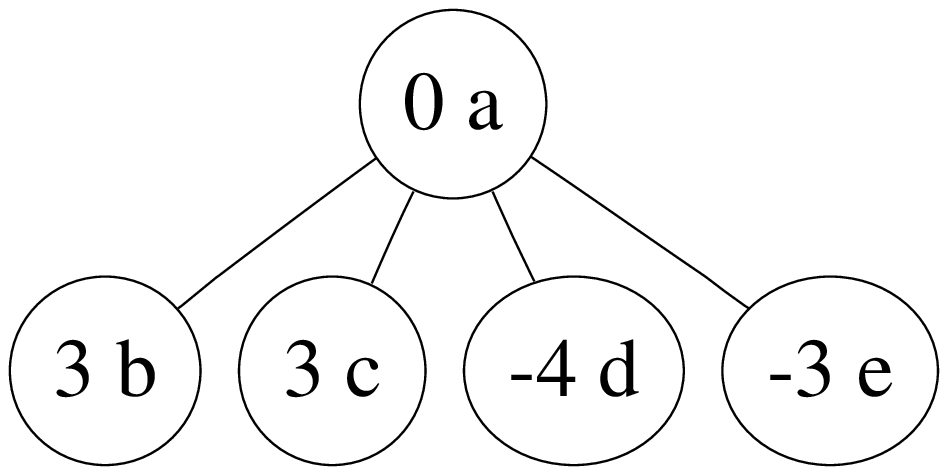}$
Embeds via Construction \ref{doubleslice}, manifold is the double branch cover of $(S^3,L_{3,1})$.

\item SFS $\left[S^2: \frac{1}{3}, \frac{1}{3}, \frac{2}{3}, -\frac{7}{5}\right]$ $H_1 = \Zed_3^2$. 
$\vec d = \begin{pmatrix} 
 0 & 0 & 0 \\ 
 0 & \frac{2}{3} & \frac{4}{3} \\
 0 & \frac{4}{3} & \frac{2}{3}
\end{pmatrix}$. $\mu = 0$. 
Surgery diagram $\includegraphics[height=1.5cm]{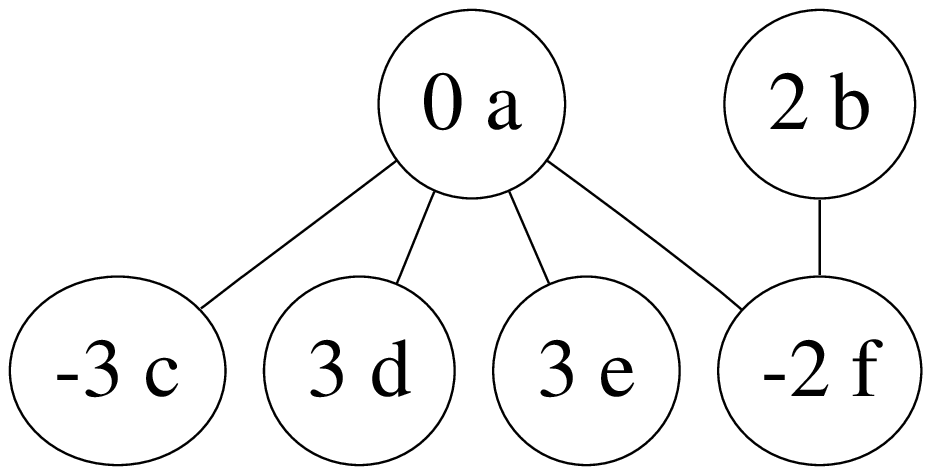}$
Embeds via Construction \ref{doubleslice}, manifold is the double branch cover of $(S^3,L_{3,2})$.   

\vskip 5mm
\centerline{$\star$ $SL_2\Real$-manifolds with infinite $H_1$ $\star$}
\vskip 5mm

All three of the manifolds below admit embeddings into $S^4$ by Lemma 3.2 of 
Crisp and Hillman \cite{CrispH}. 

\item SFS $\left[T: \frac{1}{2}\right]$, $H_1 = \Zed^2$. \ttc
\item SFS $\left[T: \frac{1}{3}\right]$, $H_1 = \Zed^2$. 
\item SFS $\left[T: \frac{1}{4}\right]$, $H_1 = \Zed^2$. 
\vskip 5mm
\centerline{$\star$ $H^2 \times \Real$ manifolds $\star$}
\vskip 5mm

\item \label{item22} SFS $\left[S^2: \frac{1}{2}, \frac{1}{2}, \frac{1}{3}, -\frac{4}{3}\right]$ 
$\Sigma_2 \times_{\Zed_6} S^1$, $H_1 = \Zed$. \ttc 
Has `deform-spun' embedding see Construction \ref{fibre-const}. 
Specifically, the genus 2 surface can be realized as a regular neighbourhood of the graph 
$G = \{ (z_1,z_2) \in \C^2 : z_1^3 \in \Real, 0 \leq z_1^3 \leq 1, z_2 = \pm \sqrt{1-|z_2|^2} \}$.
The monodromy is given by the order $6$ automorphism of $S^3$, $(z_1,z_2) \longmapsto (e^\frac{2\pi i}{3}z_1,e^{\pi i}z_2)$.

\item SFS $\left[S^2: \frac{1}{3}, \frac{1}{3}, \frac{2}{3}, -\frac{4}{3}\right]$, 
$\Sigma_2 \rtimes_{\Zed_3} S^1$, $H_1 = \Zed \oplus \Zed_3^2$. 
The surface is the same as the previous case, 
but the monodromy is given by $(z_1,z_2) \longmapsto (e^\frac{2\pi i}{3}z_1,z_2)$ 
which also allows us to realize the manifold via a deform-spun embedding.

\item \label{item24} SFS $\left[S^2: \frac{1}{2}, \frac{1}{2}, \frac{2}{5}, -\frac{7}{5}\right]$ 
$\Sigma_4 \rtimes_{\Zed_{10}} S^1$, $H_1 = \Zed$.
Consider the graph in $S^3$ given by 
$G=\{(z_1,z_2) \in \C^2 : z_1^5 \in \Real, 0 \leq z_1^5 \leq 1, z_2 = \pm \sqrt{1-|z_1|^2}\}$. 
There is a symmetry of $S^3$ of order $10$ preserving this graph 
$(z_1,z_2) \longmapsto (e^{\frac{2 \pi i}{5}}z_1,e^{\pi i}z_2)$. A surface of genus $4$ is the boundary of an equivariant regular neighbourhood of $G$ realizing the monodromy.
\vskip 5mm

\centerline{$\star$ Hyperbolic manifolds $\star$}
\vskip 5mm
\item Hyp 2.10758135 $H_1 = \Zed_5^2$.  $0$-surgery on the link $\KR{T}{10a_{114}}$.
Surgery presentation found via SnapPea.

$\includegraphics[width=3cm]{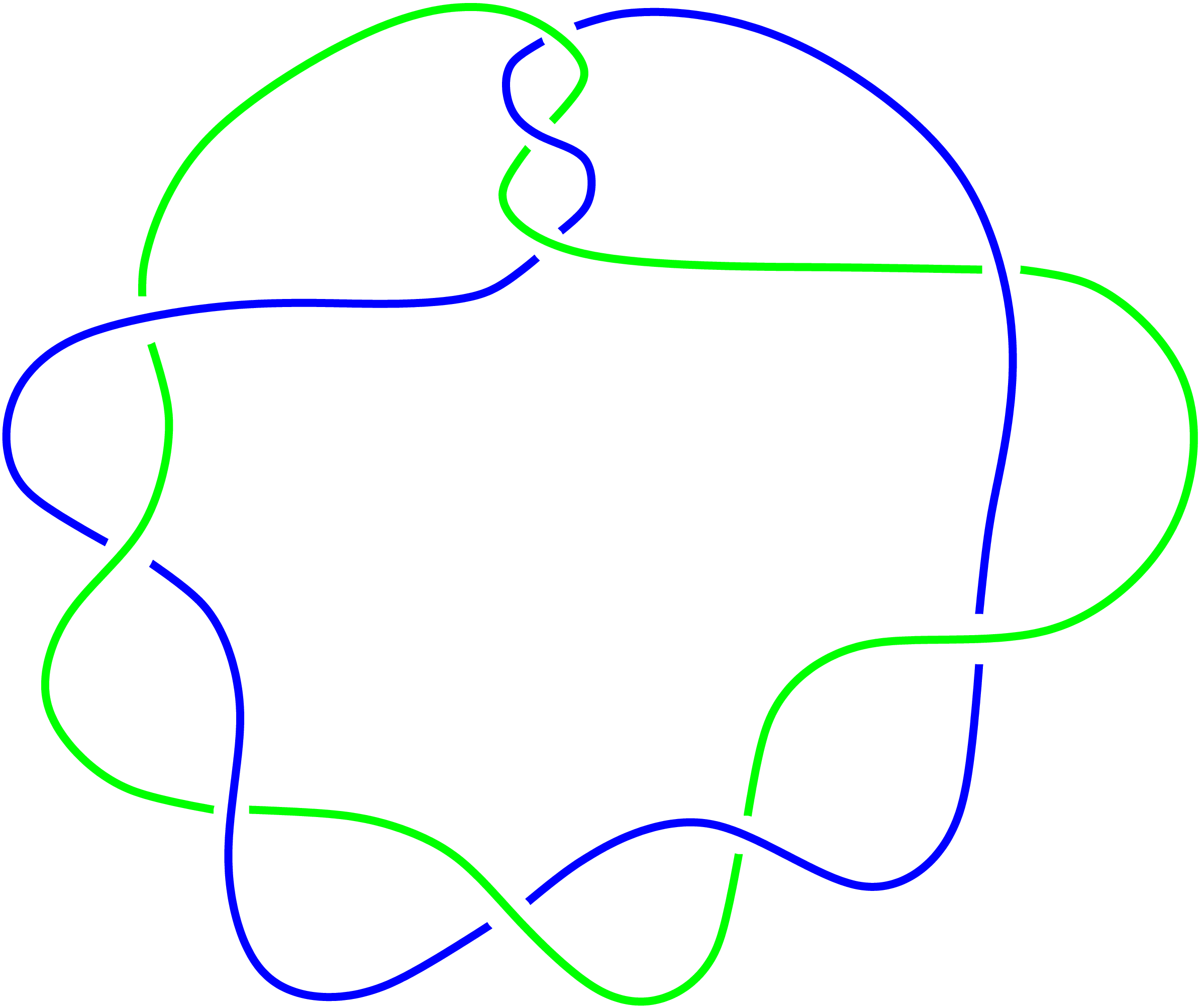}$ 

\item  Hyp 2.25976713, Homology sphere.  $0$-surgery on $\KR{T}{7a_6}$.
Surgery presentation found via SnapPea.

$\includegraphics[width=3cm]{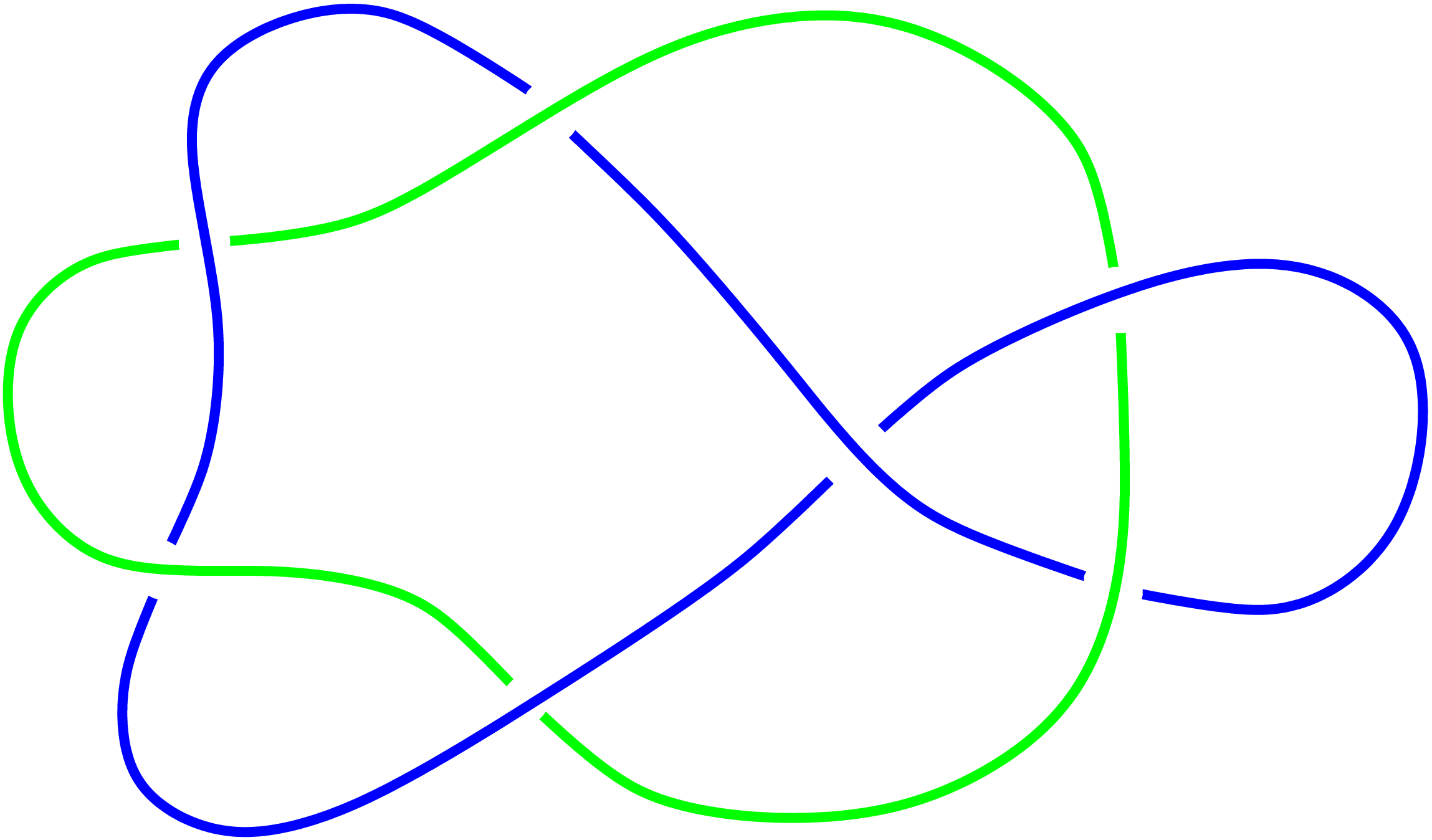}$ 

\item \label{cheg} Hyp 1.39850888, Homology sphere. 
$1$-surgery on $\KR{R}{6_1}$, also known as Stevedore's knot, which is smooth slice. 

$$ \includegraphics[width=3cm]{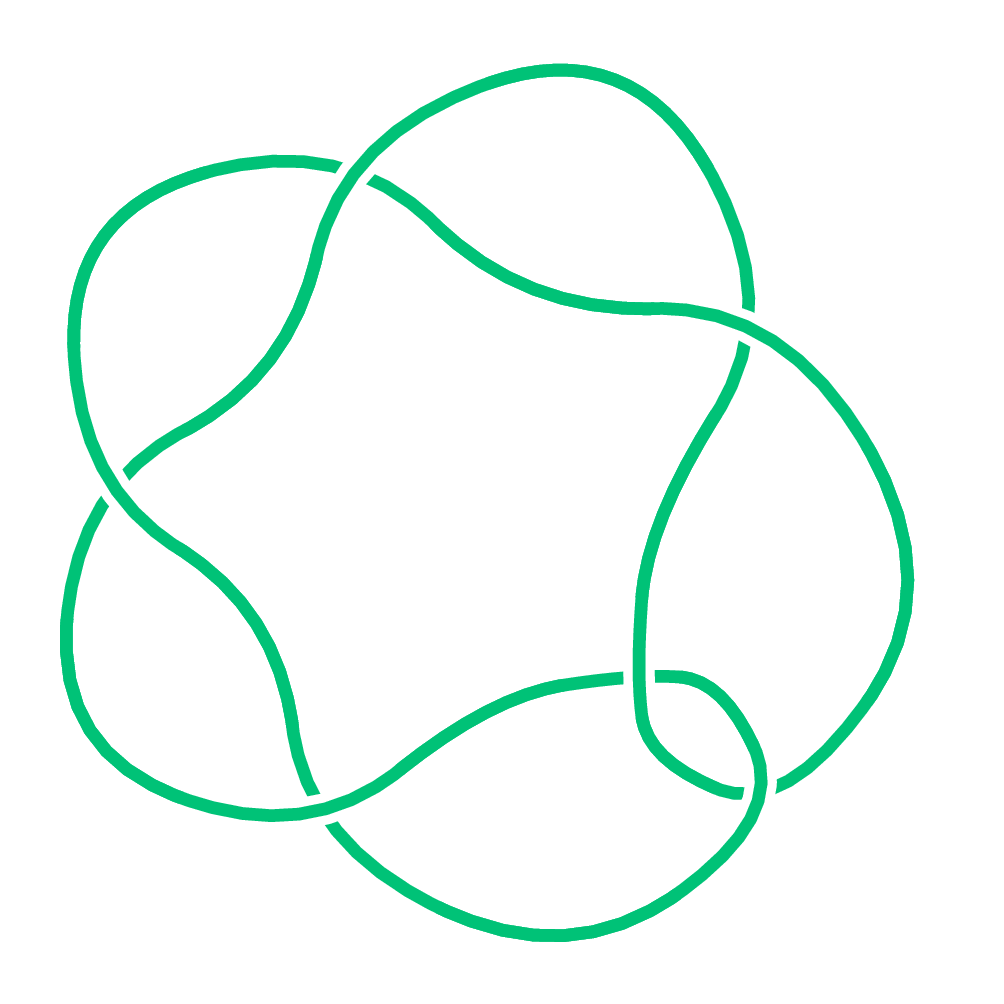} \hskip 2cm \includegraphics[width=4cm]{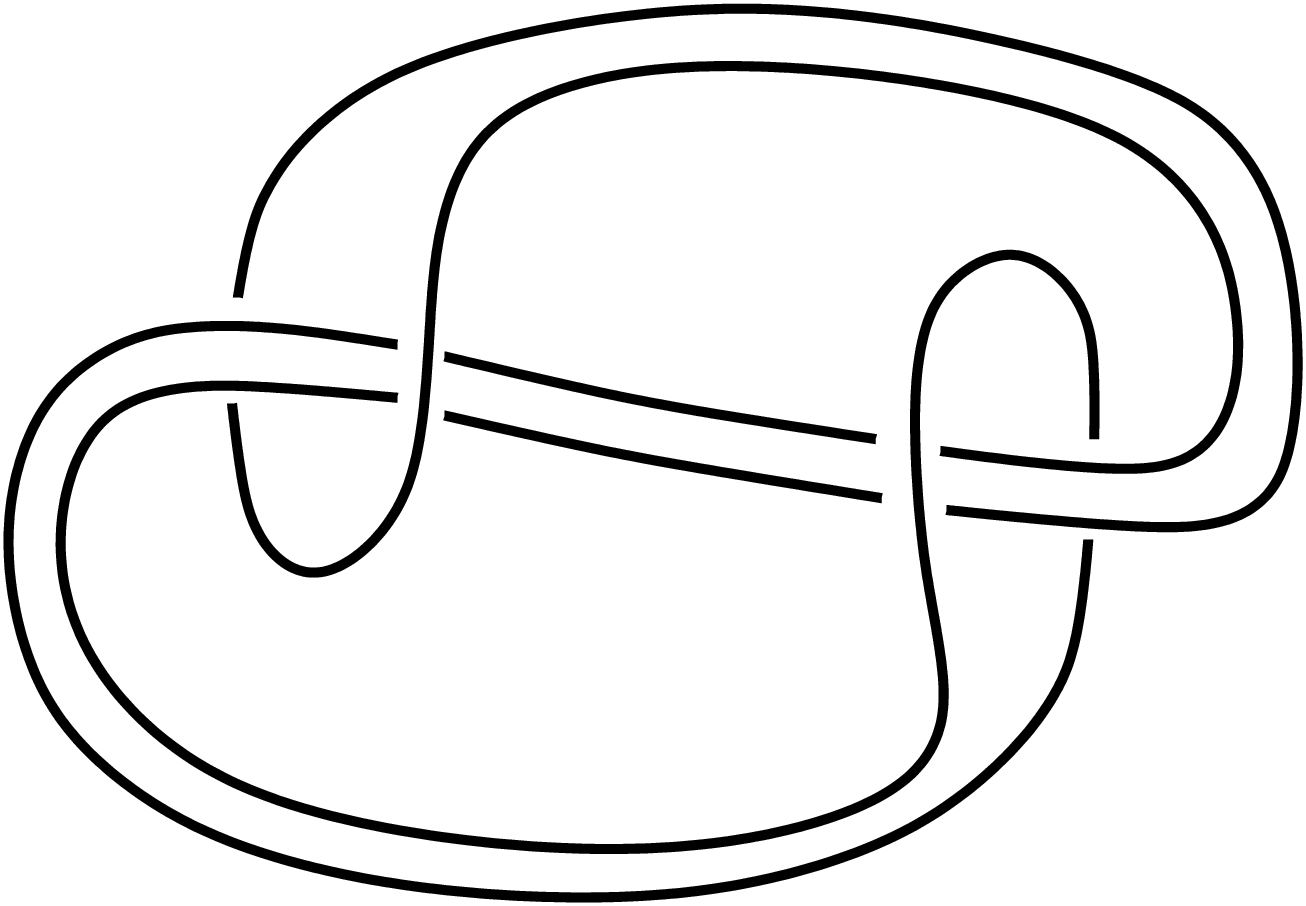} $$

By Construction \ref{tenh} this manifold embeds in a homotopy $S^4$ since it bounds a
contractible manifold $N'$.  Since $6_1$ a ribbon knot, we can apply Proposition \ref{slice-ribbon-cor}
and compute the relevant presentation of $\pi_1 N'$.  By the nature of the ribbon diagram above, 
the height function $d$ has two local minimal on the ribbon disc and one saddle point. So we have a 
presentation of the form $\langle a, b : r_1, R_1 \rangle$ where $a, b$ correspond to the
local minima of $d$ on the the ribbon disc (which also correspond to the two ribbon singularities of
the ribbon disc projected into $S^3$).  $r_1$ corresponds to the saddle, which is at the fixed point
of the symmetry of the ribbon disc, and $R_1$ to the surgery framing curve. 
So $r_1$ is the relation $a^{-1}ba = b^{-1}ab$ and $R_1$ is the relation $b=1$. 
Since $\langle a, b | a^{-1}bab^{-1}a^{-1}b, b\rangle$ is
trivializable by Andrews-Curtis moves, our manifold embeds smoothly in $S^4$. 

\item \label{emb_num} Hyp 1.91221025, Homology sphere.  
$(-1)$-surgery on $\KR{R}{8_{20}}$ which is smooth
slice, so by Construction \ref{tenh} this manifold embeds in a homotopy $S^4$.
As with item \ref{cheg} we have a ribbon diagram so we can apply Proposition \ref{slice-ribbon-cor}. 

$$\includegraphics[width=3cm]{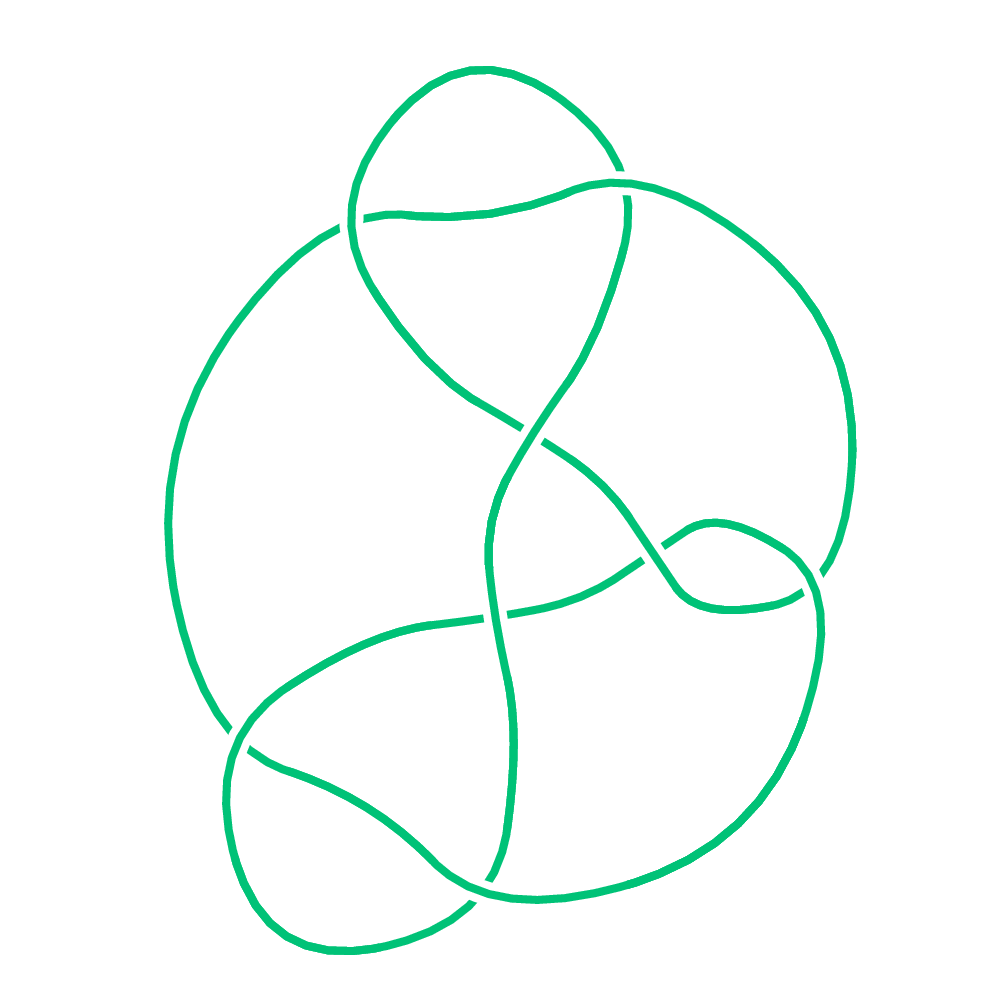} \hskip 2cm \includegraphics[width=4cm]{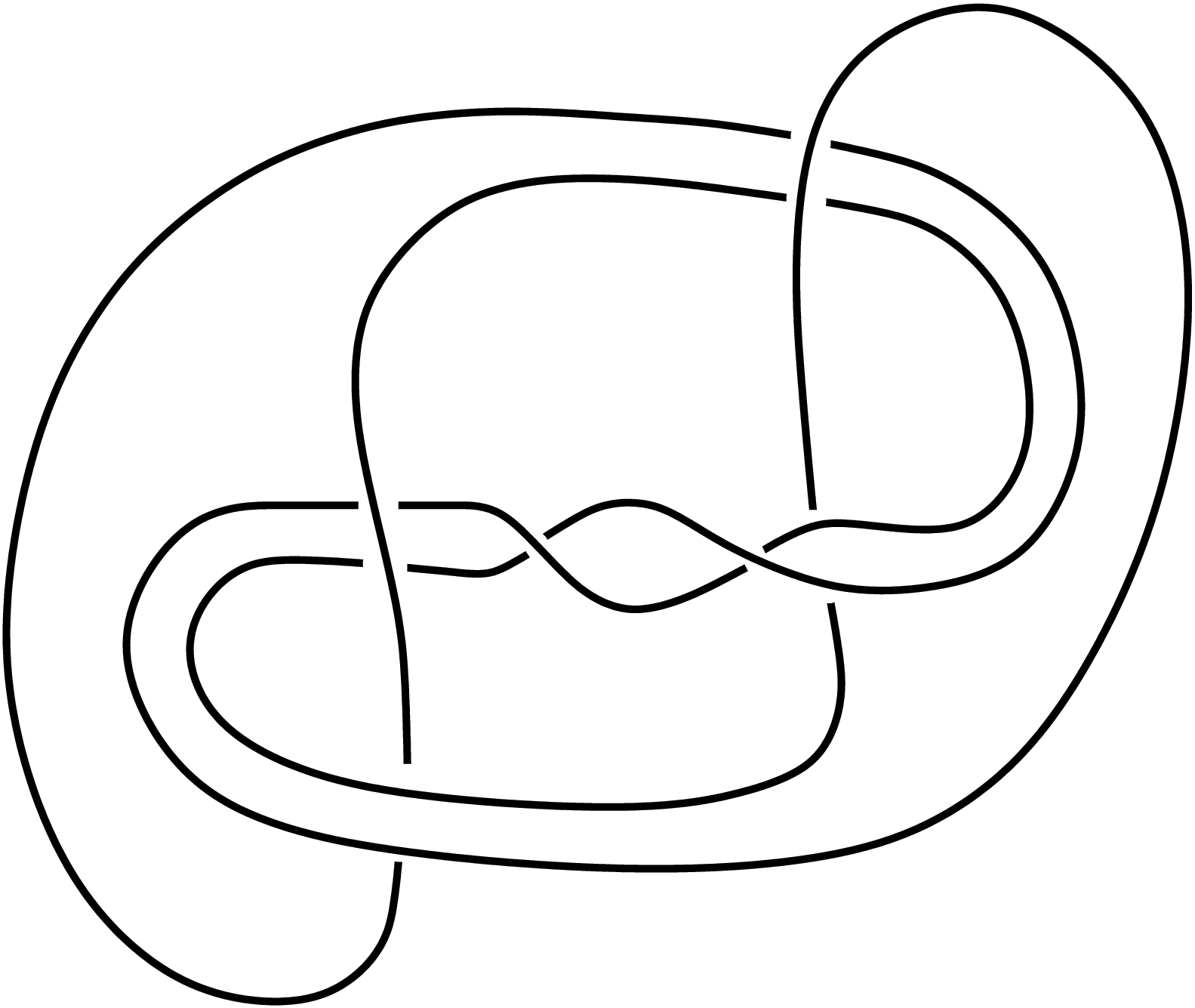} $$

This gives a similar presentation for $\pi_1 N' = \langle a, b | bab^{-1}=a^{-1}ba, b \rangle$, 
also trivializable by Andrews-Curtis moves.
\vskip 5mm

\end{enumerate}

\section{Manifolds which embed in homotopy $4$-spheres}\label{emb_htpy}

This is a list of manifolds that embed in homotopy $4$-spheres.  Likely
these homotopy $4$-spheres are diffeomorphic to $S^4$ but this has not
been determined. 

%

\begin{enumerate}
\item SFS $\left[S^2: \frac{1}{2}, \frac{1}{3}, -\frac{21}{25}\right] = \Sigma(2,3,25)$.
Although Fickle claims \cite{Fickle} that Casson and Harer \cite{CH} were the
first to show $\Sigma(2,3,25)$ bounds a contractible manifold, his Corollary 3.3 \cite{Fickle}
is the earliest written account that I have found.
\vskip 5mm
\centerline{$\star$ Hyperbolic manifolds $\star$}
\vskip 5mm

\item Hyp 1.26370924  $H_1 = \Zed_5^2$. \ttc $(-5,-5)$-surgery on $\KR{T}{5a_1}$ found via SnapPea.
$\mu = 0$, computed via the formulae in \S 4.2.3 of \cite{Sav}. 

$\includegraphics[width=3cm]{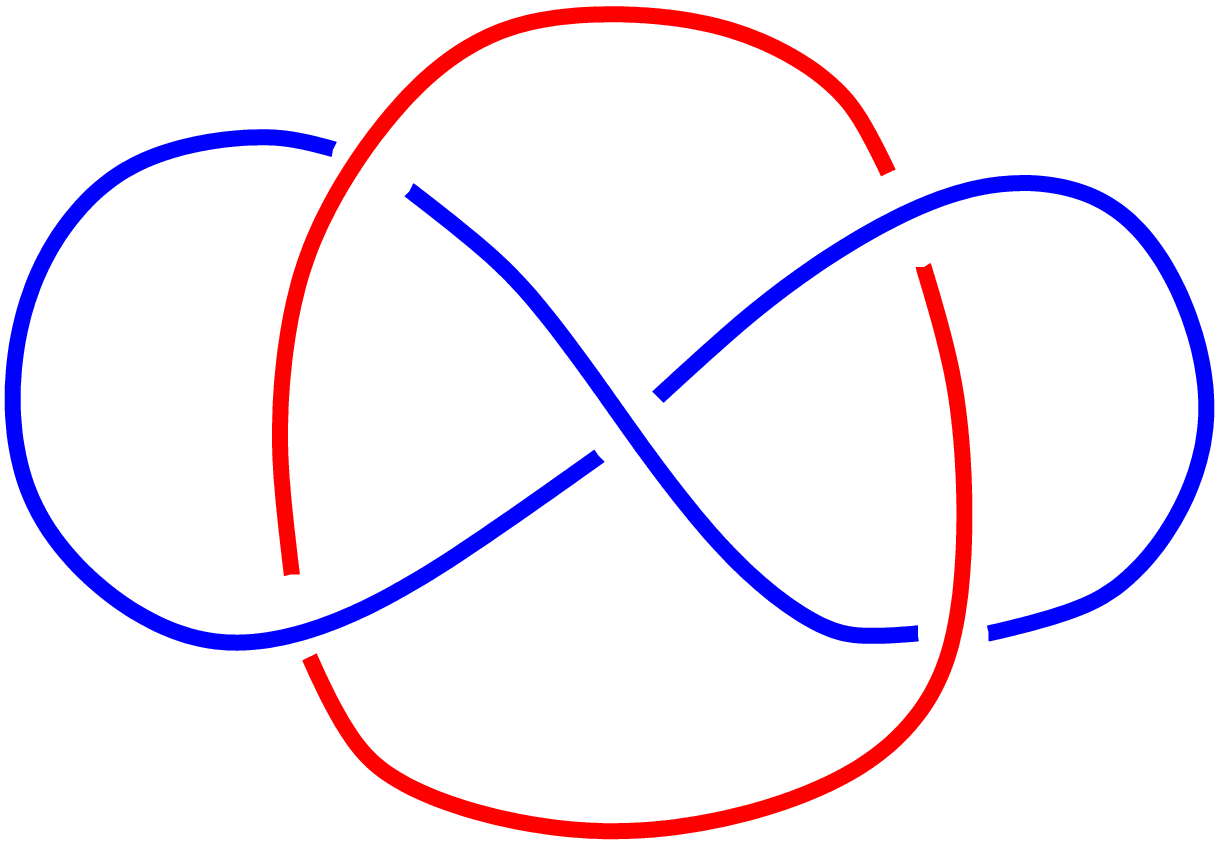}$

We use the technique of Casson and Harer \cite{CH} to embed this manifold in a homotopy $S^4$. 
A sketch is given in \S \ref{hyperbolicity}.

\vskip 5mm
\centerline{$\star$ Compound manifolds $\star$}
\vskip 5mm

\item SFS $\left[D: \frac{1}{2}, \frac{1}{3}\right]$ U/m 
      SFS $\left[D: \frac{1}{2}, \frac{2}{3}\right]$ 
$m = \begin{pmatrix}  0 & 1 \\ 1 & 0 \end{pmatrix}$ Homology sphere. \ttc

Found as a vertex-normal $3$-manifold in a $10$-pentachoron triangulated homotopy $4$-sphere with `isomorphism signature' {\tiny isoSig(): kLLLLAQQQcccdfhfjjgjgjiiiiPayaPbaaPaaaaaaadaPbcatbcaPbyaca}. 

\item \label{emb_homolsph_num} SFS $\left[D: \frac{1}{2}, \frac{1}{2}\right]$ U/m 
      SFS $\left[D: \frac{1}{3}, \frac{3}{5}\right]$ 
$m = \begin{pmatrix}  0 & 1 \\ 1 & 0\end{pmatrix}$  $H_1 = \Zed_2^2$.

Found as a vertex-normal $3$-manifold in a $6$-pentachoron triangulated homotopy $4$-sphere with `isomorphism signature' {\tiny isoSig(): gLLAQQcdbdefeeff4aYaYaYagaEaKaKaiaga}. 
\end{enumerate}

\section{Manifolds in the census known to not embed in $S^4$}\label{nonembeddable_man}

\begin{enumerate}
\item SFS $\left[S^2: \frac{1}{2}, \frac{1}{2}, -\frac{5}{6}\right] = 
SFS \left[\RProj^2 : 6 \right] = S^3 / Q_{24}$. 
Crisp-Hillman \cite{CrispH}. \ttc
\item SFS $\left[S^2: \frac{1}{2}, \frac{1}{2}, -\frac{9}{10} \right] = S^3/Q_{40}$. 
$H_1 = \Zed_2^2$. Crisp-Hillman \cite{CrispH}.
\item SFS $\left[S^2: \frac{1}{2}, \frac{1}{3}, -\frac{5}{6}\right] = 
(S^1 \times S^1) \times_{\Zed_6} S^1$. Crisp-Hillman \cite{CrispH}. \ttc
\item $(S^1 \times S^1) \rtimes_{\begin{pmatrix} 2 & 1 \\ 1 & 1 \end{pmatrix}} S^1$. 
$H_1 = \Zed$. Crisp-Hillman \cite{CrispH}. \ttc
\item SFS $\left[\RProj^2/n2: \frac{1}{2}, \frac{11}{2}\right]$. 
$H_1 = \Zed_4^2$ Nil-manifold. Crisp-Hillman \cite{CrispH}.
\item $(S^1 \times S^1) \rtimes_{\begin{pmatrix} 10 & 3 \\ 3 & 1 \end{pmatrix}} S^1$,
$H_1 = \Zed \oplus \Zed_3^2$. Crisp-Hillman \cite{CrispH}.
\item SFS $\left[D: \frac{1}{2}, \frac{1}{2}\right]$ U/m 
SFS $\left[D: \frac{1}{2}, \frac{1}{2}\right]$ 
$\begin{pmatrix}  -5 & 7 \\ -2 & 3 \end{pmatrix}$ $H_1 = \Zed_4^2$ Sol manifold.
Crisp-Hillman \cite{CrispH}.
\vskip 5mm
\centerline{$\star$ $H^2$-fibre geometry $\star$}
\vskip 5mm
These manifolds fibre over $S^1$ with fibre a hyperbolic surface, and monodromy an automorphism
of finite order.  Regina stores these manifolds via their Seifert data, see the item on computing
the monodromy from the Seifert data for details on how we compute the Alexander polynomials of these
manifolds in \S \ref{hyperbolicity}.

\item \label{it8} SFS $\left[S^2: \frac{1}{2}, \frac{1}{5}, -\frac{7}{10}\right] = \Sigma_2 \rtimes_{\Zed_{10}} S^1$. $H_1 = \Zed$
$\Delta = t^2 - t + 1 - t^{-1} + t^{-2} \neq p(t)p(t^{-1})$ 
\item SFS $\left[S^2: \frac{1}{2}, \frac{2}{7}, -\frac{11}{14}\right] = \Sigma_3 \rtimes_{\Zed_{14}} S^1$. $H_1 = \Zed$
$\Delta = t^3 - t^2 + t - 1 + t^{-1} - t^{-2} + t^{-3} \neq p(t)p(-t)$
\item SFS $\left[S^2: \frac{1}{3}, \frac{1}{4}, -\frac{7}{12}\right] = \Sigma_3 \rtimes_{\Zed_{12}} S^1$.
$H_1 = \Zed$ $\Delta = t^3-t^2+1-t^{-2}+t^{-3} \neq p(t)p(t^{-1})$
\item \label{it11} SFS $\left[S^2: \frac{1}{3}, \frac{2}{5}, -\frac{11}{15}\right] = \Sigma_4 \rtimes_{\Zed_{15}} S^1$.
$H_1 = \Zed$ $\Delta = t^4 - t^3 + t -1 + t^{-1} -t^{-3} +t^{-4} \neq p(t)p(t^{-1})$
\vskip 5mm
In a recent preprint, Jonathan Hillman \cite{Hill2} proves that $H^2 \times \Real$ manifolds
that fibre over $S^2$ must have an even number of singular fibres, generalizing 
items \ref{it8}--\ref{it11}.  He also uses the Alexander module as an obstruction. 
\vskip 5mm

\centerline{$\star$ Homology spheres with non-zero Rochlin invariant $\star$}
\vskip 5mm
These do not embed because they do not satisfy the Rochlin invariant test.
See Theorem \ref{spintests}. The $\mubar$ invariant was computed using
formula 2.4.2 in Saveliev's text \cite{Sav}.

\item SFS $\left[S^2: \frac{1}{2}, \frac{1}{3}, -\frac{4}{5}\right]$ $S^3/P_{120}$ 
Poincar\'e Dodecahedral Space. $\mubar = -1$. \ttc
\item SFS $\left[S^2: \frac{1}{2}, \frac{1}{3}, -\frac{6}{7}  \right]$ 
Brieskorn homology sphere. $\mubar = 1$, $d=0$. \ttc
\item SFS $\left[S^2: \frac{1}{2}, \frac{1}{3}, -\frac{14}{17} \right]$ 
Brieskorn homology sphere. $\mubar = -1$, $d=2$. \ttc
\item SFS $\left[S^2: \frac{1}{3}, \frac{1}{4}, -\frac{4}{7}   \right]$ 
Brieskorn homology sphere. $\mubar = -1$, $d=2$. \ttc
\item SFS $\left[S^2: \frac{1}{2}, \frac{1}{3}, -\frac{16}{19} \right]$ 
Brieskorn homology sphere. $\mubar = 1$, $d=0$.
\item SFS $\left[S^2: \frac{1}{2}, \frac{1}{3}, -\frac{24}{29} \right]$ 
Brieskorn homology sphere. $\mubar = -1$, $d=2$.
\item SFS $\left[S^2: \frac{1}{2}, \frac{1}{5}  -\frac{12}{17} \right]$ 
Brieskorn homology sphere. $\mubar = 1$, $d=0$.
\item SFS $\left[S^2: \frac{1}{3}, \frac{1}{4}  -\frac{10}{17} \right]$ 
Brieskorn homology sphere. $\mubar = 1$,  $d=0$.
\vskip 5mm

\centerline{$\star$ Brieskorn homology spheres with non-zero $d$-invariant $\star$}
\vskip 5mm
These manifolds fail the $d$-invariant test, see Theorem \ref{spintests}.

\item SFS $\left[S^2: \frac{1}{2}, \frac{1}{3}, -\frac{9}{11}\right]$, $\mubar = 0$, $d=2$. \ttc
\item SFS $\left[S^2: \frac{1}{2}, \frac{1}{3}, -\frac{19}{23}\right]$, $\mubar = 0$, $d=2$.
\item SFS $\left[S^2: \frac{1}{2}, \frac{2}{7}, -\frac{7}{9}\right]$, $\mubar = 0$, $d=2$.
\item SFS $\left[S^2: \frac{1}{3}, \frac{2}{5}, -\frac{8}{11}\right]$, $\mubar = 0$, $d=2$.
\item SFS $\left[S^2: \frac{1}{2}, \frac{1}{5}, -\frac{9}{13}\right]$, $\mubar = 0$, $d=2$.
\item SFS $\left[S^2: \frac{1}{3}, \frac{1}{4}, -\frac{11}{19}\right]$, $\mubar = 0$, $d=2$.
\item SFS $\left[S^2: \frac{1}{3}, \frac{2}{7}, -\frac{8}{13}\right]$, $\mubar = -2$, $d=2$.
\vskip 5mm

\centerline{$\star$ Rational homology spheres which do not satisfy the $\vec d$ test $\star$}
\vskip 5mm

\item SFS $\left[S^2: \frac{1}{3}, \frac{2}{3}, -\frac{5}{6}\right]$ 
$H_1 = \Zed_3^2$. \ttc 
$\vec d = 
\begin{pmatrix} 
 0 & 2 & 2 \\ 
 0 & 2/3 & 4/3 \\
 0 & 4/3 & 2/3 
\end{pmatrix}$ 
see Corollary \ref{vectest}.

\item SFS $\left[S^2: \frac{1}{3}, \frac{1}{3}, -\frac{7}{12}\right]$ 
$H_1 = \Zed_3^2$. \ttc 
$\vec d = 
\begin{pmatrix} 
 2 & 0 & 0 \\ 
 0 & 4/3 & 2/3 \\
 0 & 2/3 & 4/3 
\end{pmatrix}$ see Corollary \ref{vectest}.

\item SFS $\left[S^2: \frac{1}{5}, \frac{2}{5}, -\frac{2}{5}\right]$  $H_1 = \Zed_5^2$.
$\vec d = \begin{pmatrix} 
 0 & 0 & 0 & 0 & 0 \\ 
 0 & 2/5 & 4/5 & 6/5 & 8/5 \\
 2 & 4/5 & 8/5 & 2/5 & 6/5 \\
 2 & 6/5 & 2/5 & 8/5 & 4/5 \\
 0 & 8/5 & 6/5 & 4/5 & 2/5
\end{pmatrix}$ 
see Corollary \ref{vectest}.

\item SFS $\left[S^2: \frac{1}{3}, \frac{1}{3}, -\frac{19}{30}\right]$ $H_1 = \Zed_3^2$. 
$\vec d = \begin{pmatrix} 
 2 & 0 & 0 \\ 
 0 & 4/3 & 2/3 \\
 0 & 2/3 & 4/3
\end{pmatrix}$ see Corollary \ref{vectest}.

\item SFS $\left[S^2: \frac{1}{2}, \frac{1}{2}, \frac{1}{3}, -\frac{13}{10}\right]$ $H_1 = \Zed_2^2$.  
$\vec d = \begin{pmatrix} 
 2 & 0 \\ 
 1 & 0 
\end{pmatrix}$ see Corollary \ref{vectest}.

\item SFS $\left[S^2: \frac{1}{3}, \frac{1}{3}, -\frac{13}{21}\right]$ $H_1 = \Zed_3^2$.
$\vec d = \begin{pmatrix} 
 2 & 0 & 0 \\ 
 0 & 4/3 & 2/3 \\
 0 & 2/3 & 4/3
\end{pmatrix}$ see Corollary \ref{vectest}.
\vskip 5mm

\centerline{$\star$ Rational homology spheres that do not satisfy that $\vec \mu$-test $\star$}
\vskip 5mm

\item SFS $\left[S^2: \frac{1}{2}, \frac{1}{6}, -\frac{7}{10}\right]$ $H_1 = \Zed_2^2$.  
Characteristic links: $(\{a,b,c\},\{a,b,d\},\{e\}, \{c,d,e\})$, 
$\vec \mu = (0,\frac{1}{2}, -\frac{1}{2}, 0)$,
$\vec d = \begin{pmatrix}  0 & 0 \\  -1 & 0  \end{pmatrix}$, 
surgery diagram: \includegraphics[height=1.5cm]{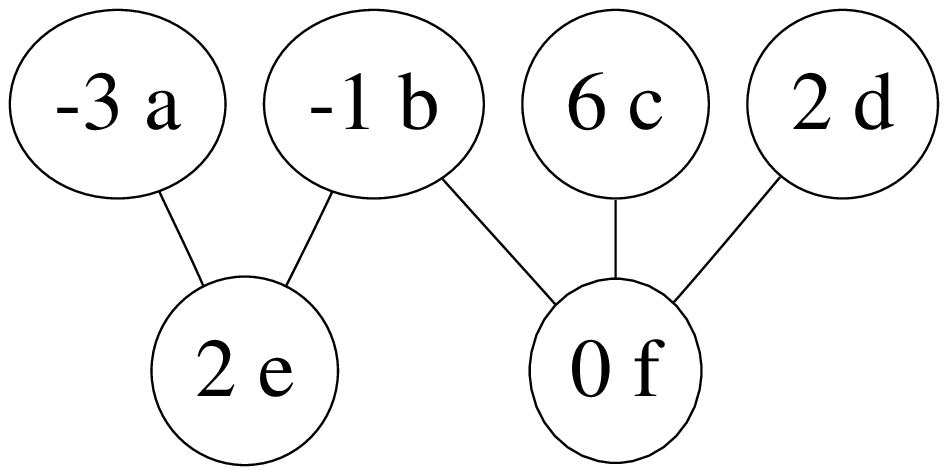}.

\item SFS $\left[S^2: \frac{1}{3}, \frac{1}{3}, -\frac{5}{6}\right]$ $H_1 = \Zed_3^2$. \ttc
Characteristic link $\{a,d,e\}$, $\mu = -\frac{3}{4}$,
$\vec d = \begin{pmatrix} 0 & 0 & 0 \\  \frac{2}{3} & 0 & -\frac{2}{3} \\ \frac{2}{3} & -\frac{2}{3} & 0 \\
\end{pmatrix}$, 
surgery diagram: \includegraphics[height=1.5cm]{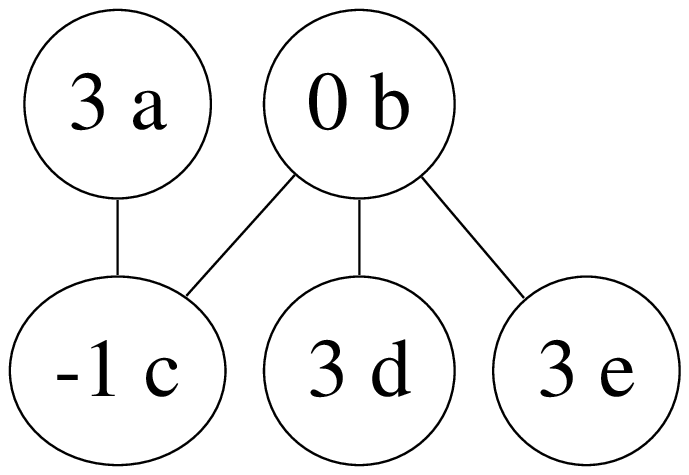}.

\item SFS $\left[S^2: \frac{1}{3}, \frac{1}{3}, -\frac{11}{15}\right]$ $H_1 = \Zed_3^2$. 
Characteristic link $\{c,e,f\}$, $\mu = -\frac{3}{4}$, 
$\vec d = \begin{pmatrix}  0 & 0 & 0 \\  0 & -\frac{2}{3} & \frac{2}{3} \\ 0 & \frac{2}{3} & -\frac{2}{3} \\
\end{pmatrix}$, 
surgery diagram: \includegraphics[height=1.5cm]{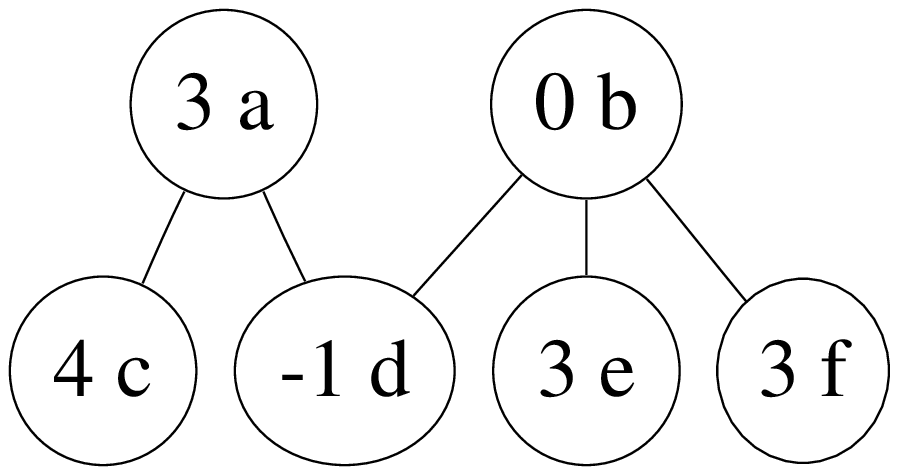}.

\item SFS $\left[S^2: \frac{1}{3}, \frac{1}{3}, -\frac{17}{24}\right]$ $H_1 = \Zed_3^2$, 
Characteristic link $\{a,b,f,g\}$, $\mu = 1$,
$\vec d = \begin{pmatrix} 
 0 & 0 & 0 \\  0 & -\frac{2}{3} & \frac{2}{3} \\ 0 & \frac{2}{3} & -\frac{2}{3} \\
\end{pmatrix}$,
surgery diagram: \includegraphics[height=1.5cm]{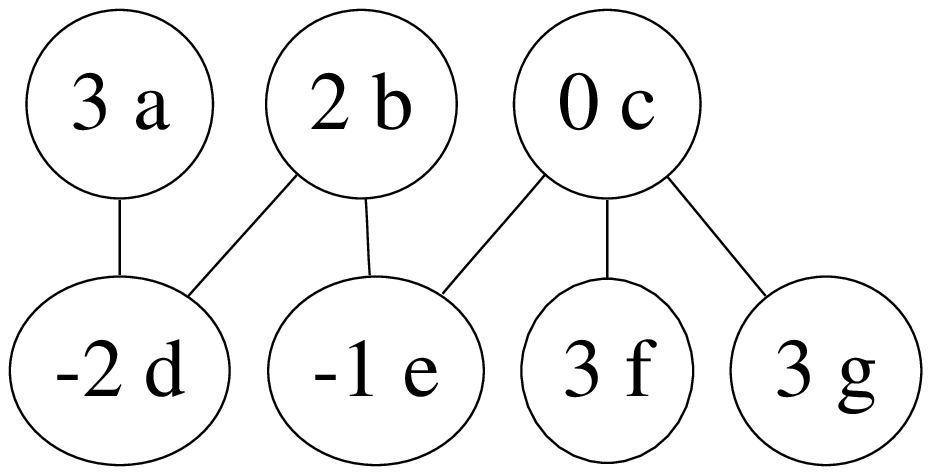}.

\vskip 5mm

\centerline{$\star$ Other rational homology spheres $\star$}
\vskip 5mm

\item SFS $\left[\RProj^2/n2: \frac{1}{4}, -\frac{1}{4}\right]$ $H_1 = \Zed_8^2$. Crisp-Hillman Proposition 1.2 \cite{CrispH}.
\vskip 5mm
\centerline{$\star$ Hyperbolic manifolds $\star$}
\vskip 5mm

\item Hyp 1.73198278, Homology sphere. $\mu = 1$. 
$+\frac{1}{3}$-surgery on $\KR{R}{4_1}$ (found via SnapPea).  $\mu$ is computed using the surgery
formula (Theorem 2.8 of \cite{Sav}). 

{
\psfrag{1/3}[tl][tl][1][0]{$\frac{1}{3}$}
$\includegraphics[width=3cm]{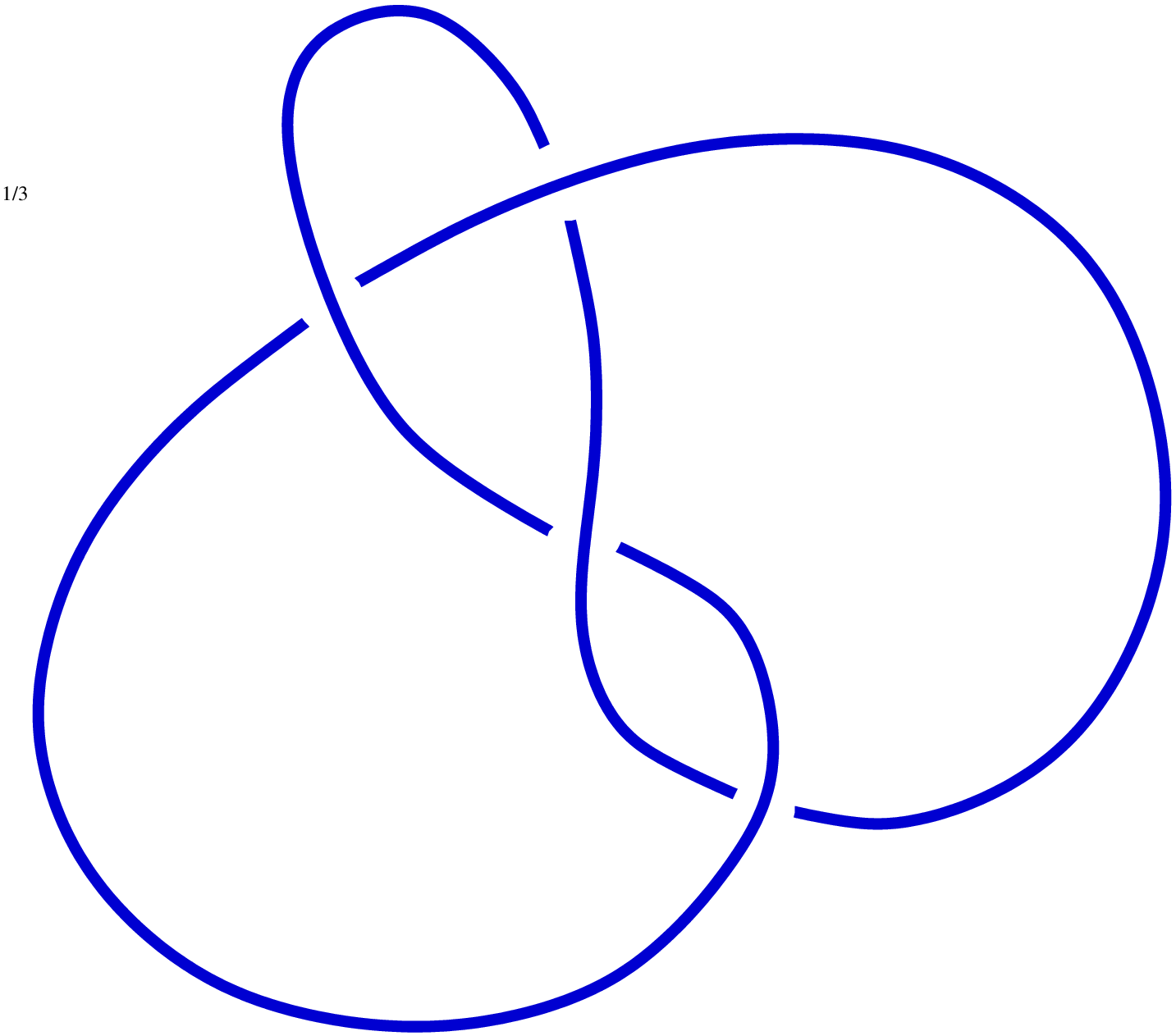}$ 
}
\vskip 5mm

\centerline{$\star$ Manifolds with non-trivial JSJ-decompositions $\star$}
\vskip 5mm
These manifolds are all of the form $SFS \left[ A: \frac{\alpha}{\beta} \right] / 
\begin{pmatrix} a & b \\ c & d \end{pmatrix}$ where $ad-bc = -1$. These manifolds have
$b_1 = 1$ if and only if the polynomial $\beta t^2 + ( (d-a)\beta - b\alpha)t + \beta$
does not have a zero at $t=1$, moreover, if $b_1=1$, this polynomial is the Alexander
polynomial of the corresponding covering space. Checking that this polynomial has the form
$rp(t)p(t^{-1})$ where $p(t)$ is a rational Laurent polynomial and $r$ is rational amounts to
determining if the number $((a-d)+\frac{b\alpha}{\beta})^2-4$ is a rational squared.
These five manifolds do not embed since their Alexander polynomials do not satisfy the Kawauchi
condition. See Theorem \ref{KawCond}.

\item SFS $\left[A: \frac{1}{2}\right]$ / $\begin{pmatrix}  0 & 1 \\ 1 & -1\end{pmatrix}$ $H_1 = \Zed$. \ttc
\item SFS $\left[A: \frac{1}{3}\right]$ / $\begin{pmatrix}  0 & 1 \\ 1 & -2\end{pmatrix}$ $H_1 = \Zed$.
\item SFS $\left[A: \frac{2}{3}\right]$ / $\begin{pmatrix}  0 & 1 \\ 1 & -1\end{pmatrix}$ $H_1 = \Zed$.
\item SFS $\left[A: \frac{1}{4}\right]$ / $\begin{pmatrix}  0 & 1 \\ 1 & -2\end{pmatrix}$ $H_1 = \Zed$.
\item SFS $\left[A: \frac{3}{4}\right]$ / $\begin{pmatrix}  0 & 1 \\ 1 & -1\end{pmatrix}$ $H_1 = \Zed$.

$\mubar$ is computed for the examples below using the splicing additivity 
formula for $\mubar$, Proposition 2.16 from \cite{Sav}.

\item SFS $\left[D: \frac{1}{2}, \frac{1}{3}\right]$ U/m 
      SFS $\left[D: \frac{1}{2}, \frac{5}{7}\right]$ 
$m = \begin{pmatrix}  0 & 1 \\ 1 & 0 \end{pmatrix}$  Homology sphere.
$\Sigma(2,3,5) \splice_{5,5} \Sigma(2,5,7)$, $\mubar = -1$.

\item SFS $\left[D: \frac{1}{2}, \frac{1}{3}\right]$ U/m 
      SFS $\left[D: \frac{2}{3}, \frac{3}{5}\right]$ 
$m = \begin{pmatrix}  -1 & 2 \\ 0 & 1 \end{pmatrix}$  Homology sphere.
$\Sigma(2,3,5) \splice_{5,4} \Sigma(3,4,5)$, $\mubar = -1$.

\item SFS $\left[D: \frac{1}{2}, \frac{2}{3}\right]$ U/m 
      SFS $\left[D: \frac{1}{2}, \frac{4}{11}\right]$ 
$m = \begin{pmatrix}  0 & 1 \\ 1 & 0 \end{pmatrix}$  Homology sphere.
$\Sigma(2,7,11) \splice \Sigma(2,3,19)$ 
Issa-McCoy obstruction \cite{Issa}. 

\item SFS $\left[D: \frac{1}{2}, \frac{2}{3}\right]$ U/m 
      SFS $\left[D: \frac{2}{3}, \frac{1}{5}\right]$ 
$m = \begin{pmatrix}  0 & 1 \\ 1 & 0 \end{pmatrix}$  Homology sphere.
$\Sigma(3,5,7) \splice \Sigma(2,3,13)$ 
Issa-McCoy obstruction \cite{Issa}.

\item SFS $\left[D: \frac{1}{2}, \frac{3}{5}\right]$ U/m 
      SFS $\left[D: \frac{2}{3}, \frac{1}{4}\right]$ 
$m = \begin{pmatrix}  0 & 1 \\ 1 & 0 \end{pmatrix}$  Homology sphere.
$\Sigma(2,5,11) \splice \Sigma(3,4,11)$ 
Issa-McCoy obstruction \cite{Issa}. 

In the next few examples we need to compute the Alexander polynomials of some
graph manifolds.  The underlying Seifert-fibred manifolds are all of the type
$SFS \left[ D: \frac{a}{b}, \frac{c}{d}\right]$.  An elementary computation shows that
$$ H_1 SFS \left[ D: \frac{a}{b}, \frac{c}{d}\right] \simeq \Zed \oplus \Zed_{GCD(b,d)}.$$
These manifolds fibre over $S^1$ -- the horizontal incompressible surface is the fibre. 
Moreover, since these manifolds fibre over a disc with two singular fibres, the monodromy
can be realized as the covering transformation of a surface such that the quotient orbifold
is a disc with two cone points. This gives an immediate Mayer-Vietoris computation of the
Alexander polynomial, considering it as the order ideal of the homology of the fibre (of the
fibring over $S^1$).

\begin{lem}\label{alexlem}
Consider a manifold $M \cup_T N$ which is the union of two submanifolds $M$ and $N$ along a 
common boundary torus $T$.  Assume $M \cup_T N$ is a rational homology $S^1 \times S^2$,
and both $M$ and $N$ are rational homology $S^1 \times D^2$ manifolds.
$$\Delta_{M\cup_T N}(t) = \frac{\Delta_{M}(t^p) \Delta_{N}(t^q)(t-1)^2}{(t^p-1)(t^q-1)}$$
where $coker(H_1 M \to fH_1 (M\cup_T N)) = \Zed_p$ and $coker(H_1 N \to fH_1 (M\cup_T N)) = \Zed_q$. 
$p$ and $q$ have a simpler computation since $coker(H_1 T \to fH_1 M) = \Zed_q$
and $coker(H_1 T \to fH_1 N) = \Zed_p$.
Moreover, 
$$ \Delta SFS \left[ D: \frac{a}{b}, \frac{c}{d}\right] = \frac{ (t^{LCM(b,d)}-1)(t-1) }{ (t^{b'}-1)(t^{d'}-1) } $$
where $b' = \frac{b}{GCD(b,d)}, d'=\frac{d}{GCD(b,d)}$.
\end{lem}

The relevant non-embedding result is Theorem \ref{KawCond}. 

\item SFS $\left[D: \frac{1}{2}, \frac{1}{2}\right]$ U/m 
      SFS $\left[D: \frac{1}{2}, \frac{1}{3}\right]$ 
$m = \begin{pmatrix}  1 & 5 \\ 1 & 4\end{pmatrix}$  $H_1 = \Zed $.
$\Delta(t) = t^4-t^2+1$.
\item SFS $\left[D: \frac{1}{2}, \frac{1}{2}\right]$ U/m 
      SFS $\left[D: \frac{1}{2}, \frac{1}{3}\right]$ 
$m = \begin{pmatrix}  5 & 1 \\ 4 & 1\end{pmatrix}$  $H_1 = \Zed$.
$\Delta(t) = t^4-t^2+1$.
\item  SFS $\left[D: \frac{1}{2}, \frac{1}{3}\right]$ U/m 
      SFS $\left[D: \frac{1}{2}, \frac{7}{10}\right]$ 
$m = \begin{pmatrix}  0 & 1 \\ 1 & 0\end{pmatrix}$ $H_1 = \Zed$.
$\Delta = (t^4-t^2+1)(t^4-t^3+t^2-t+1)$. 

\vskip 5mm
\centerline{$\star$ Fibres over $S^1$ with reducible monodromy $\star$}
\vskip 5mm

\item SFS $\left[D: \frac{1}{2}, \frac{1}{2}\right]$ U/m 
      SFS $\left[D: \frac{2}{5}, \frac{3}{5}\right]$ 
$m = \begin{pmatrix}  0 & 1 \\ 1 & 0\end{pmatrix}$ $\Sigma_4 \rtimes S^1$, $H_1 = \Zed$.
The monodromy is reducible with reduction system a union of $5$ circles separating
$\Sigma_4$ into two $5$-punctures spheres.  Perhaps the easiest way to describe the monodromy
is that it differs from the monodromy of item \ref{item24} \S \ref{embeddable_man} 
by a single Dehn twist about a reduction curve.  The Alexander polynomial for this
manifold is the same as item \ref{item24} \S \ref{embeddable_man}, so it does not
provide an obstruction to embedding. Alternatively, the monodromy extends over a handlebody
thus this manifolds bounds a genus 4 handlebody bundle over $S^1$ which must be a homology
$S^1 \times D^3$.  The obstruction to embedding is a variant of
the Crisp-Hillman Theorem \ref{chthm}. If this manifold embeds in $S^4 = V_1 \cup_M V_2$, 
then $V_1$ is a homology $S^1 \times D^3$ and $V_2$ is a homology $S^2 \times D^2$. Replace
$V_1$ with $V_1'$ the corresponding handlebody bundle over $S^1$, $W = V_1' \cup_M V_2$ is
therefore also a homology $S^4$, but it contains a Klein bottle with normal Euler class $\pm 2$,
contradicting Theorem \ref{chthm}.

\item SFS $\left[D: \frac{1}{2}, \frac{1}{2}\right]$ U/m 
      SFS $\left[D: \frac{1}{3}, \frac{2}{3}\right]$ 
$m = \begin{pmatrix}  0 & 1 \\ 1 & 0\end{pmatrix}$ $\Sigma_2 \rtimes S^1$, $H_1 = \Zed$. \ttc 
The monodromy is reducible, the reduction system of $3$ curves 
separates the genus $2$ surface into two $3$-punctures spheres.  The monodromy
differs from the monodromy of item \ref{item22} \S \ref{embeddable_man} by
a single Dehn twist about a reduction curve. Again the Alexander polynomial is the
same as in item \ref{item22} \S \ref{embeddable_man} so it is no obstruction
to embedding.  This does not embed for essentially the same reason as the previous
example, only in this case we use the appropriate genus $2$ handlebody bundle over $S^1$.

\item SFS $\left[D: \frac{1}{2}, \frac{1}{2}\right]$ U/m 
      SFS $\left[D: \frac{1}{3}, \frac{2}{3}\right]$ 
$m = \begin{pmatrix}  -2 & 3 \\ -1 & 2\end{pmatrix}$  $\Sigma_2 \rtimes S^1$, $H_1 = \Zed$.
 The monodromy is reducible with a reduction system of $3$ curves separating
the surface into two pairs of pants.  The monodromy differs from the monodromy of item \ref{item22} \S
\ref{embeddable_man} by the cube of a Dehn twist along one of the reduction curves. Thus the
manifold bounds a handlebody bundle over $S^1$.  Notice this bundle contains a Klein bottle with
normal Euler class $W_2 = \pm 6$, which does not embed in a homology $S^4$ by Theorem \ref{chthm}.
$\Delta(t) = (t^2-t+1)^2$ 

\item SFS $\left[D: \frac{1}{2}, \frac{1}{2}\right]$ U/m 
      SFS $\left[D: \frac{1}{4}, \frac{3}{4}\right]$ 
$m = \begin{pmatrix}  -1 & 2 \\ 0 & 1\end{pmatrix}$ $\Sigma_2 \rtimes S^1$, $H_1 = \Zed \oplus \Zed_2^2$.
 The monodromy is reducible with reduction system $4$ curves separating
the surface into two $4$-punctured spheres.  Like the previous examples, this bundle bounds
a handlebody bundle over $S^1$, which in this case contains a Klein bottle with normal Euler class
$\pm 2$, and so this $3$-manifold does not embed in $S^4$ by Theorem \ref{chthm}.
$\Delta = (t^2+1)^2$  

\vskip 5mm

\centerline{$\star$ Compound rational homology spheres $\star$}
\vskip 5mm

These manifolds are primarily the union of two Seifert-fibred manifolds that
fibre over a disc, with at most $3$ singular fibres.  We compute the
$\vec \mu$-invariant via the Kaplan algorithm (see Theorem 5.7.14 of
\cite{GS}). We do not compute the $\vec d$-invariant as at present there
is no simple way to compute $\vec d$ for these manifolds.   To apply the Kaplan
algorithm we need an integral surgery diagram to start with.  There is
a rather simple way to construct surgery presentations for these manifolds, 
see \S \ref{hyperbolicity}.
 
\item \label{firsteg} SFS $\left[D: \frac{1}{2}, \frac{1}{2}\right]$ U/m 
      SFS $\left[D: \frac{1}{2}, \frac{1}{3}\right]$ 
$m = \begin{pmatrix}  -3 & 4 \\ -2 & 3\end{pmatrix}$  $H_1 = \Zed_2^2$ 
We follow the techniques of \S \ref{hyperbolicity} to construct a surgery
presentation for this manifold. If we label the components of the surgery link 
left-to-right we get
$$\includegraphics[height=2cm]{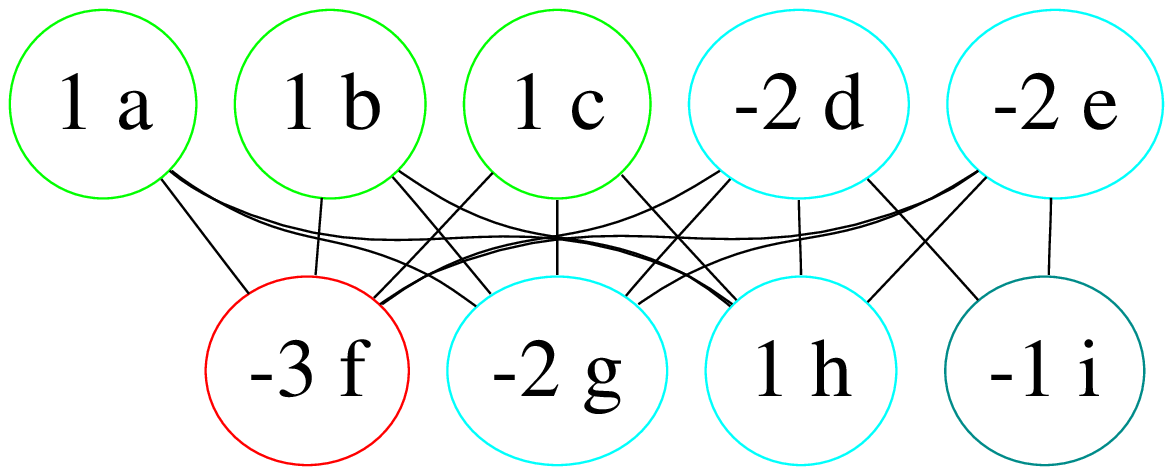}$$
as the graph of framing/linking numbers.  The four $\spin$-structures on this manifold correspond to the
characteristic links, which are given by
$$( \{f, g, h, i\}, \{d, e, f, g, h, i\}, \{a, b, c, d, f, h\}, \{a, b, c, e, f, h\}).$$
We apply Kaplan's algorithm to construct surgery presentations of the $\spin$
$4$-manifolds bounding each of these $\spin$ $3$-manifolds, which will allow us
to compute $\vec \mu$. In the order the characteristic links are listed for the above case, 
$\vec \mu = (\frac{1}{2}, 1, 0, 0)$. So this fails the Rochlin vector test (Corollary \ref{vectest}).

\item \label{2ndeg} SFS $\left[D: \frac{1}{2}, \frac{1}{2}\right]$ U/m 
      SFS $\left[D: \frac{1}{2}, \frac{2}{3}\right]$ 
$m = \begin{pmatrix}  -3 & 4 \\ -2 & 3\end{pmatrix}$  $H_1 = \Zed_2^2$. 

Characteristic links 
$(\{f,h,j\},$ $\{d,e,f,h,j\}, \{a,b,c,d,f,g,h\}, \{a,b,c,e,f,g,h\} )$, 
$\vec \mu = (\frac{3}{8}, -\frac{3}{8}, \frac{7}{8}, \frac{7}{8} )$.
Surgery presentation framing/linking matrix for item \ref{2ndeg} 
$\includegraphics[height=1.5cm]{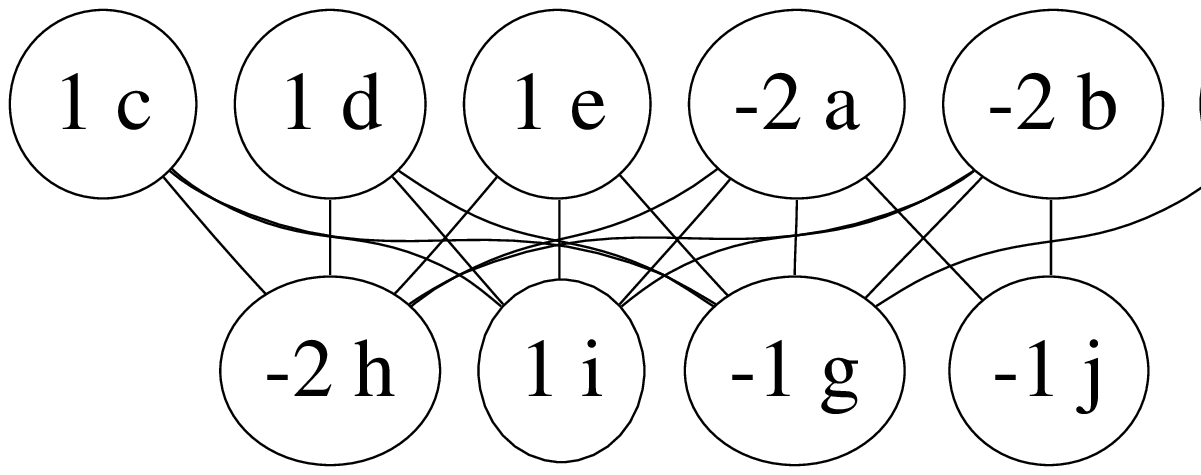}$

\item SFS $\left[D: \frac{1}{2}, \frac{1}{2}\right]$ U/m 
      SFS $\left[D: \frac{1}{2}, \frac{1}{3}\right]$ 
$m = \begin{pmatrix}  1 & 4 \\ 1 & 3\end{pmatrix}$  $H_1 = \Zed_2^2$.

Characteristic links 
$(\{a,b,d,e\}, \{a,c,d,e\}, \{e,f,g,h\}, \{b,c,e,f,g,h\} )$,
$\vec \mu = (-\frac{1}{8}, \frac{3}{8}, 1, 1)$.

Surgery diagram: $\includegraphics[height=1.5cm]{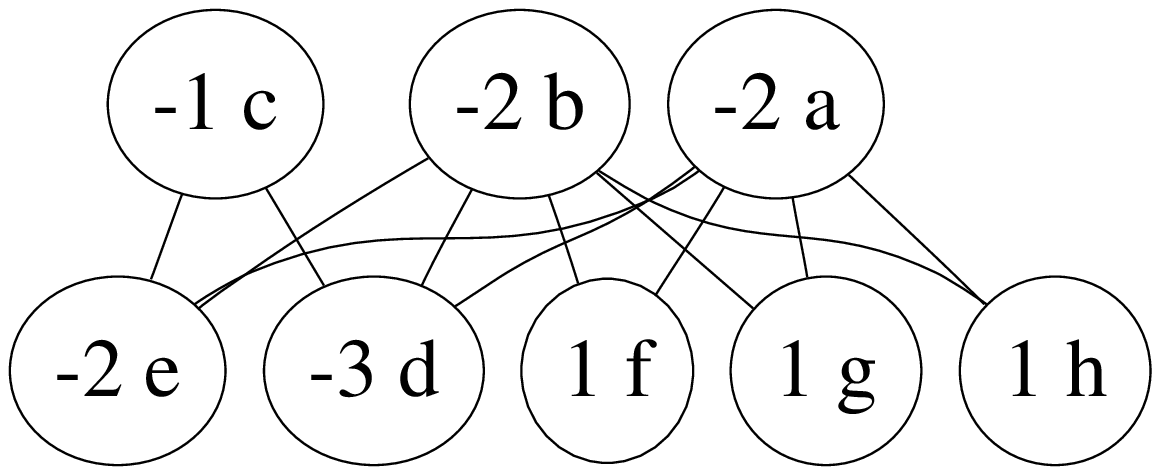}$

\item SFS $\left[D: \frac{1}{2}, \frac{1}{2}\right]$ U/m 
      SFS $\left[D: \frac{1}{2}, \frac{3}{5}\right]$ 
$m = \begin{pmatrix}  -3 & 4 \\ -2 & 3\end{pmatrix}$ $H_1 = \Zed_2^2$.

Characteristic links $( \{g,h,i,j\}, \{e,f,g,h,i,j\}, \{b,c,d,e,g,i\}, \{b,c,d,f,g,i\})$,
$\vec \mu = (\frac{1}{2},1,\frac{3}{8}, \frac{3}{2})$.
Surgery diagram: $\includegraphics[height=1.5cm]{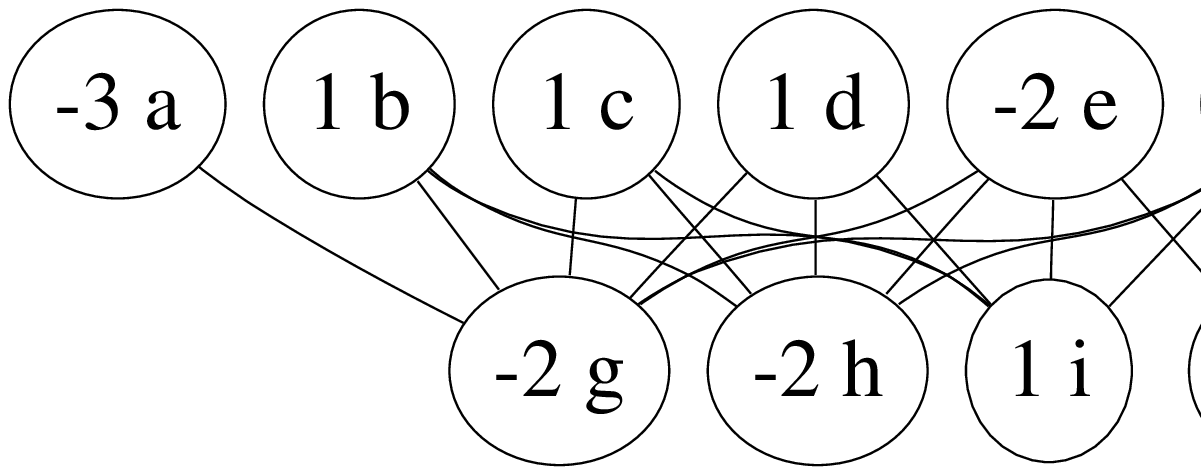}$

\item SFS $\left[D: \frac{1}{2}, \frac{1}{2}\right]$ U/m 
      SFS $\left[D: \frac{1}{2}, \frac{1}{6}\right]$ 
$m = \begin{pmatrix}  1 & 1 \\ 1 & 0\end{pmatrix}$ $H_1 = \Zed_4^2$.

Characteristic links 
$( \{c,d\}, \{c,e\}, \{a,b,c,d\}, \{a,b,c,e\} )$, $\vec \mu = (0,0,1, 1)$.
Surgery diagram: $\includegraphics[height=1.5cm]{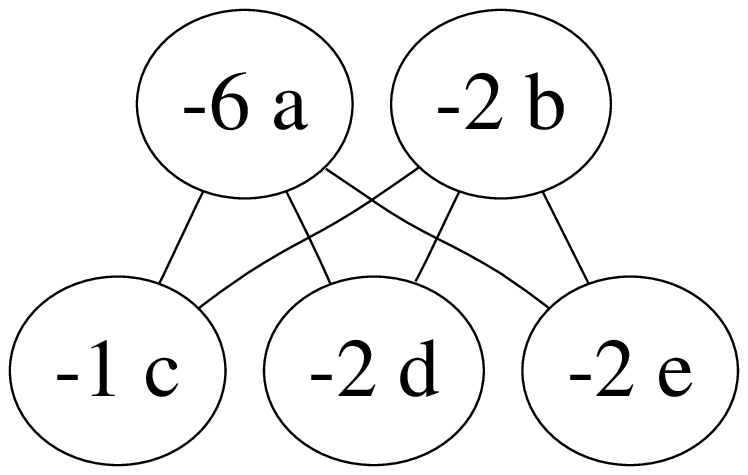}$

\item SFS $\left[D: \frac{1}{2}, \frac{1}{2}\right]$ U/m 
      SFS $\left[D: \frac{1}{3}, \frac{1}{3}\right]$ 
$m = \begin{pmatrix}  2 & 1 \\ 1 & 0\end{pmatrix}$  $H_1 = \Zed_6^2$. 

Characteristic links 
$(\{d,e,f\}, \{c,e,f\}, \{a,b,d,e,f\}, \{a,b,c,e,f\} )$, $\vec \mu = (0,0,\frac{3}{4}, \frac{3}{4})$.

Surgery diagram: $\includegraphics[height=1.5cm]{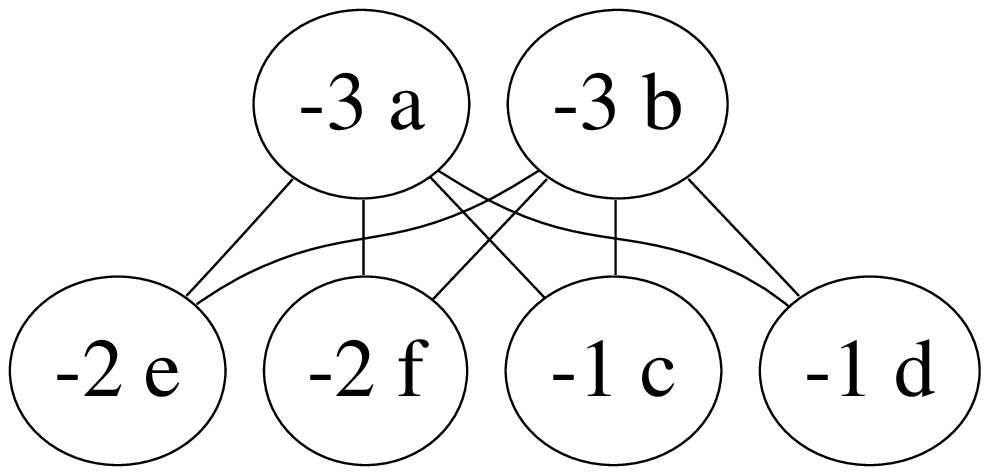}$

\item SFS $\left[D: \frac{1}{2}, \frac{1}{2}\right]$ U/m 
      SFS $\left[D: \frac{1}{3}, \frac{3}{4}\right]$ 
$m = \begin{pmatrix}  -1 & 2 \\ 0 & 1\end{pmatrix}$  $H_1 = \Zed_2^2$. 

Characteristic links $(\{b,c,f,g\}, \{b,c,d,e,f,g\}, \{a,c,d,g,h\}, \{a,c,e,g,h\})$,
$\vec \mu = (\frac{1}{2}, 1, 0, 0)$.

Surgery diagram: $\includegraphics[height=1.5cm]{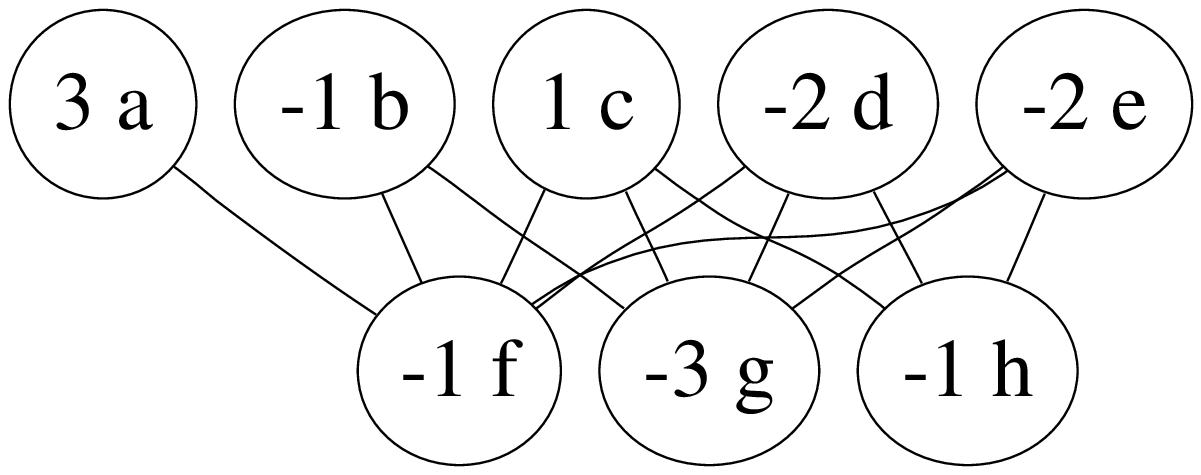}$

\item SFS $\left[D: \frac{1}{2}, \frac{1}{2}\right]$ U/m 
      SFS $\left[D: \frac{1}{3}, \frac{5}{7}\right]$ 
$m = \begin{pmatrix}  0 & 1 \\ 1 & 0\end{pmatrix}$   $H_1 = \Zed_2^2$.

Characteristic links 
$(\{a,b,c\}, \{d\}, \{e\}, \{a,b,c,d,e\})$, $\vec \mu = (\frac{1}{2}, 0, 0, 1)$.

Surgery diagram: $\includegraphics[height=1.5cm]{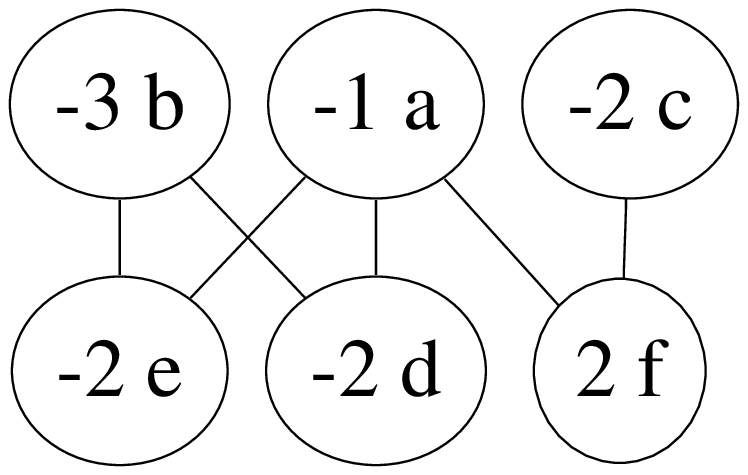}$

\item SFS $\left[D: \frac{1}{2}, \frac{1}{2}\right]$ U/m 
      SFS $\left[D: \frac{2}{3}, \frac{3}{8}\right]$ 
$m = \begin{pmatrix}  -1 & 2 \\ 0 & 1\end{pmatrix}$ $H_1 = \Zed_2^2$.

Characteristic links $(\{a,b,c,d\}, \{a,b,c,d,e,f\}, \{b,d,e,g,i\}, \{b,d,f,g,i\})$,
$\vec \mu = (-\frac{3}{4}, -\frac{1}{4}, 0, 0)$.
Surgery diagram: $\includegraphics[height=1.5cm]{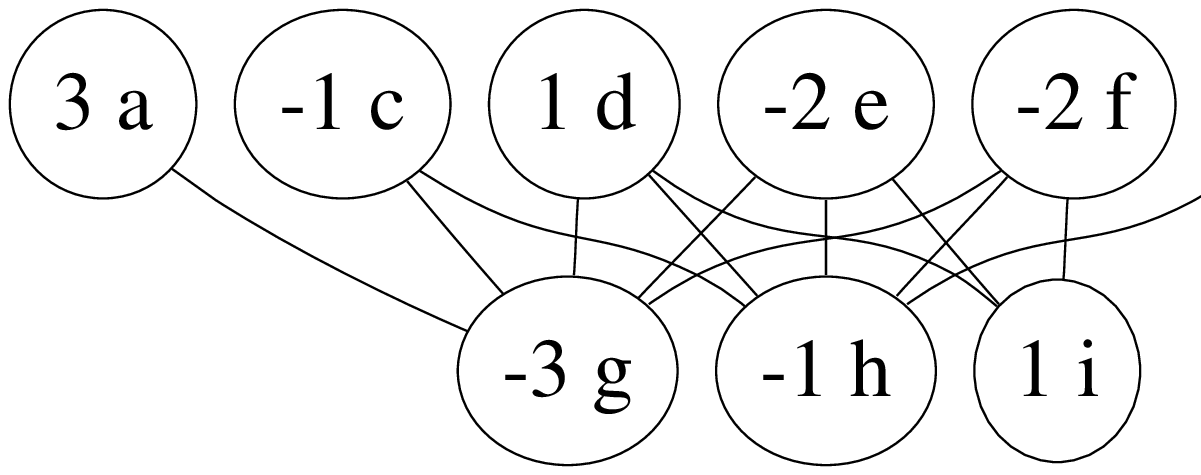}$

\item SFS $\left[D: \frac{1}{2}, \frac{2}{3}\right]$ U/m 
      SFS $\left[D: \frac{1}{2}, \frac{1}{2}, \frac{1}{2}\right]$ 
$m = \begin{pmatrix}  1 & 1 \\ 1 & 2\end{pmatrix}$ $H_1 = \Zed_2^2$.

Characteristic links $(\{a,b,g,h,i\}, \{a,b,d,e,g,h,i\}, \{a,b,d,f,g,h,i\}, \{a,b,e,f,g,h,i\})$,
$\vec \mu = (-\frac{1}{2}, 1, 1, 0)$.
Surgery diagram: $\includegraphics[height=1.5cm]{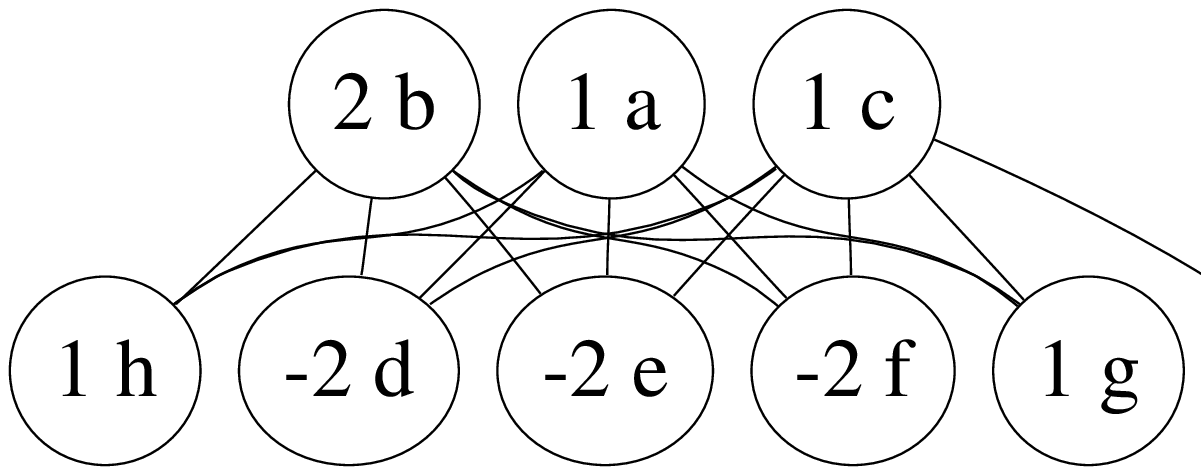}$
\vskip 5mm

\item SFS $\left[D: \frac{1}{2}, \frac{1}{2}\right]$ U/m 
      SFS $\left[D: \frac{2}{3}, \frac{4}{11}\right]$ 
$m = \begin{pmatrix}  0 & 1 \\ 1 & 0\end{pmatrix}$ $H_1 = \Zed_2^2$.
Issa-McCoy obstruction \cite{Issa}.
Characteristic links 

$(\{a,b\}, \{c\}, \{d\}, \{a,b,c,d\})$.

Surgery diagram: $\includegraphics[height=1.5cm]{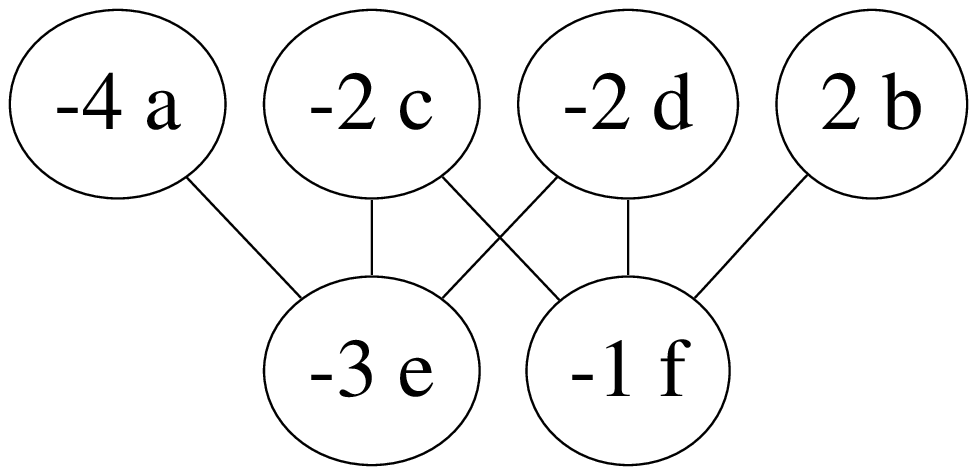}$

\item \label{nonemb_num} SFS $\left[D: \frac{1}{2}, \frac{4}{7}\right]$ U/m 
      SFS $\left[D: \frac{1}{3}, \frac{2}{3}\right]$ 
$m = \begin{pmatrix}  0 & 1 \\ 1 & 0 \end{pmatrix}$ $H_1 = \Zed_3^2$.
Issa-McCoy obstruction \cite{Issa}.
Characteristic link $\{b\}$.

Surgery diagram: $\includegraphics[height=1.5cm]{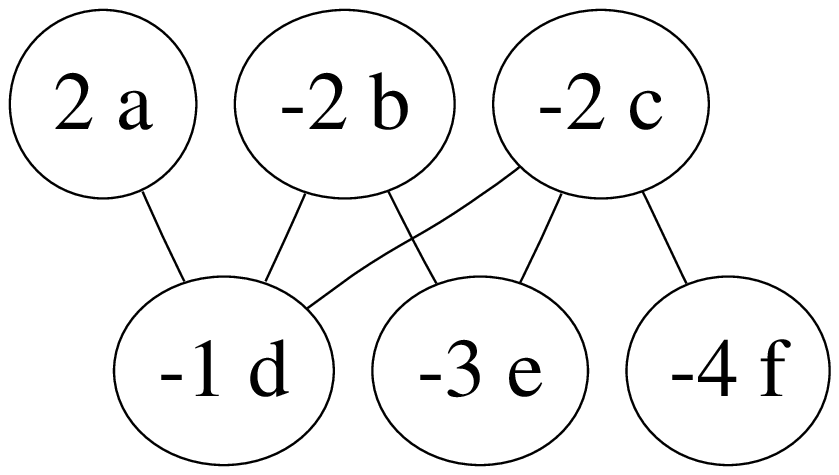}$

\end{enumerate}

\section{Manifolds for which embeddability is not known}\label{unknown_embeddable_man}
\vskip 5mm
\centerline{$\star$ $SL_2\Real$-manifolds with finite $H_1$ $\star$}
\vskip 5mm

\begin{enumerate}

\item SFS $\left[\RProj^2/n2: \frac{3}{5}, \frac{3}{5}\right]$  $H_1 = \Zed_{10}^2$. 
$\vec d$ not computed.
Characteristic links 

$(\{a,c,d\}, \{b,c,d\}, \{e,f\}, \{a,b,e,f\})$,
$\vec \mu = (-\frac{1}{2}, 0,0,0 )$.
Surgery diagram $\includegraphics[height=1.5cm]{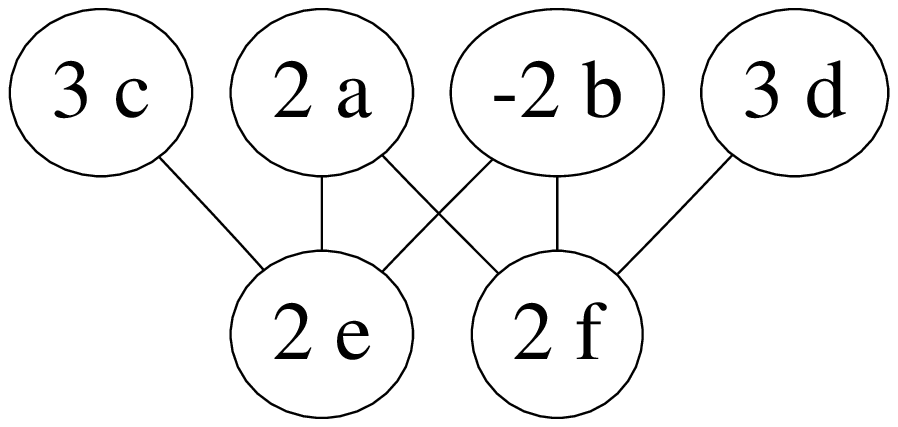}$
\vskip 5mm

\centerline{$\star$ Hyperbolic manifolds $\star$}
\vskip 5mm

These manifolds are uniquely identified in Burton's census \cite{BBurton} by their
volumes. The Rochlin invariant is given from a surgery presentation via
Theorem 2.13 \cite{Sav}.  See \S \ref{hyperbolicity} for notes on how surgery
presentations are found.  The Rochlin invariant is computed as described in \S \ref{Known results}. 
A brief description of the calculation is given below. See Theorem \ref{spintests}.

\item Hyp 1.96273766 $H_1 = \Zed_7^2$. Initial surgery
presentation on $\KR{R}{8^3_8}$ found via SnapPea. The first reduction eliminates the
unknotted component with framing number $-1/3$ via a Rolfsen twist on that
component. A second move creates integral surgery via slam-dunk move on component
with framing number $\frac{14}{3}$. 

{
\psfrag{k38}[tl][tl][1][0]{$8^3_8$}
\psfrag{s1}[tl][tl][0.7][0]{$\frac{5}{3}$}
\psfrag{s2}[tl][tl][0.7][0]{$-\frac{1}{3}$}
\psfrag{s3}[tl][tl][0.7][0]{$\frac{2}{1}$}
$\includegraphics[width=4cm]{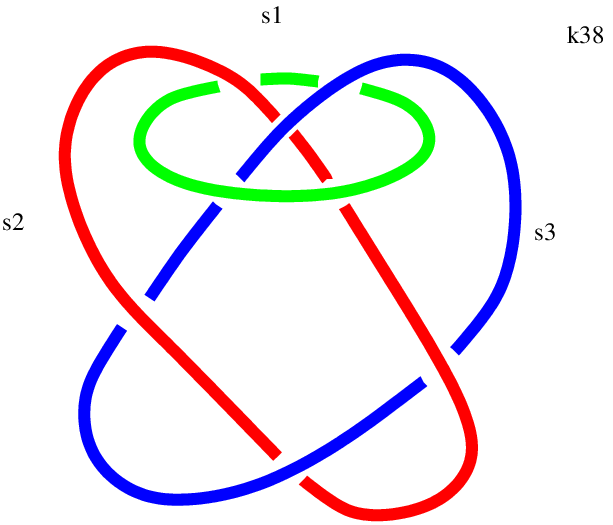}$
\hskip 1cm
}
{
\psfrag{14/3}[tl][tl][0.8][0]{$\frac{14}{3}$}
\psfrag{14}[tl][tl][1][0]{$14$}
$\includegraphics[width=4cm]{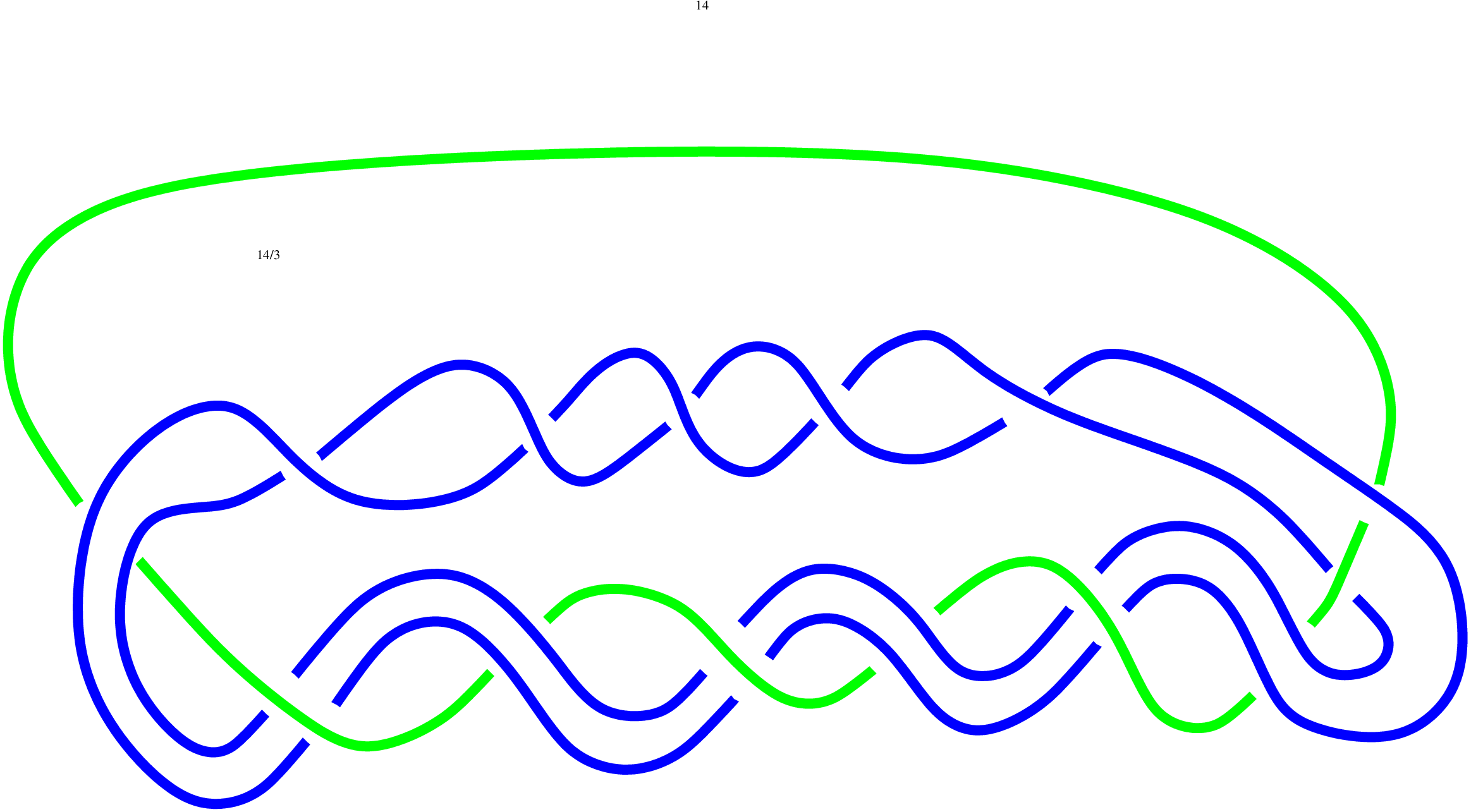}$ 
\hskip 1cm
}
{
\psfrag{5}[tl][tl][1][0]{$5$}
\psfrag{14}[tl][tl][1][0]{$14$}
\psfrag{3}[tl][tl][1][0]{$3$}
$\includegraphics[width=4cm]{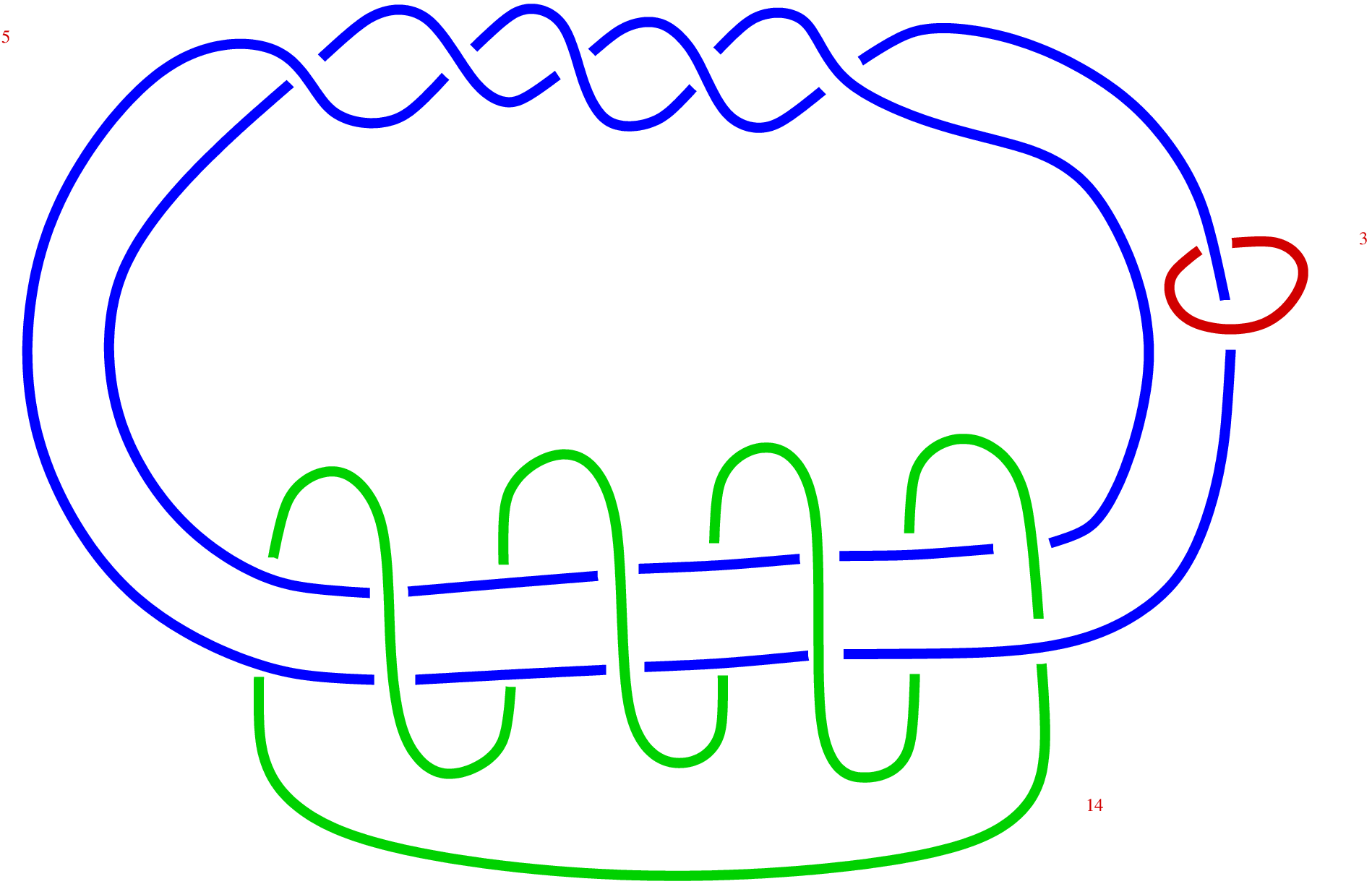}$
}
To which we apply the Kaplan algorithm to get the presentation:

{
\psfrag{5}[tl][tl][1][0]{$4$}
\psfrag{14}[tl][tl][1][0]{$14$}
\psfrag{2}[tl][tl][1][0]{$-2,-2$}
$\includegraphics[width=4cm]{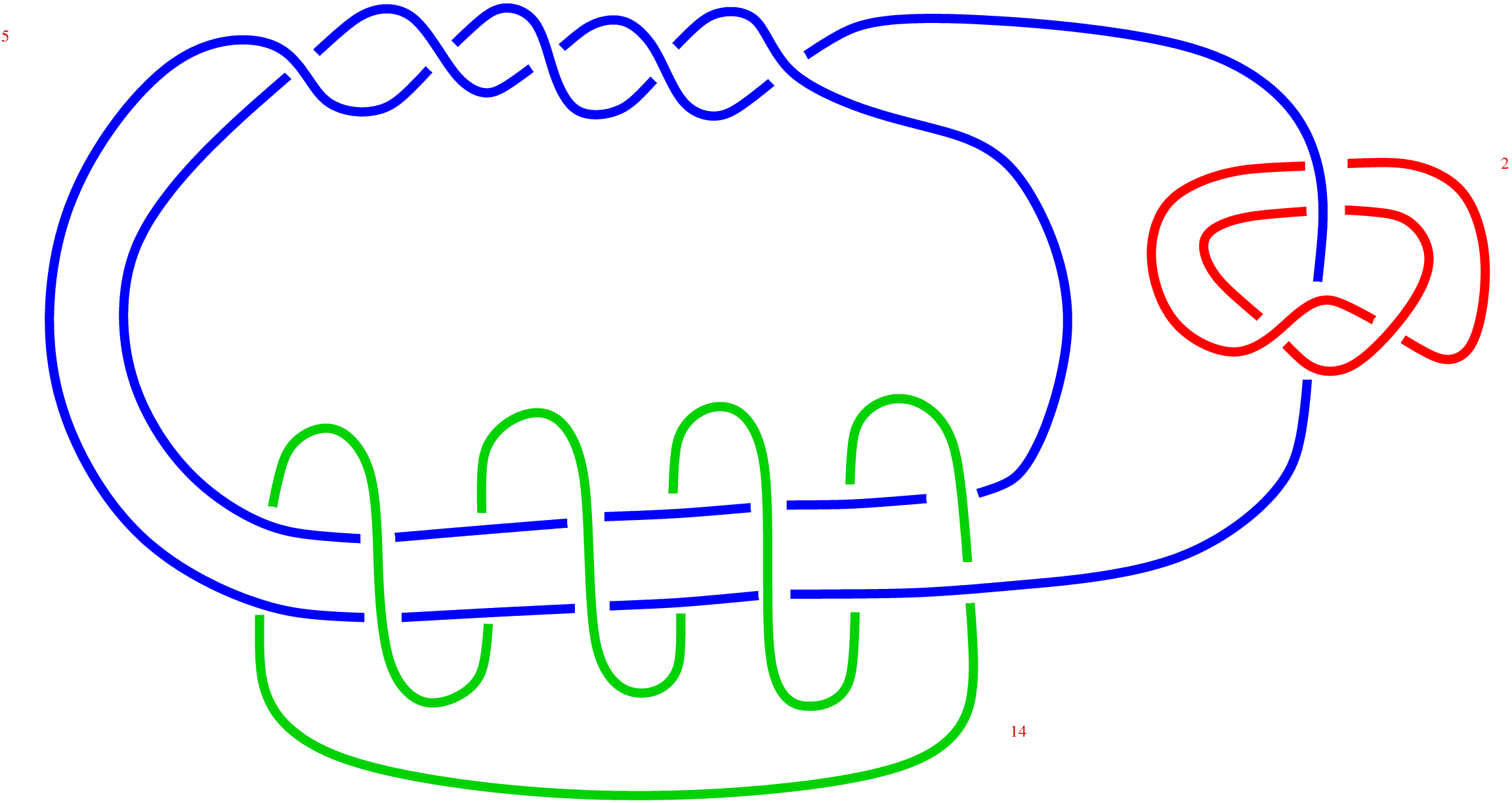}$
}
\hskip 2cm
$\includegraphics[width=4cm]{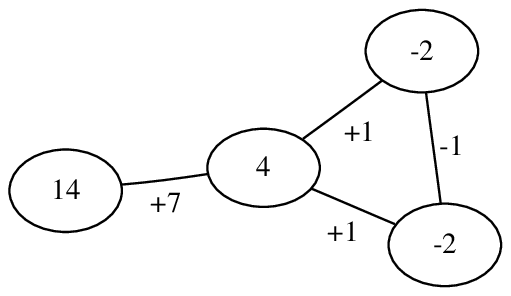}$

The graph consists of the framing/linking numbers. The characteristic polynomial
of the intersection product is $t^4-14t^3-16t^2+49$, thus the signature is
zero and $\mu = 0$. 

\item \label{lastgeometric} Hyp 2.22671790, Homology sphere. $\mu = 0$. 
$+\frac{1}{2}$-surgery on $\KR{R}{5_2}$ found
via Snappea. $\mu$ computed via Theorem 2.10 in \cite{Sav}. 

$\includegraphics[width=3cm]{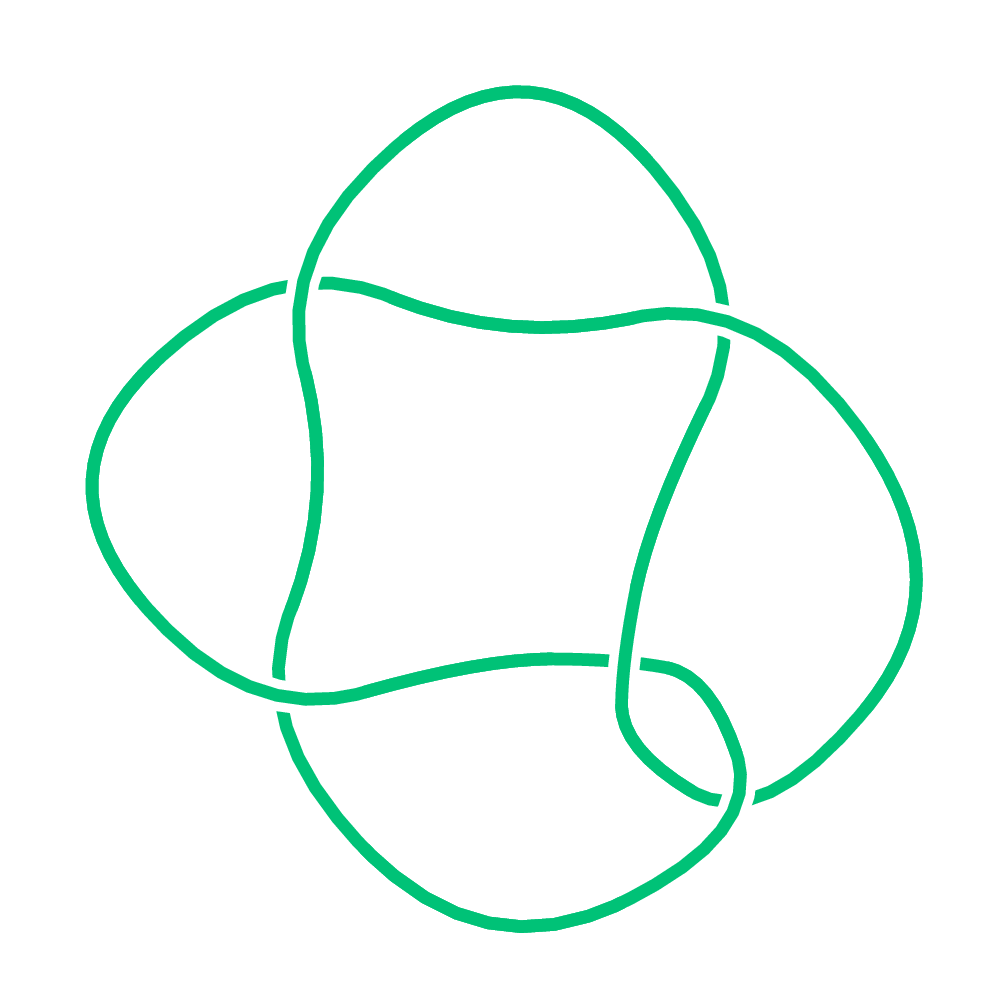}$
\vskip 5mm

\centerline{$\star$ Graph manifold with single non-separating torus in JSJ $\star$}
\vskip 5mm

These manifolds are all of the form $SFS \left[ A: \frac{\alpha}{\beta} \right] / 
\begin{pmatrix} a & b \\ c & d \end{pmatrix}$ where $ad-bc = -1$. These manifolds have
$\rank(H_1) =1$ if and only if the polynomial $\beta t^2 + ( (d-a)\beta - b\alpha)t + \beta$
does not have $1$ as a root, moreover this is the Alexander polynomial in this case. Thus
the three manifolds below all have Alexander polynomial $\Delta = 2t^2-5t+2$,
which satisfies Kawauchi's Theorem \ref{KawCond}.  Unfortunately, 
$2t^2 - 5t +2 = (2t-1)(t-2)$ so all signature invariants are zero for these manifolds.

\item SFS $\left[A: \frac{1}{2}\right]$ / $\begin{pmatrix}  0 & 1 \\ 1 & -2\end{pmatrix}$ $H_1 = \Zed$. 
\ttc
\item SFS $\left[A: \frac{1}{2}\right]$ / $\begin{pmatrix}  2 & 5 \\ 1 & 2\end{pmatrix}$ $H_1 = \Zed$
\item SFS $\left[A: \frac{1}{2}\right]$ / $\begin{pmatrix}  -1 & 3 \\ 1 & -2\end{pmatrix}$ $H_1 = \Zed$
\item SFS $\left[A: \frac{1}{2}\right]$ / $\begin{pmatrix}  1 & -3 \\ -1 & 2\end{pmatrix}$ $H_1 = \Zed \oplus \Zed_3^2$. $\Delta = 2t^2 + 5t + 2 = (2t+1)(t+2)$ also satisfies Theorem \ref{KawCond} and
has trivial signature invariants.

\item SFS $\left[A: \frac{1}{3}\right]$ / $\begin{pmatrix}  -1 & 3 \\ 1 & -2\end{pmatrix}$ $H_1 = \Zed^2$. 
No tests have been performed for this manifold.

\vskip 5mm

\centerline{$\star$ Compound homology spheres $\star$}
\vskip 5mm

Their splicing decomposition is listed and $\mubar$ is computed using
splicing additivity (Proposition 2.16 in \cite{Sav}).

\item SFS $\left[D: \frac{1}{2}, \frac{2}{5}\right]$ U/m 
      SFS $\left[D: \frac{1}{2}, \frac{3}{5}\right]$ 
$m = \begin{pmatrix}  0 & 1 \\ 1 & 0\end{pmatrix}$ Homology sphere.

$\Sigma(2,5,9) \splice -\Sigma(2,5,9)$ $\mubar=0$

\item SFS $\left[D: \frac{1}{2}, \frac{2}{3}\right]$ U/m 
      SFS $\left[D: \frac{1}{2}, \frac{2}{5}\right]$ 
$m = \begin{pmatrix}  -2 & 3 \\ -1 & 2 \end{pmatrix}$ Homology sphere.

$\Sigma(2,5,9) \splice \Sigma(2,3,7)$ $\mubar=0$

\item SFS $\left[D: \frac{1}{2}, \frac{2}{3}\right]$ U/m 
      SFS $\left[D: \frac{1}{4}, \frac{3}{5}\right]$ 
$m = \begin{pmatrix}  0 & 1 \\ 1 & 0 \end{pmatrix}$  Homology sphere.

$\Sigma(2,3,7) \splice -\Sigma(4,5,7)$ $\mubar=0$

\item SFS $\left[D: \frac{1}{2}, \frac{1}{3}\right]$ U/m 
      SFS $\left[A: \frac{1}{2}\right]$ U/n 
      SFS $\left[D: \frac{1}{2}, \frac{2}{3}\right]$ 
$m = \begin{pmatrix}  -1 & 1 \\ 1 & 0 \end{pmatrix}, 
 n = \begin{pmatrix}  0 & 1 \\ 1 & 1 \end{pmatrix}$ Homology sphere.
$\Sigma(2,3,11) \splice S^3(L) \splice -\Sigma(2,3,11)$ where this indicates 
splicing over the link $L$ in $S^3$ which is the union of two regular 
fibres in the `$(2,1)$-fibring' of $S^3$. 
$\mubar = 0$ 
\vskip 5mm

\centerline{$\star$ Compound rational homology spheres $\star$}
\vskip 5mm

See \S \ref{nonembeddable_man} for details on how the Rochlin 
vector $\vec \mu$ is computed for these manifolds.

\item SFS $\left[D: \frac{1}{2}, \frac{1}{2}\right]$ U/m 
      SFS $\left[D: \frac{1}{2}, \frac{2}{5}\right]$ 
$m = \begin{pmatrix}  -3 & 4 \\ -2 & 3\end{pmatrix}$  $H_1 = \Zed_2^2$.

Characteristic links $( \{i,j\}, \{e,f,i,j\}, \{b,c,d,e,i,h\}, \{b,c,d,f,g,i\})$,
$\vec \mu = ( 0, \frac{1}{2}, 0, 0 )$.  This manifold can be thought of as the $(5,2)$-torus knot complement union the orientable $S^1$-bundle over a M\"obius band.  There is a natural embedding of such a manifold into $S^4$, since the $(5,2)$-torus knot bounds a Klein bottle, and the orientable $I$-bundle over the Klein bottle is the $S^1$-bundle over the M\"obius band. The above gluing map does not produce this manifold. 

Surgery diagram: $\includegraphics[height=1.5cm]{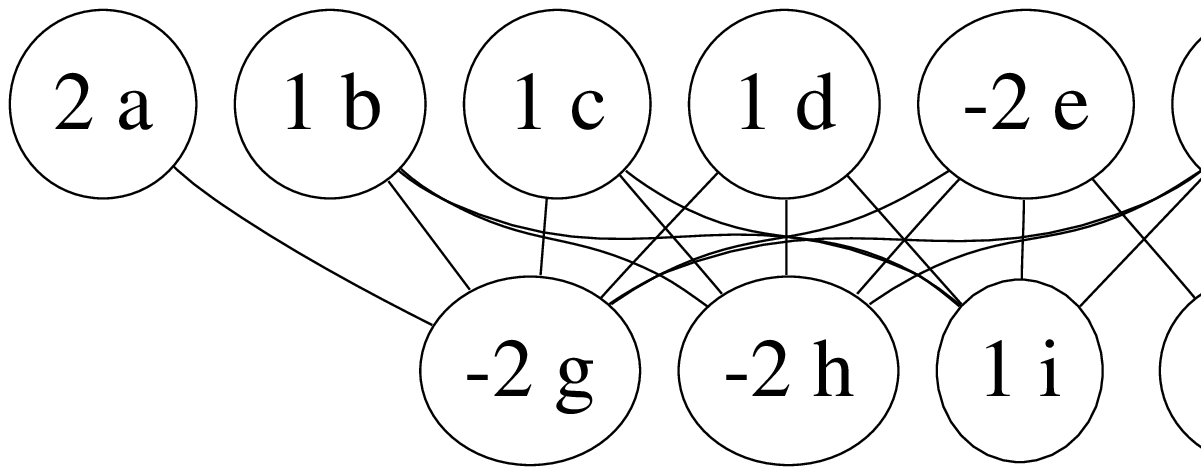}$

\item SFS $\left[D: \frac{1}{2}, \frac{1}{2}\right]$ U/m 
      SFS $\left[D: \frac{2}{3}, \frac{1}{4}\right]$ 
$m = \begin{pmatrix}  -1 & 2 \\ 0 & 1\end{pmatrix}$  $H_1 = \Zed_2^2$.

Characteristic links $(\{a,b,c\}, \{a,b,c,d,e\}, \{a,c,d,g,h\}, \{a,c,e,g,h\})$,
$\vec \mu = (-\frac{1}{2},0,0,0)$.

Surgery diagram: $\includegraphics[height=1.5cm]{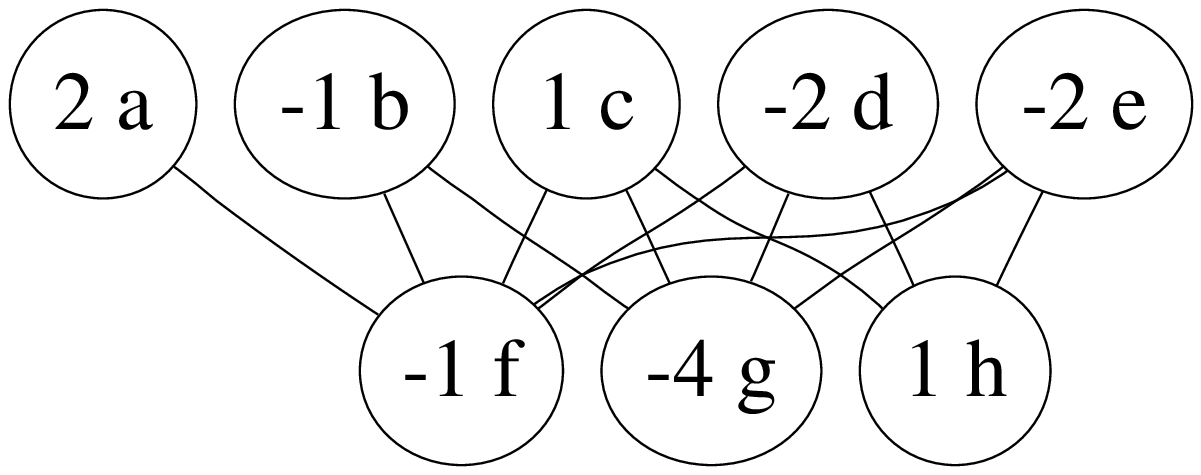}$

\item SFS $\left[D: \frac{1}{2}, \frac{1}{2}\right]$ U/m 
      SFS $\left[D: \frac{1}{3}, \frac{2}{3}\right]$ 
$m = \begin{pmatrix}  -2 & 5 \\ -1 & 3\end{pmatrix}$ $H_1 = \Zed_6^2$.

Characteristic links 
$(\{b,c,d,e,f,g,h\}, \{g,h,i\}, \{g,h,j\}, \{b,c,d,e,f,g,h,i,j\})$, $\vec \mu = (0, 0, 0,\frac{1}{2})$.

Surgery diagram: $\includegraphics[height=1.5cm]{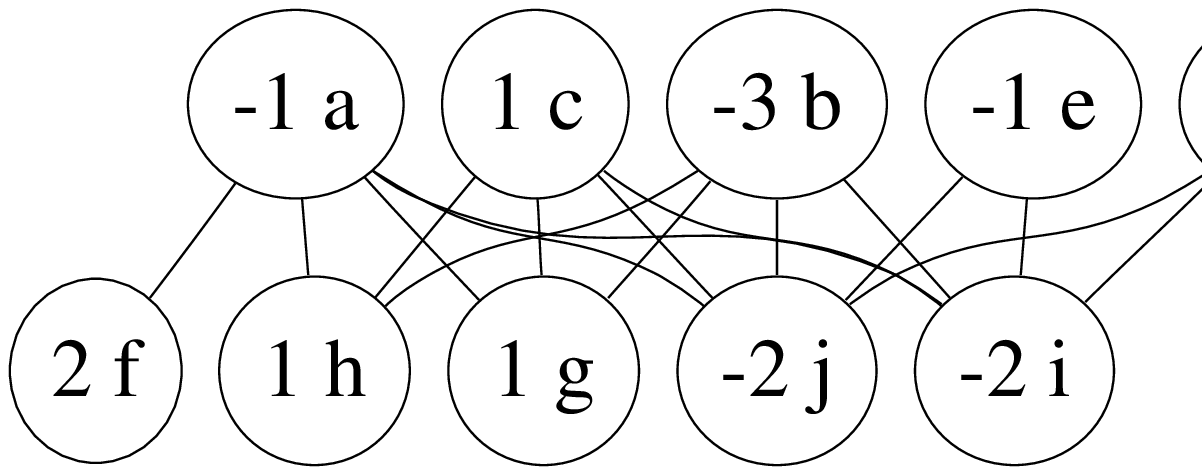}$

\item SFS $\left[D: \frac{1}{2}, \frac{1}{2}\right]$ U/m 
      SFS $\left[D: \frac{1}{3}, \frac{5}{8}\right]$ 
$m = \begin{pmatrix}  -1 & 2 \\ 0 & 1\end{pmatrix}$ $H_1 = \Zed_2^2$.

Characteristic links $(\{b,c,f,g,h\}, \{b,c,d,e,f,g,h\}, \{a,c,d,h,i\}, \{a,c,e,h,i\})$,
$\vec \mu = (0, \frac{1}{2}, 0, 0)$.

Surgery diagram: $\includegraphics[height=1.5cm]{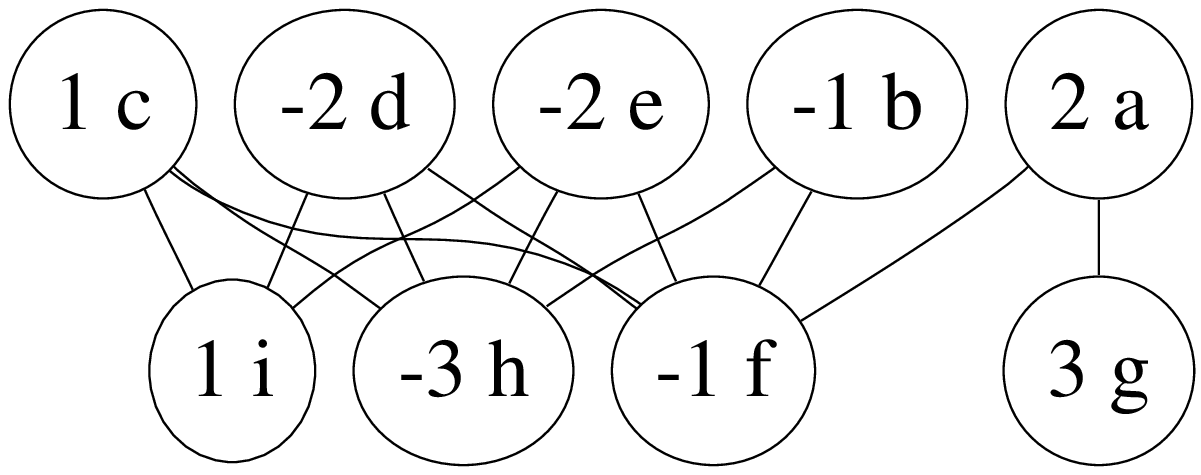}$

\item SFS $\left[D: \frac{1}{2}, \frac{1}{2}\right]$ U/m 
      SFS $\left[D: \frac{2}{3}, \frac{2}{5}\right]$ 
$m = \begin{pmatrix}  -2 & 3 \\ -1 & 2\end{pmatrix}$ $H_1 = \Zed_2^2$.

Characteristic links $(\{b,c,d,e\}, \{b,c,d,f\}, \{a,i,j\}, \{a,e,f,i,j\})$,
$\vec \mu = (0, 0, 0, \frac{1}{2})$.

Surgery diagram: $\includegraphics[height=1.5cm]{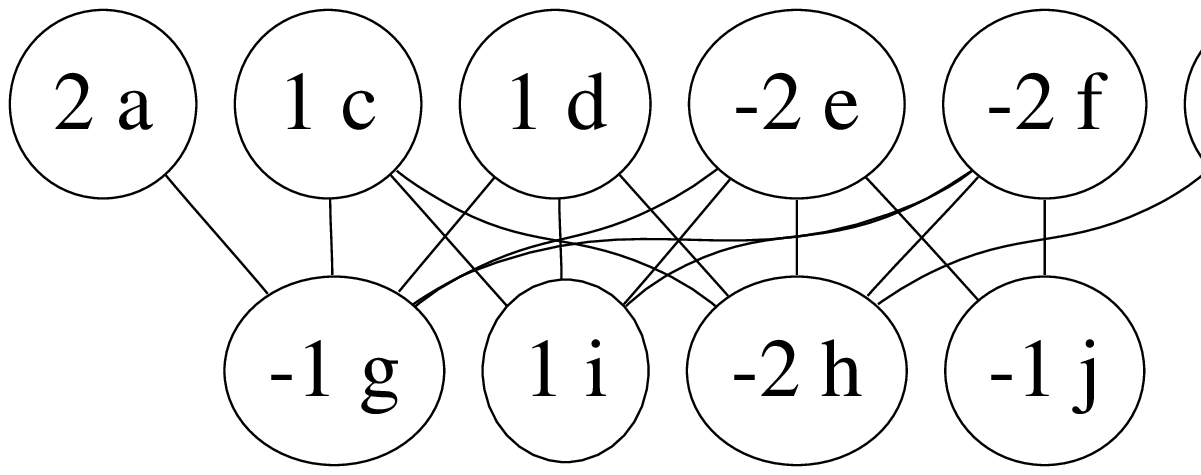}$

\item SFS $\left[D: \frac{1}{2}, \frac{1}{2}\right]$ U/m 
      SFS $\left[D: \frac{2}{5}, \frac{4}{7}\right]$ 
$m = \begin{pmatrix}  0 & 1 \\ 1 & 0\end{pmatrix}$ $H_1 = \Zed_2^2$.

Characteristic links 
$(\phi, \{a,b,c\}, \{a,b,d\}, \{c,d\})$, $\vec \mu = (\frac{1}{2}, 0,0,0)$.

Surgery diagram: $\includegraphics[height=1.5cm]{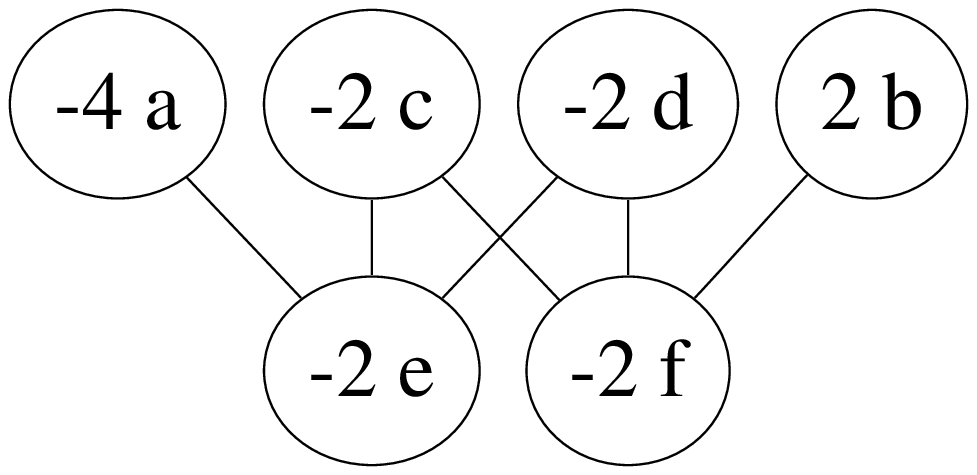}$

\item SFS $\left[D: \frac{1}{2}, \frac{1}{3}\right]$ U/m 
      SFS $\left[D: \frac{1}{3}, \frac{2}{3}\right]$ 
$m = \begin{pmatrix}  -2 & 3 \\ -1 & 2\end{pmatrix}$ $H_1 = \Zed_3^2$.

Characteristic link $\{b,c,e\}, \mu = 0$.

Surgery diagram: $\includegraphics[height=1.5cm]{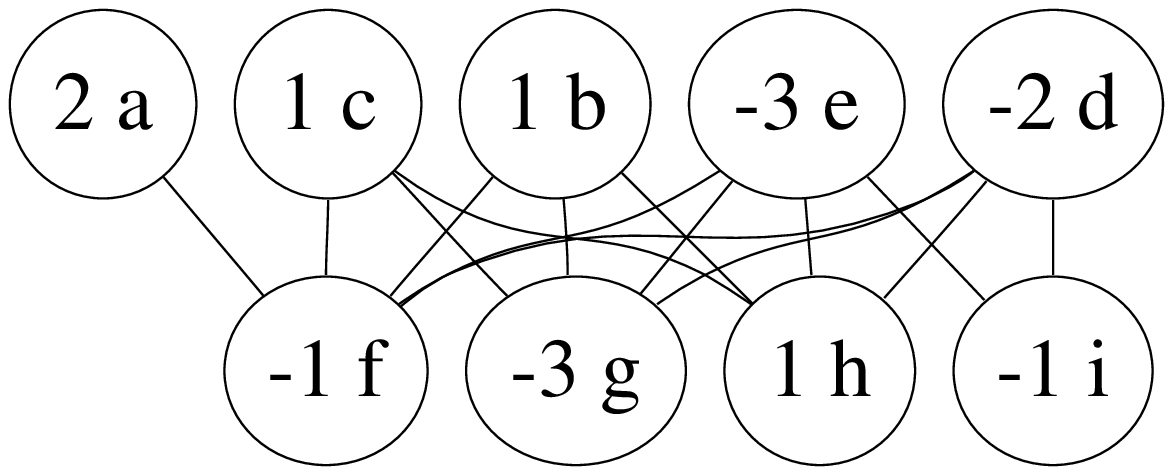}$

\item SFS $\left[D: \frac{1}{2}, \frac{1}{3}\right]$ U/m 
      SFS $\left[D: \frac{2}{5}, \frac{3}{5}\right]$ 
$m = \begin{pmatrix}  0 & 1 \\ 1 & 0\end{pmatrix}$ $H_1 = \Zed_5^2$.

Characteristic link $\{a,c,d\}$, $\mu = 0$.

Surgery diagram: $\includegraphics[height=1.5cm]{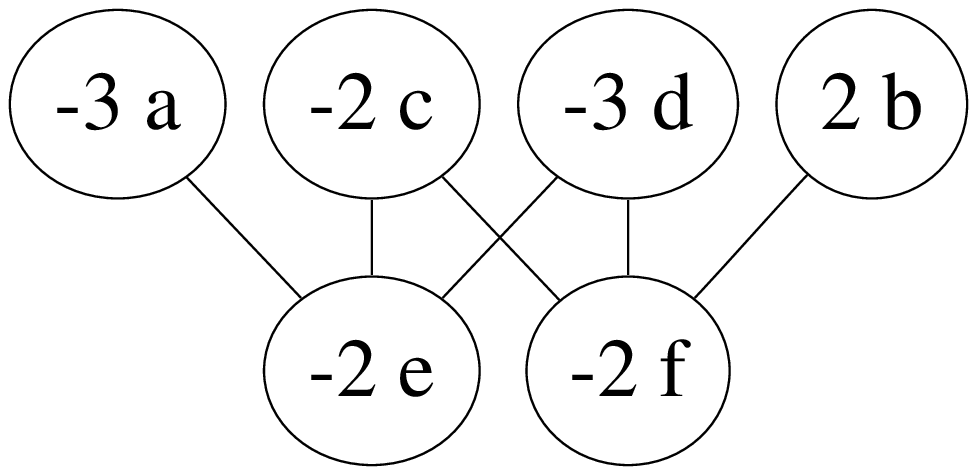}$

\item SFS $\left[D: \frac{1}{2}, \frac{2}{3}\right]$ U/m 
      SFS $\left[D: \frac{1}{3}, \frac{2}{3}\right]$ 
$m = \begin{pmatrix}  -2 & 3 \\ -1 & 2 \end{pmatrix}$ $H_1 = \Zed_3^2$.

Characteristic link $\{b,c,d,j\}$, $\mu = 0$.

Surgery diagram: $\includegraphics[height=1.5cm]{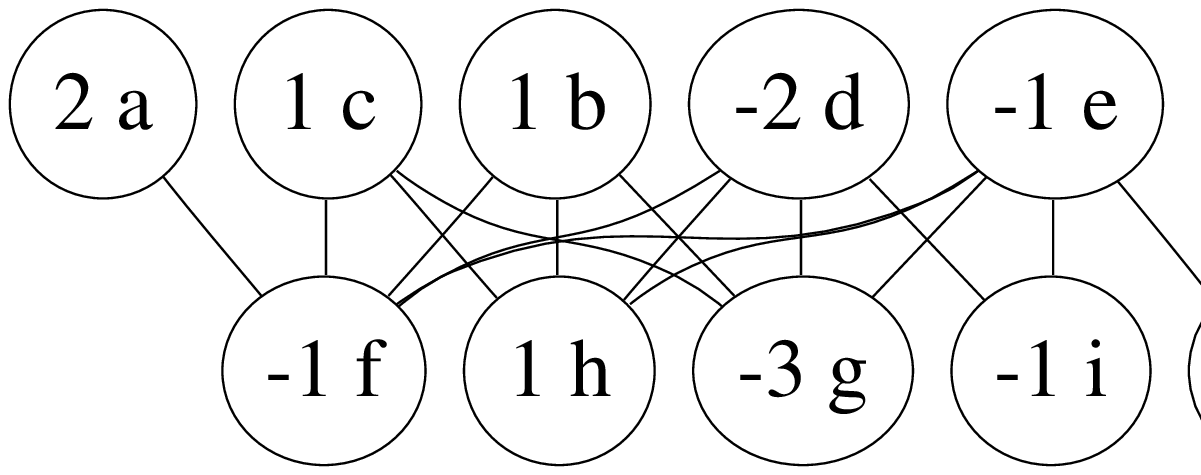}$

\item SFS $\left[D: \frac{1}{2}, \frac{2}{5}\right]$ U/m 
      SFS $\left[D: \frac{1}{3}, \frac{2}{3}\right]$ 
$m = \begin{pmatrix}  0 & 1 \\ 1 & 0\end{pmatrix}$ $H_1 = \Zed_3^2$.

Characteristic link $\{b\}$, $\mu = 0$.

Surgery diagram: $\includegraphics[height=1.5cm]{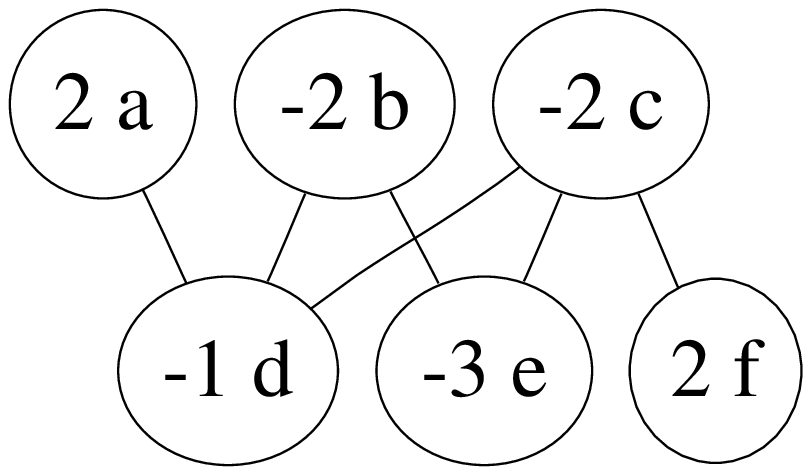}$

\item SFS $\left[D: \frac{1}{2}, \frac{3}{5}\right]$ U/m 
      SFS $\left[D: \frac{1}{3}, \frac{2}{3}\right]$ 
$m = \begin{pmatrix}  -1 & 2 \\ 0 & 1\end{pmatrix}$  $H_1 = \Zed_3^2$.

Characteristic link $\{a,b,d,f,g\}$, $\mu = 0$.

Surgery diagram: $\includegraphics[height=1.5cm]{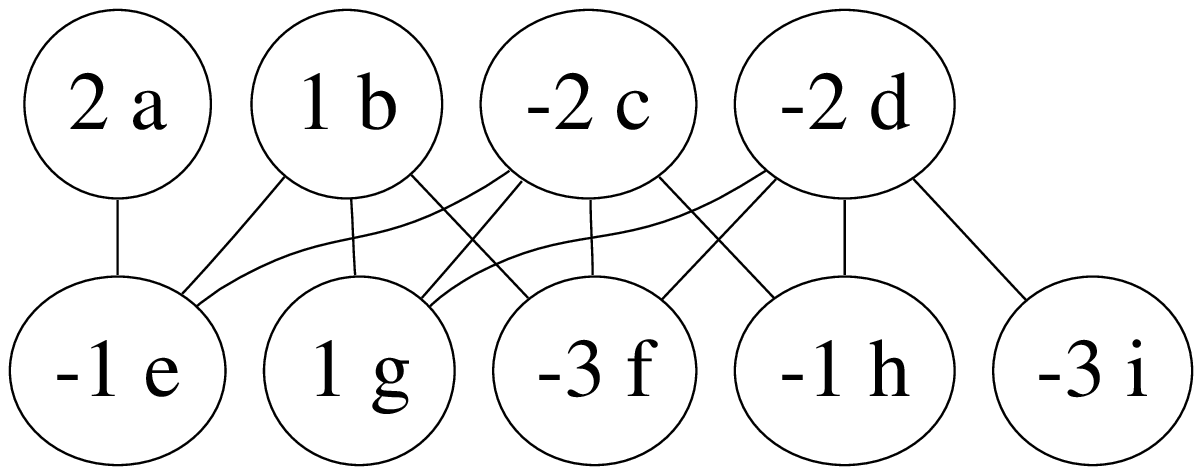}$

\item SFS $\left[D: \frac{1}{3}, \frac{1}{3}\right]$ U/m 
      SFS $\left[D: \frac{2}{3}, \frac{2}{3}\right]$ 
$m = \begin{pmatrix}  0 & 1 \\ 1 & 0\end{pmatrix}$  $H_1 = \Zed_3^2$.

Characteristic link $\{a,b,c,d\} $, $\mu = 0$.

Surgery diagram: $\includegraphics[height=1.5cm]{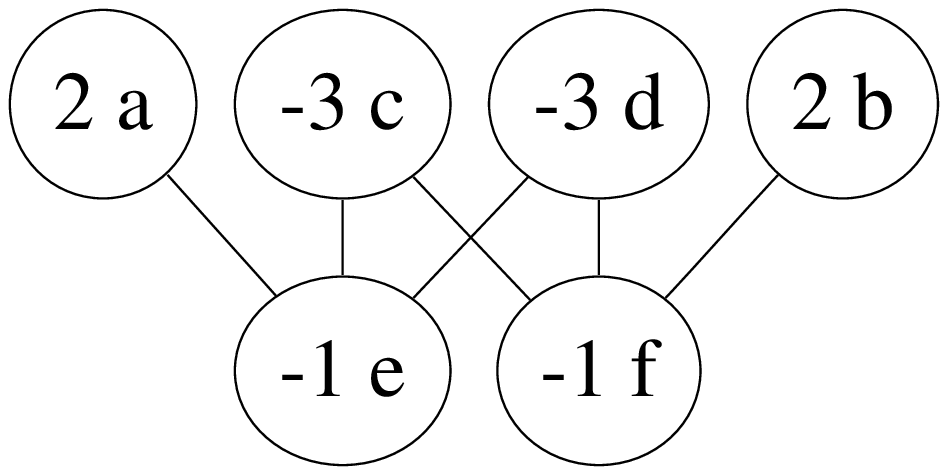}$

\item SFS $\left[D: \frac{1}{2}, \frac{1}{2}\right]$ U/m 
      SFS $\left[D: \frac{1}{2}, \frac{1}{2}, \frac{1}{3}\right]$ 
$m = \begin{pmatrix}  1 & 1 \\ 1 & 2\end{pmatrix}$ $H_1 = \Zed_4^2$.

Characteristic links $(\{a,b,d,e,g,h\}, \{a,c,d,e,g,h\}, \{a,b,d,f,g,h\}, \{a,c,d,f,g,h\} )$,
$\vec \mu = (0,0,0,0)$.

Surgery diagram: $\includegraphics[height=1.5cm]{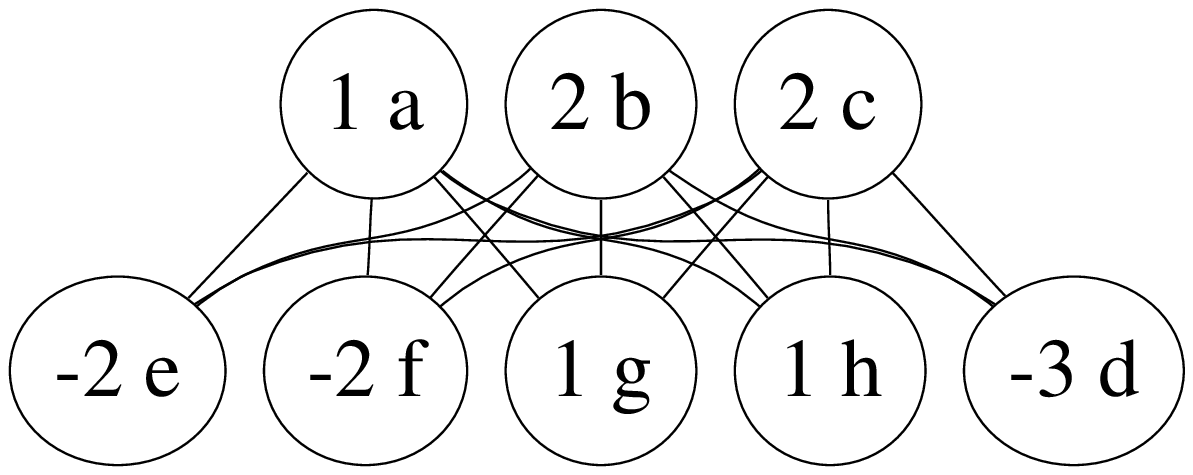}$

\item SFS $\left[D: \frac{1}{2}, \frac{1}{3}\right]$ U/m 
      SFS $\left[D: \frac{1}{2}, \frac{1}{2}, \frac{1}{2}\right]$ 
$m = \begin{pmatrix}  1 & 1 \\ 1 & 2\end{pmatrix}$ $H_1 = \Zed_2^2$.

Characteristic links $(\{a,c,g,h\}, \{a,c,d,e,g,h\}, \{a,c,d,f,g,h\}, \{a,c,e,f,g,h\})$,
$\vec \mu = (-\frac{1}{2}, 0,0,0)$.

Surgery diagram: $\includegraphics[height=1.5cm]{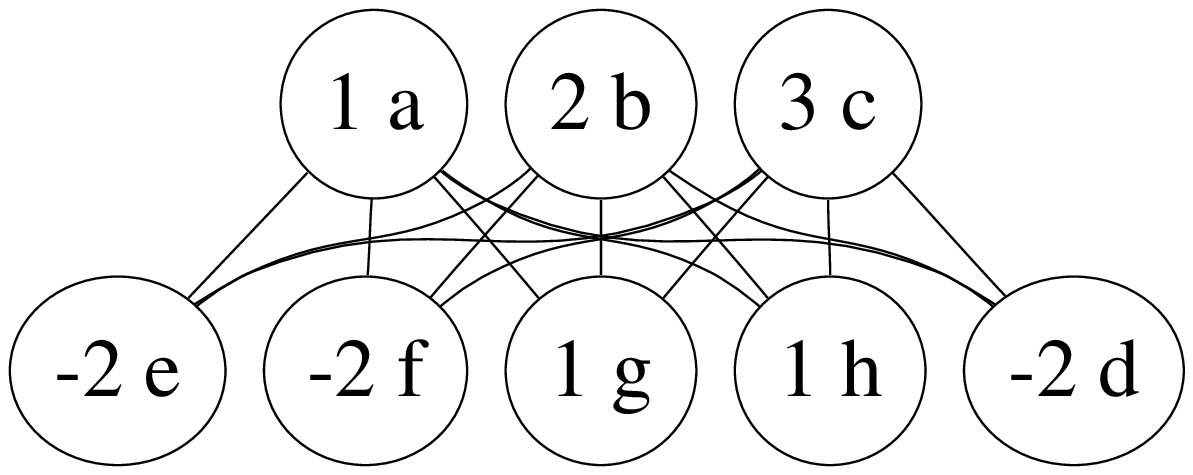}$

\item SFS $\left[D: \frac{1}{3}, \frac{2}{3}\right]$ U/m 
      SFS $\left[D: \frac{1}{2}, \frac{1}{2}, \frac{1}{2}\right]$ 
$m = \begin{pmatrix}  0 & 1 \\ 1 & 0\end{pmatrix}$ $H_1 = \Zed_6^2$.

Characteristic links $(\{c\}, \{d\}, \{e\}, \{c,d,e\})$,
$\vec \mu = (0,0,0,\frac{1}{2})$.

Surgery diagram: $\includegraphics[height=1.5cm]{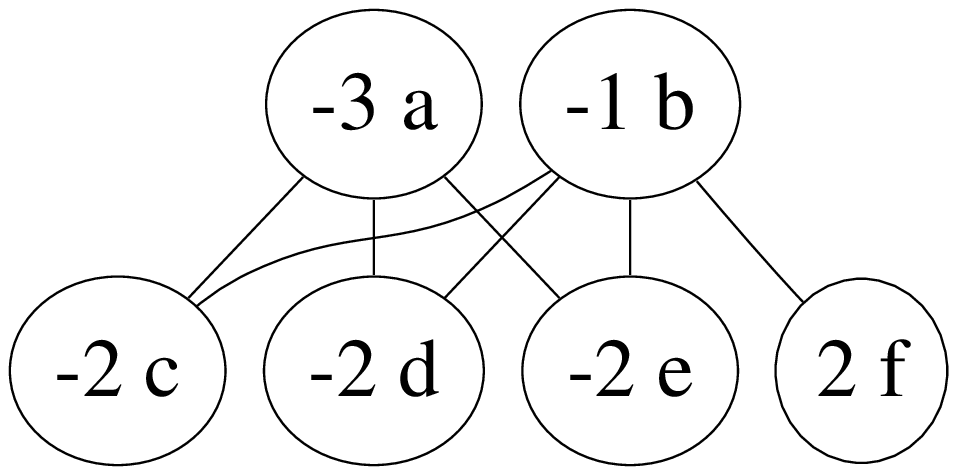}$

\item SFS $\left[D: \frac{2}{3}, \frac{2}{3}\right]$ U/m 
      SFS $\left[D: \frac{1}{2}, \frac{1}{2}, \frac{1}{2}\right]$ 
$m = \begin{pmatrix}  -1 & 1 \\ 1 & 0\end{pmatrix}$ $H_1 = \Zed_6^2$.

Characteristic links $(\{d\}, \{e\}, \{f\}, \{d,e,f\})$,
$\vec \mu = (0,0,0,\frac{1}{2})$.

Surgery diagram: $\includegraphics[height=1.5cm]{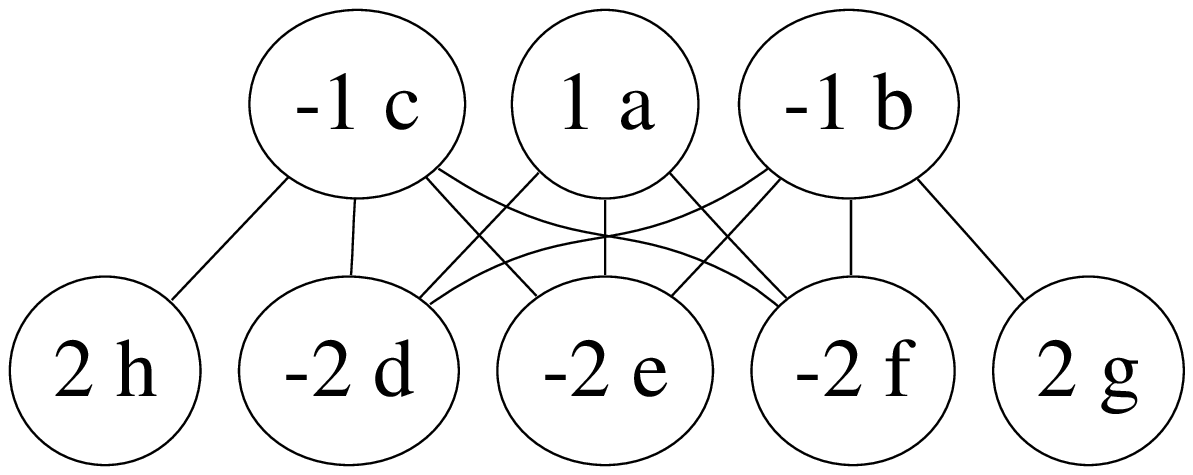}$

\item SFS $\left[D: \frac{1}{2}, \frac{1}{2}\right]$ U/m 
      SFS $\left[A: \frac{1}{2}\right]$ U/n 
      SFS $\left[D: \frac{1}{2}, \frac{1}{2}\right]$ 
$m = \begin{pmatrix} 1 & 1 \\ 1 & 0 \end{pmatrix}, 
 n = \begin{pmatrix}  0 & 1 \\ 1 & -1\end{pmatrix}$ $H_1 = \Zed_8^2$.
 
Characteristic links $(\{a,c,f,g\}, \{b,c,f,g\}, \{a,c,f,h\}, \{b,c,f,h\})$,
$\vec \mu = (0,0,0,0)$. 

Surgery diagram: $\includegraphics[height=1.5cm]{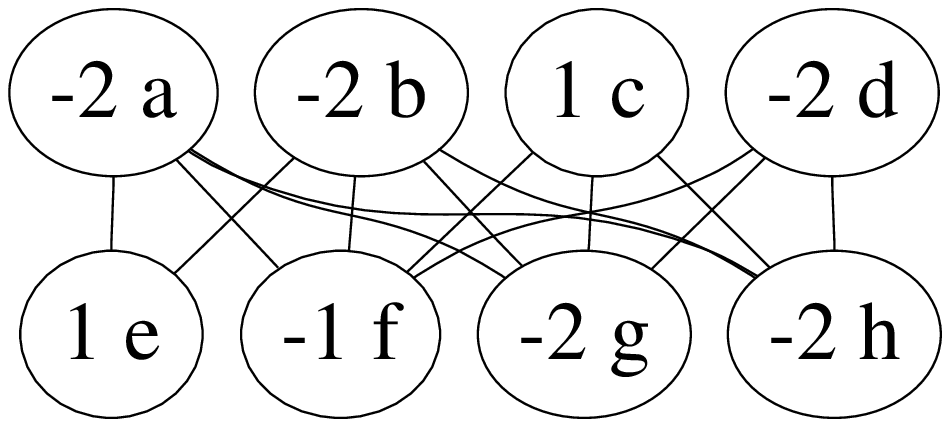}$ 

\item SFS $\left[D: \frac{1}{2}, \frac{1}{2}\right]$ U/m 
      SFS $\left[A: \frac{1}{2}\right]$ U/n 
      SFS $\left[D: \frac{1}{2}, \frac{1}{2}\right]$ 
$m = \begin{pmatrix}  1 & 1 \\ 1 & 0 \end{pmatrix}, 
 n = \begin{pmatrix}  1 & -1 \\ 0 & 1\end{pmatrix}$ $H_1 = \Zed_8^2$.
 
Characteristic links $(\{a,c,d,e,i,j\}, \{b,c,d,e,i,j\}, \{a,c,d,e,i,k\}, \{b,c,d,e,i,k\})$,
$\vec \mu = (0,0,0,0)$. 

Surgery diagram: $\includegraphics[height=1.5cm]{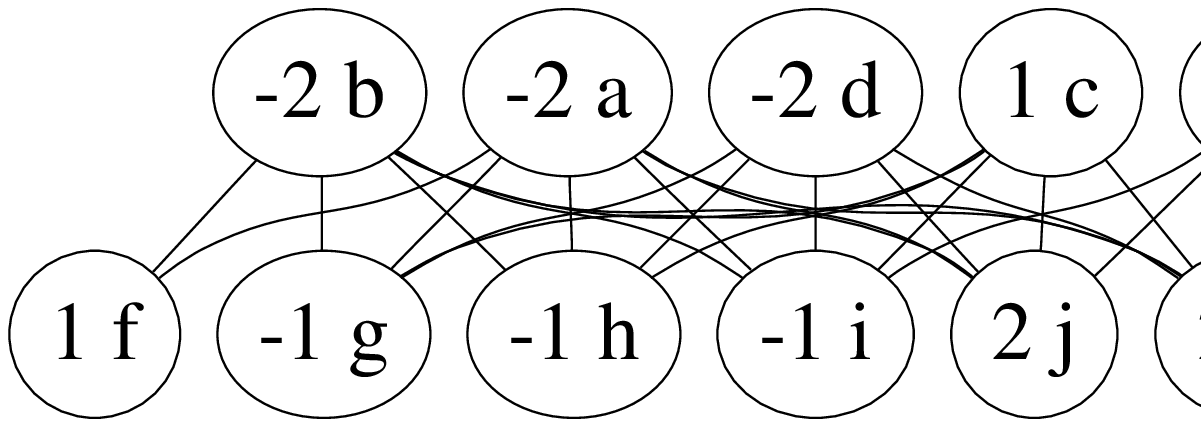}$

\item SFS $\left[D: \frac{1}{2}, \frac{1}{2}\right]$ U/m 
      SFS $\left[A: \frac{1}{3}\right]$ U/n 
      SFS $\left[D: \frac{1}{2}, \frac{1}{2}\right]$ 
$m = \begin{pmatrix}  0 & -1 \\ 1 & 0 \end{pmatrix}, 
 n = \begin{pmatrix}  0 & 1 \\ 1 & 0\end{pmatrix}$ $H_1 = \Zed_4^2$.
Characteristic links $(\{a,b,f,h\}, \{a,c,f,h\}, \{a,b,g,h\}, \{a,c,g,h\})$,
$\vec \mu = (0,0,0,0)$. 

Surgery diagram: $\includegraphics[height=1.5cm]{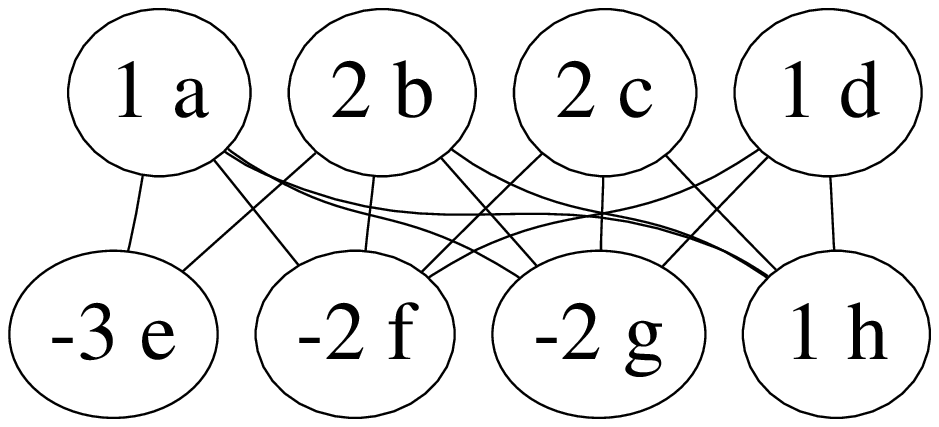}$

\item SFS $\left[D: \frac{1}{2}, \frac{1}{2}\right]$ U/m 
      SFS $\left[A: \frac{2}{3}\right]$ U/n 
      SFS $\left[D: \frac{1}{2}, \frac{1}{2}\right]$ 
$m = \begin{pmatrix}  0 & 1 \\ 1 & 0 \end{pmatrix}, 
 n = \begin{pmatrix}  1 & 1 \\ 1 & 0 \end{pmatrix}$ $H_1 = \Zed_4^2$. 
 
Characteristic links $(\{a,c,e,g\}, \{b,c,e,g\}, \{a,c,f,g\}, \{b,c,f,g\})$,
$\vec \mu = (0,0,0,0)$. 

Surgery diagram: $\includegraphics[height=1.5cm]{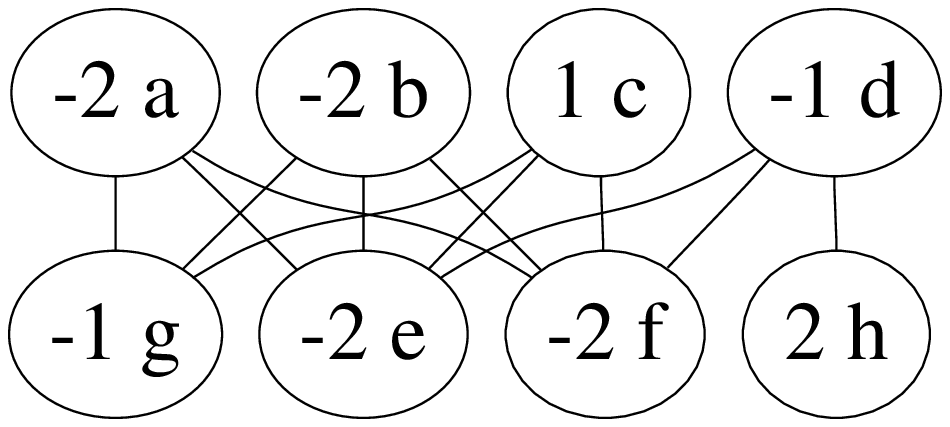}$

\vskip 5mm

\centerline{$\star$ Compound manifolds, $H_1$ infinite $\star$}
\centerline{ fibres over $S^1$ }
\vskip 5mm

\item SFS $\left[D: \frac{1}{2}, \frac{1}{2}\right]$ U/m 
      SFS $\left[D: \frac{1}{3}, \frac{2}{3}\right]$ 
$m = \begin{pmatrix}  -1 & 2 \\ 0 & 1\end{pmatrix}$  $H_1 = \Zed$. \ttc
$\Sigma_2 \rtimes S^1$.  The monodromy is reducible, differing
from the monodromy in {item \ref{item22} \S
\ref{embeddable_man}} by the square of a Dehn twist about a reduction
curve.  So although $\Delta(t) = (t^2-t+1)^2$, the monodromy extends over
a handlebody thus there are no signature obstructions to embedding in $S^4$.

\item SFS $\left[D: \frac{1}{2}, \frac{1}{2}\right]$ U/m 
      SFS $\left[D: \frac{1}{3}, \frac{2}{3}\right]$ 
$m = \begin{pmatrix}  -3 & 4 \\ -2 & 3\end{pmatrix}$ $H_1 = \Zed$.
$\Sigma_2 \rtimes S^1$.  The monodromy is reducible, differing
from the monodromy in {item \ref{item22} \S
\ref{embeddable_man}} by the $4$-th power of a Dehn twist about a reduction
curve.  So although $\Delta(t) = (t^2-t+1)^2$ the monodromy extends over a
handlebody thus there are no signature obstructions to embedding in $S^4$.

\item SFS $\left[D: \frac{1}{2}, \frac{1}{2}\right]$ U/m 
      SFS $\left[D: \frac{2}{5}, \frac{3}{5}\right]$ 
$m = \begin{pmatrix}  -1 & 2 \\ 0 & 1\end{pmatrix}$ $H_1 = \Zed$.
$\Sigma_4 \rtimes S^1$.  The monodromy is reducible and is given by the
composition of the map $(z_1,z_2) \longmapsto (e^{\frac{4\pi i}{5}}z_1,e^{\pi i} z_2)$ composed
with the square of a Dehn twist about a reduction
curve, thinking of the surface $\Sigma_4$ as in item \ref{item24} of \S \ref{embeddable_man}.
 So although $\Delta = (t^4-t^3+t^2-t+1)^2$ the monodromy extends over
a handlebody thus there are no signature obstructions to embedding in $S^4$.
\vskip 5mm
\centerline{$\star$ Compound manifolds, $H_1$ infinite $\star$}
\centerline{ do not fibre over $S^1$ }
\vskip 5mm

In order to compute the following Alexander polynomials we need to 
extend Lemma \ref{alexlem} by:
$$ \Delta SFS \left[ D: \frac{a}{b}, \frac{c}{d}, \frac{e}{f}\right] = 
\frac{ (t^{LCM(b,d,f)}-1)^2(t-1) }{ (t^{b'}-1)(t^{d'}-1)(t^{f'}-1) } $$
where $b' = \frac{LCM(b,d,f)}{b}, d'=\frac{LCM(b,d,f)}{d}, f'=\frac{LCM(b,d,f)}{f}$.

\item SFS $\left[D: \frac{1}{2}, \frac{1}{2}\right]$ U/m 
      SFS $\left[D: \frac{1}{2}, \frac{1}{2}, \frac{1}{2}\right]$ 
$m = \begin{pmatrix}  1 & 1 \\ 1 & 2\end{pmatrix}$ $H_1 = \Zed \oplus \Zed_2^2$.
$\Delta = (t^2+1)^2$.  The homology of the universal $\Zed$-cover has presentation
$\Zed[\Zed]/(t^2+1) \oplus \Zed[\Zed]/(t^2+1)$.  If we represent the generators by $a$ and $b$ then
$\langle a,a\rangle=\langle b,b\rangle=0$ and $\langle a,b\rangle = \frac{1}{t^2+1}$, 
which has all signatures equal to zero.

\item \label{unknown_num} SFS $\left[D: \frac{1}{2}, \frac{1}{2}\right]$ U/m 
      SFS $\left[D: \frac{1}{2}, \frac{1}{2}, \frac{1}{2}\right]$ 
$m = \begin{pmatrix}  -1 & 3 \\ -1 & 4\end{pmatrix}$  $H_1 = \Zed \oplus \Zed_2^2$.
$\Delta = (t^2+1)^2$.  Exactly as in the previous case, all signatures are zero.

\end{enumerate}

\section{Notes on computations \& notation}\label{hyperbolicity}

In order to deal with all the manifolds in the census efficiently, extensive use of computers 
was made while writing this paper. 

\begin{itemize}
\item The census of prime $3$-manifolds admitting a triangulation with 11 or less tetrahedra was 
created independently by Ben Burton, Sergei Matveev \cite{Matv} and also Bruno Martelli and Carlo 
Petronio.  Burton's software Regina \cite{BBurton} allows for relatively easy navigation of the census.  
We use the word {\it triangulated} to mean the smooth/PL manifold has been given a compatible 
{\it unordered delta complex} structure. Precisely, denote the $n$-simplex by 
$\Delta^n = \{(x_0,\cdots,x_n) \in \Real^{n+1} : x_i \geq 0 \ \forall i \text{ and } x_0+x_1+\cdots+x_n = 1\}$. 
 Given $i \in \{0,1,\cdots,n\}$ the $i$-th face map of $\Delta^n$ is the map 
 $f_i : \Delta^{n-1} \to \Delta^n$ given by 
 $$f_i(x_0,\cdots,x_{n-1}) = (x_0,x_1,\cdots,x_{i-1}, 0, x_i, x_{i+1}, \cdots, x_{n-1}).$$ 
 Given a permutation $\sigma \in \Sigma(\{0,1,\cdots,n\})$, the induced automorphism of 
 $\Delta^n$ is given by $\sigma_* : \Delta^n \to \Delta^n$, 
 $\sigma_*(x_0,x_1,\cdots,x_n) = (x_{\sigma(1)}, x_{\sigma(2)}, \cdots, x_{\sigma(n)})$.  
 An {\it unordered delta complex} is a CW-complex $X$ such that the domains of the attaching maps 
 are the boundaries of simplices (rather than discs), 
 $\phi : \partial \Delta^n \to X^{(n-1)}$, and for each $i$, the composite satisfies
 $\phi \circ f_i = \Phi \circ \sigma_*$ where $\Phi : \Delta^{n-1} \to X^{(n-1)}$ is a characteristic
  map of the $(n-1)$-skeleton, and $\sigma \in \Sigma(\{0,1,\cdots,n-1\})$ is some permutation. 
  If $\sigma$ is always the identity permutation, this would be an {\it ordered delta complex}. 

\item Surgery presentations for the closed hyperbolic $3$-manifolds in the census were created
using programs built from SnapPea \cite{Weeks} and Morwen Thistlethwaite's tables of knots and 
links. SnapPea allows one to drill a selection of geodesics out of a hyperbolic $3$-manifold, 
computing the canonical polyhedral decomposition on the resulting hyperbolic manifolds.  
The procedure used to find surgery presentations for closed hyperbolic 3-manifolds is to 
`randomize' the initial triangulation via a sequence of Pachner moves.  SnapPea then drills out
an initial curve in the  1-skeleton of the triangulation, resulting in a 1-cusped hyperbolic
manifold.  If that manifold is in the census of knots, the procedure terminates with a knot
surgery diagram.  If not, SnapPea is employed to give a list of drillable curves in the dual
1-skeleton of the cusped triangulation.  The software then systematically drills out up to
two additional geodesics, and then searches for the manifold in Thistlethwaite's table of 
hyperbolic link complements. SnapPea's isometry-checking routines determine the filling slopes
if a match is found among the link tables. 
\item Alexander polynomials of knots and smooth $4$-ball genus of many knots in the knot tables
can be looked up on Cha and Livingston's web page \cite{ChaLiv}. 
\item The Oszv\'ath-Szab\'o `$d$-invariant / correction term' for the Seifert-fibred rational homology spheres
in the census were computed using software written by Brendan Owens and Sa\v{s}o Strle.
\item The computation of the hyperbolicity of the torsion linking form was implemented by the
author in Regina since version 4.4. Details are given below.
\item Knots and links from tables are referred to via the notation $\KR{X}{C}$ where
$X$ indicates the table name $X=R$ indicates $C$ is taken from the Rolfsen table, 
$X=T$ indicates $C$ is taken from the Thistlethwaite table.  For example, $\KR{T}{2a_1}$
is the Hopf link and $\KR{R}{3_1}$ indicates the trefoil knot. A convenient place to 
view these tables is the Knot Atlas \cite{KA}. 
\end{itemize}

{\bf NCellularData} is a Regina class which implements the computation of the torsion linking 
form of a $3$-manifold and also tests its hyperbolicity via Kawauchi and Kojima's classification 
of symmetric bilinear forms on finite abelian groups taking values in $\Rat/\Zed$ \cite{KK}.
Given two elements $[v],[w] \in \tau H_1(M,\Zed)$, the torsion linking form
$\langle [v], [w]\rangle \in \Rat/\Zed$ is an intersection number. 
A multiple of $[v]$ is zero, $n[v] = 0$ for some $n \in \Zed$, so $nv = \partial S$ for 
some $2$-chain $S$.  Perturb $S$ and $w$ to intersect transversely, and let $m \in \Zed$ be 
the signed (algebraic) intersection number of $S$ and $w$.  The torsion linking form
is defined as $\langle [v],[w]\rangle = \frac{m}{n} \in \Rat/\Zed$.

The way this is implemented in Regina is to consider $v$ and $w$ as simplicial chains
in the simplicial chain-complex of $M$ coming from the triangulation.  
$M$ has a dual polyhedral-complex where the $i$-cells of the dual complex correspond to the 
$(3-i)$-cells of the triangulation. This is Poincar\'e's proof of his duality theorem \cite{ST}, 
simplified using CW-complexes.  For example, a $2$-cell in the dual polyhedral decomposition
corresponds to an edge $e$ of the triangulation.  Moreover, the $2$-cell is an $n$-gon, 
and the $n$-gon is a union of quadrilaterals, one quadrilateral for each time a tetrahedron
contains the edge $e$ ($e$ can be contained in a tetrahedron more than once since the triangulation
is semi-simplicial).  So the $2$-cells of the dual polyhedral decomposition intersect the
$1$-cells of the triangulation transversely. We homotope the identity map on $M$ to
be a cellular map from the triangulation to the dual polyhedral decomposition (this is the
core of the algorithm).  This allows us to express $v$ in the simplicial homology of the 
triangulation of $M$, and $w$ in the cellular homology of the dual polyhedral decomposition.  
So now $S$ is a simplicial $2$-chain and $w$ is a dual $1$-chain intersecting transversely, 
allowing for the computation of the intersection product via $\Zed$-linear algebra.

\vskip 5mm

The torsion linking form is stored as a square matrix of rational numbers, where the rows and 
columns are indexed by the invariant factors of $H_1(M,\Zed)$. The Kawauchi-Kojima classification 
of torsion linking forms \cite{KK} takes as input this matrix and determines hyperbolicity via 
linear-algebraic manipulations of the matrix. 
\vskip 5mm
\centerline{$\star$ Notation -- Regina's naming conventions for $3$-manifolds $\star$}
\vskip 5mm

For Seifert-fibred manifolds, Regina's notation is essentially the same as
Orlik's book \cite{Orlik}.  Given a surface $\Sigma$ let $M_\Sigma$ denote an orientable
$S^1$-bundle over $\Sigma$ with a section.   
The manifold $SFS \left[ \Sigma : \frac{a_1}{b_1}, \cdots, \frac{a_k}{b_k}\right]$ 
is obtained from $M_\Sigma$ by doing surgery on $k$ fibres in $M_\Sigma$, using filling
slopes $\frac{a_1}{b_1}, \cdots, \frac{a_k}{b_k}$ (slope zero being the slope of
the section).  If $\Sigma$ has boundary, the curves in $\partial M_\Sigma$ corresponding to the section 
will be denoted `$o$', and the curves corresponding to the fibre is denoted `$f$'. 

Only a few types of graph manifolds appear in the $11$-tetrahedron census. The underlying
graphs, if non-trivial, are of the form:
$$ \includegraphics[height=0.8cm]{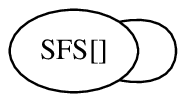}  \hskip 1cm
   \includegraphics[height=0.8cm]{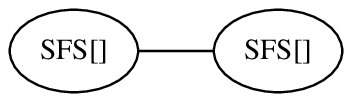}  \hskip 1cm
   \includegraphics[height=0.8cm]{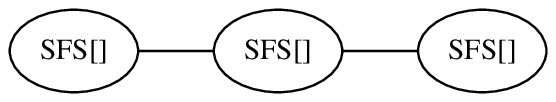}.$$

Meaning they have at most three vertices: one-vertex graphs have a single edge, and the
remaining two graph types are linear. Regina's convention for naming these manifolds are:

\begin{itemize}
\item Manifolds with a single non-separating torus, in this case Regina uses the
notation $M / \begin{pmatrix} a & b \\ c & d\end{pmatrix}$ where $M$ is a Seifert
fibred manifold with two boundary tori.  This indicates that we glue to two boundary tori together
so that $f_2$ is identified with $af_1+bo_1$, and $o_2$ is identified with
$cf_1+do_1$. 
\item Manifolds with a single separating torus are denoted $M_1$ U/m $M_2,  
m = \begin{pmatrix} a & b \\ c & d\end{pmatrix}$ where $M_1$ and $M_2$ 
follow the notation for Seifert-fibred manifolds above. 
The matrix $m$ indicates that $\partial M_2$ is glued to $\partial M_1$ by a map
that identifies $f_2$ with $af_1+bo_1$ and
$o_2$ with $cf_1+do_1$.  
\item The remaining class of manifolds have the form $M_1$ U/m $M_2$ U/n $M_3$, 
$m = \begin{pmatrix} a_1 & b_1 \\ c_1 & d_1 \end{pmatrix}$
$n = \begin{pmatrix} a_2 & b_2 \\ c_2 & d_2 \end{pmatrix}$.
The matrices $m$ and $n$ denote the gluing maps 
$m : \partial M_2 \to \partial M_1$ and $n : \partial M_3 \to \partial M_2$,
precisely $m(f_2) = a_1f_1 + b_1o_1$, $m(o_2) = c_1f_1 + d_1o_1$ and
$n(f_3) = a_2f_2+b_2o_2$, $n(o_3)=c_2f_2+d_2o_2$.
\end{itemize}
There are two other classes of manifolds assigned special names by Regina:
\begin{itemize}
\item Hyperbolic manifolds are named in a somewhat ad-hoc way.  The first part of
such a manifold's name is the initial 8 terms of the decimal expansion of the volume of the 
manifold, followed by the invariant factor decomposition of its first homology
group.  If this data does not uniquely identify the manifold in the census, an additional
identifier of the shortest geodesic length is given, suitably rounded.
\item If the manifold fibres over $S^1$ with fibre a torus, the manifold is denoted
by the notation $T \times I / \begin{pmatrix} a & b \\ c & d\end{pmatrix}$ where the matrix
describes the monodromy (assuming the tori are parametrized so as to be parallel). 
In these notes such manifolds are denoted
$(S^1 \times S^1) \rtimes_{\begin{pmatrix} a & b \\ c & d\end{pmatrix}} S^1$.
\end{itemize}

\vskip 5mm
\centerline{$\star$ Surgery presentations of graph manifolds $\star$}
\vskip 5mm

The technique used to construct surgery presentations is relatively primitive but effective. 
Lickorish's proof that $3$-manifolds have surgery presentations had a key idea about gluings of manifolds.  
Let $M$ and $N$ be disjoint $3$-manifolds and $f : \partial M \to \partial N$ a diffeomorphism. 
Let $M \cup_f N$ be the manifold obtained by gluing $\partial M$ to $\partial N$
along $f$.  Let $c$ be a curve in $\partial N$, and let $D_c : \partial N \to \partial N$
be the positive Dehn twist about $c$, then $M \cup_{D_c \circ f} N \simeq M \cup_f N'$ where 
$N'$ is the manifold obtained from $N'$ by doing a $\pm 1$-Dehn surgery along a curve $c'$ in the
interior of $N$, parallel to $c$. 

For example, consider item \ref{firsteg} of \S \ref{nonembeddable_man}. The manifolds
SFS $\left[D: \frac{1}{2}, \frac{1}{2}\right]$ and 
SFS $\left[D: \frac{1}{2}, \frac{1}{3}\right]$ should be thought of as the `Dehn surgery' 
for the partially framed links
\vskip 5mm
\centerline{
{
\psfrag{T}[tl][tl][0.8][0]{$-2$}
\psfrag{R}[tl][tl][0.8][0]{$-2$}
$\includegraphics[width=4cm]{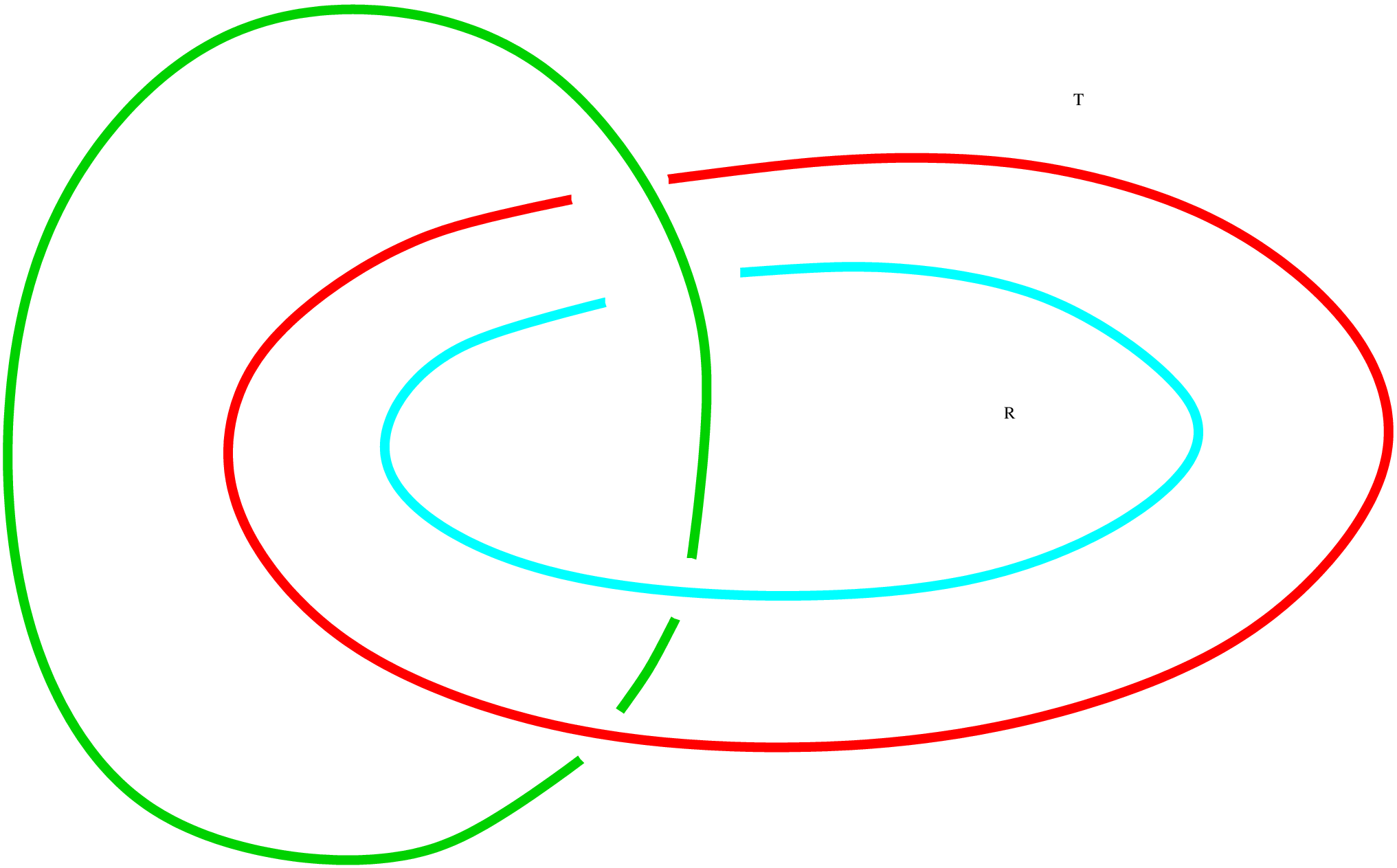}$
}
\hskip 5mm and \hskip 5mm
{
\psfrag{T}[tl][tl][0.8][0]{$-3$}
\psfrag{R}[tl][tl][0.8][0]{$-2$}
$\includegraphics[width=4cm]{ssf.eps}$
}
}

respectively, in the sense that the unlabelled (green) curves are drilled but not filled.  
Let $f_2, o_2$ and $f_1, o_1$ denote the fibre and base boundary curves for the manifolds
SFS $\left[D: \frac{1}{2}, \frac{1}{2}\right]$ and 
SFS $\left[D: \frac{1}{2}, \frac{1}{3}\right]$ respectively.  Then as 
a map from the boundary of the first manifold to the boundary of the second
our gluing map has the form $\begin{pmatrix}  -3 & -2 \\ 4 & 3\end{pmatrix}$ (i.e. the transpose
of the matrix in item \ref{firsteg}) where the bases curves in the domain are $\{f_1,o_1\}$
and in the range $\{f_2,o_2\}$.  Now we `splice' the above two surgery diagrams together
using Lickorish's idea -- in particular we compare this union of two solid tori to the genus one
Heegaard splitting of $S^3$.  So multiply the gluing map on the left by $\begin{pmatrix}  0 & 1 \\ 1 & 0\end{pmatrix}$
and write this matrix as a product of row/column operations: 
$$\begin{pmatrix}  0 & 1 \\ 1 & 0\end{pmatrix}\begin{pmatrix}  -3 & -2 \\ 4 & 3\end{pmatrix}
= \begin{pmatrix}  4 & 3 \\ -3 & -2\end{pmatrix} = \begin{pmatrix}  1 & -1 \\ 0 & 1\end{pmatrix}
\begin{pmatrix} 1 & 0 \\ -3 & 1\end{pmatrix}\begin{pmatrix}  1 & 1 \\ 0 & 1\end{pmatrix}.$$ 
We think of this product as $D_{m_2}^{-1} \circ D_{l_2}^3 \circ D_{m_2}$ i.e. a product of powers
of positive Dehn twists about the standard meridians and longitudes in a solid torus in the standard genus $1$
Heegaard splitting of $S^3$.  This gives us the `spliced' surgery presentation for
the manifold in item \ref{firsteg}.

{
\psfrag{a}[tl][tl][0.8][0]{$-2$}
\psfrag{b}[tl][tl][0.8][0]{$-2$}
\psfrag{c}[tl][tl][0.8][0]{$1$}
\psfrag{d}[tl][tl][0.8][0]{$1$}
\psfrag{e}[tl][tl][0.8][0]{$1$}
\psfrag{f}[tl][tl][0.8][0]{$-3$}
\psfrag{g}[tl][tl][0.8][0]{$-2$}
\psfrag{h}[tl][tl][0.8][0]{$1$}
\psfrag{i}[tl][tl][0.8][0]{$-1$}
$$\includegraphics[width=8cm]{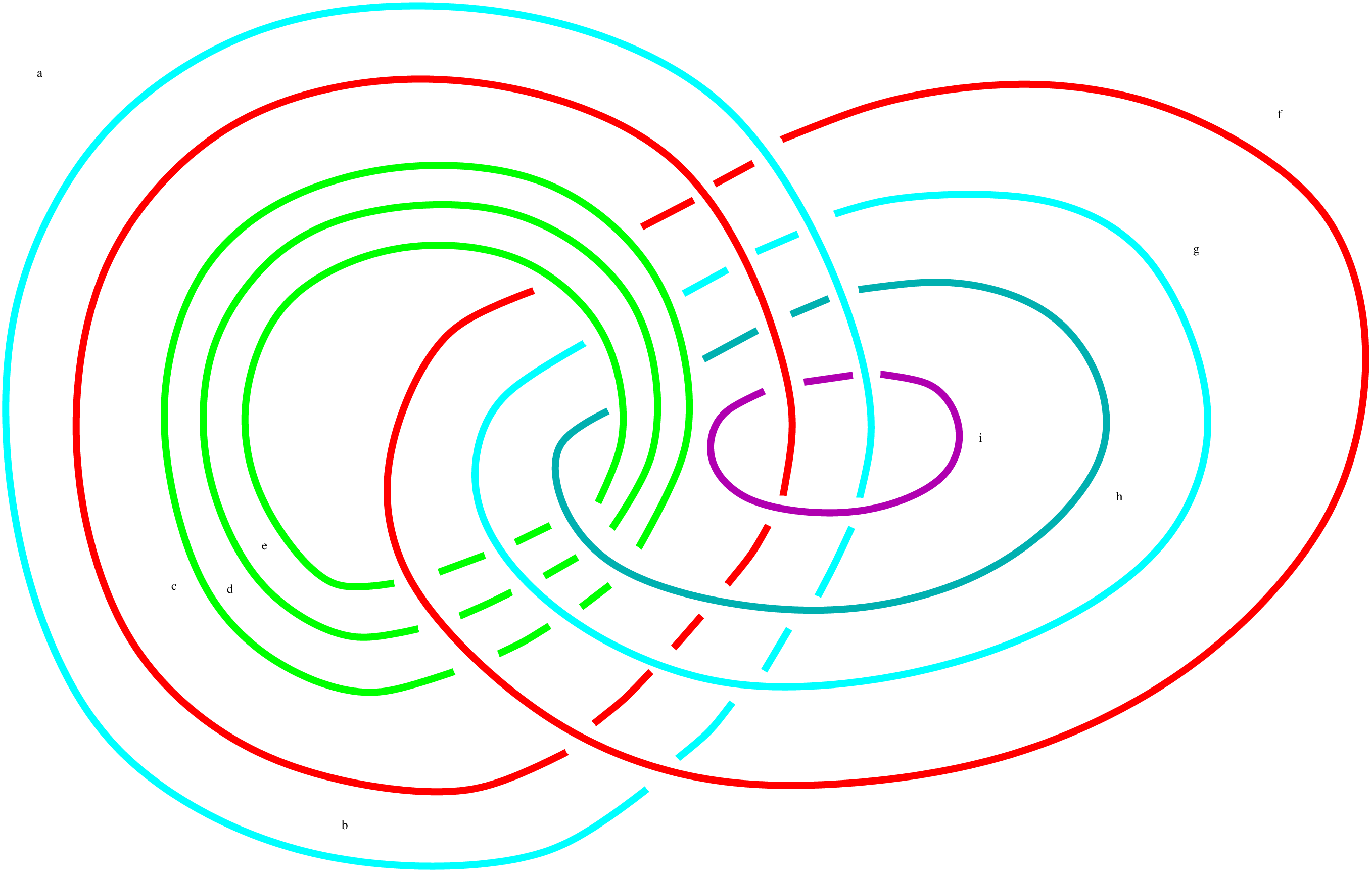}$$
\centerline{SFS $\left[D: \frac{1}{2}, \frac{1}{2}\right]$ U/m 
      SFS $\left[D: \frac{1}{2}, \frac{1}{3}\right]$ 
$m = \begin{pmatrix}  -3 & 4 \\ -2 & 3\end{pmatrix}$}
}

A similar computation shows that the manifolds
SFS $\left[D: \frac{a_1}{b_1}, \frac{a_2}{b_2}\right]$ U/m 
      SFS $\left[D: \frac{a_3}{b_3}, \frac{a_4}{b_4}\right]$ 
$m = \begin{pmatrix}  \alpha & \beta \\ \gamma & \delta\end{pmatrix}$
all have integral surgery presentations along links $L$ which decompose
into a union of disjoint sub-links $L = L_1 \sqcup L_2$ where
$L_1 \subset \Real^2 \times \{0\}$ is a collection of nested circles of various
radii centered around points in $\{0\} \times \Real \times \{0\}$ and
$L_2 \subset \{0\} \times \Real^2$ is also a collection
of nested circles, centered around points in $\{0\}\times \Real \times \{0\}$. 
Such surgery presentations are perhaps most easily represented via
a graph, analogous to a plumbing diagram, which represents the framing/linking
matrix:
$$\includegraphics[height=3cm]{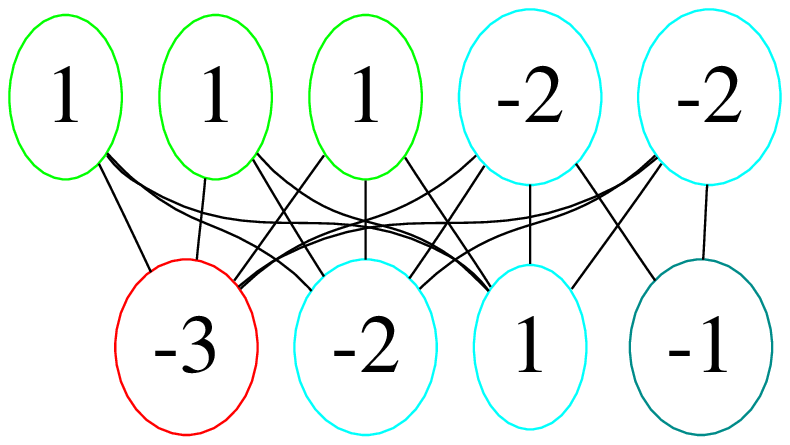}.$$

\centerline{$\star$ Computing the monodromy from the Seifert data $\star$}

These are the fibre bundles over $S^1$ with fibre a closed surface of genus $g \geq 2$, such that 
the monodromy is a finite-order diffeomorphism of the surface.  Denote such a manifold by 
$\Sigma_g \rtimes_{\Zed_n} S^1$ where $n$ is the order of the monodromy. Precisely, if $f : \Sigma_g \to \Sigma_g$ 
denotes the monodromy, $\Sigma_g \rtimes_{\Zed_n} S^1 = (\Sigma_g \times S^1)/\Zed_n$ where $\Zed_n$ acts on
$\Sigma_g \times S^1$ by $e^{\frac{2\pi ik}{n}}.(x,z) = (f^{(k)}(x),e^{\frac{2\pi ik}{n}}z)$ 
where we make the identification $\Zed_n \equiv \{e^{\frac{2\pi i k}{n}} : k \in \Zed \}$.
These manifolds are all Seifert fibred -- the fibring being covered by the product fibring of 
$\Sigma_g \times S^1$. The fibre $\Sigma_g$ is the unique horizontal incompressible surface, thus these manifolds 
all have the form SFS $\left[ \Sigma_m : \frac{a_1}{b_1}, \cdots, \frac{a_k}{b_k} \right]$ where $\sum_{i=1}^k \frac{a_i}{b_i} = 0$. 
Thus $n = LCM\{b_1,b_2,\cdots,b_k\}$. $k$ is the number of non-free orbits of $\Zed_n$ acting on $\Sigma_g$
and $\chi(\Sigma_g) = n(\chi(\Sigma_m) + \sum_{i=1}^k (\frac{1}{b_i} - 1))$. The numbers $b_i$ give the cone angles
$2\pi/b_i$ for the singular orbits of $\Zed_n$ acting on $\Sigma_g$.  For example, items
\ref{it8} through \ref{it11} in \S \ref{nonembeddable_man} all have the form
SFS $\left[S^2 : \frac{\alpha_1}{\beta_1}, \frac{\alpha_2}{\beta_2}, \frac{\alpha_3}{\beta_1\beta_2}\right]$ where
$GCD(\beta_1,\beta_2)=1$.  The obstruction to show these manifolds
do not embed (in any homology sphere) is the Alexander polynomial.  For these manifolds an efficient way of
computing the Alexander polynomial is by constructing an equivariant CW-decomposition of the fibre -- and to consider the Alexander polynomial to be the order ideal of the homology of the fibre as a
module $\lau$-module, where $\Zed$ acts via the monodromy. Since
the base space is $S^2$ with three singular points, consider it to be a square with an identification made to
the edges. This square lifts to a CW-decomposition of the fibre, and in this case the cell structure reduces to
one with $\beta_1+\beta_2$ 0-cells, $\beta_1\beta_2$ 1-cells and a single $2$-cell.  The monodromy has a fixed
point which is the centre of the 2-cell, and the remaining singular points are the 0-skeleton, allowing a rather
direct computation of the Alexander polynomial. Checking that the Alexander polynomial does not have the
form $p(t)p(t^{-1})$ can be done readily by using computer algebra software (such as Pari) to compute the roots in $\C$. See Theorem \ref{KawCond}.
\vskip 5mm
\centerline{$\star$ A technique of Casson and Harer $\star$}
\vskip 5mm
In their paper Casson and Harer \cite{CH} demonstrate a technique to find contractible $4$-manifolds
bounding $3$-manifolds. We show here how this technique allows us to find embeddings of a certain class of 
$3$-manifolds in homotopy $4$-spheres. Take for example manifold \ref{cheg} from the list in \S 
\ref{embeddable_man}, this is $(-5,-5)$-surgery on the Whitehead link. 

{
\psfrag{f1}[tl][tl][0.8][0]{$-5$}
\psfrag{f2}[tl][tl][0.8][0]{$-5$}
\psfrag{f3}[tl][tl][0.8][0]{$-1$}
\psfrag{t5}[tl][tl][0.8][0]{$-1$}
\psfrag{t4}[tl][tl][0.8][0]{$-4$}
\psfrag{s1}[tl][tl][0.8][0]{Step 1}
\psfrag{bd}[tl][tl][0.8][0]{Step 2}
\psfrag{s3}[tl][tl][0.8][0]{Step 3}
$$\includegraphics[width=14cm]{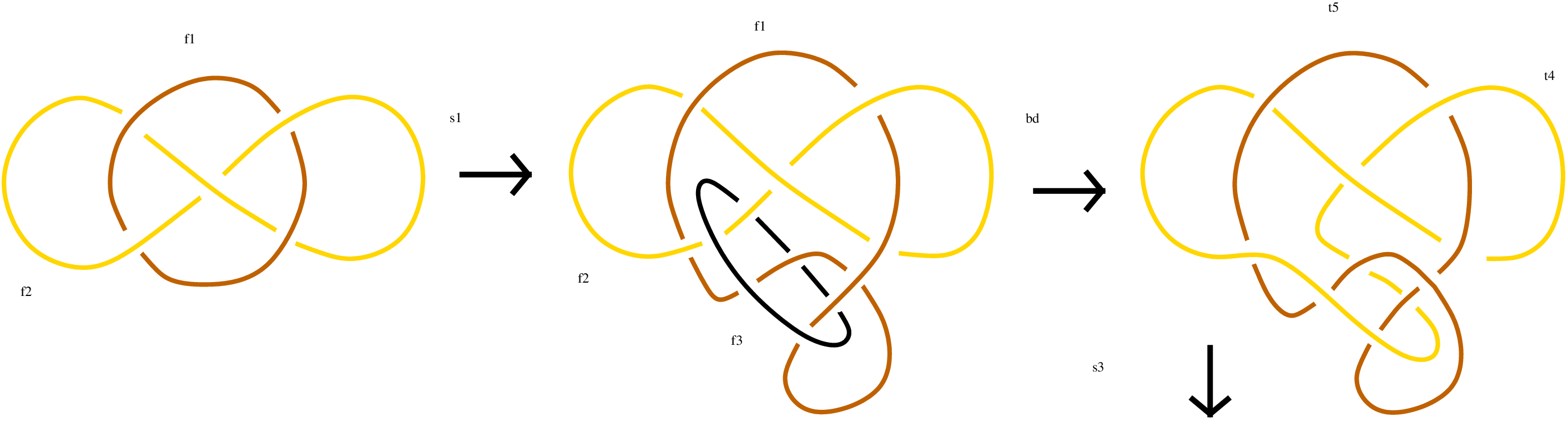}$$
}
{
\psfrag{f1}[tl][tl][0.8][0]{$-4$}
\psfrag{f2}[tl][tl][0.8][0]{$-1$}
\psfrag{f3}[tl][tl][0.8][0]{$0$}
\psfrag{s4}[tl][tl][0.8][0]{Step 4}
\psfrag{s5}[tl][tl][0.8][0]{Step 5}
$$\includegraphics[width=14cm]{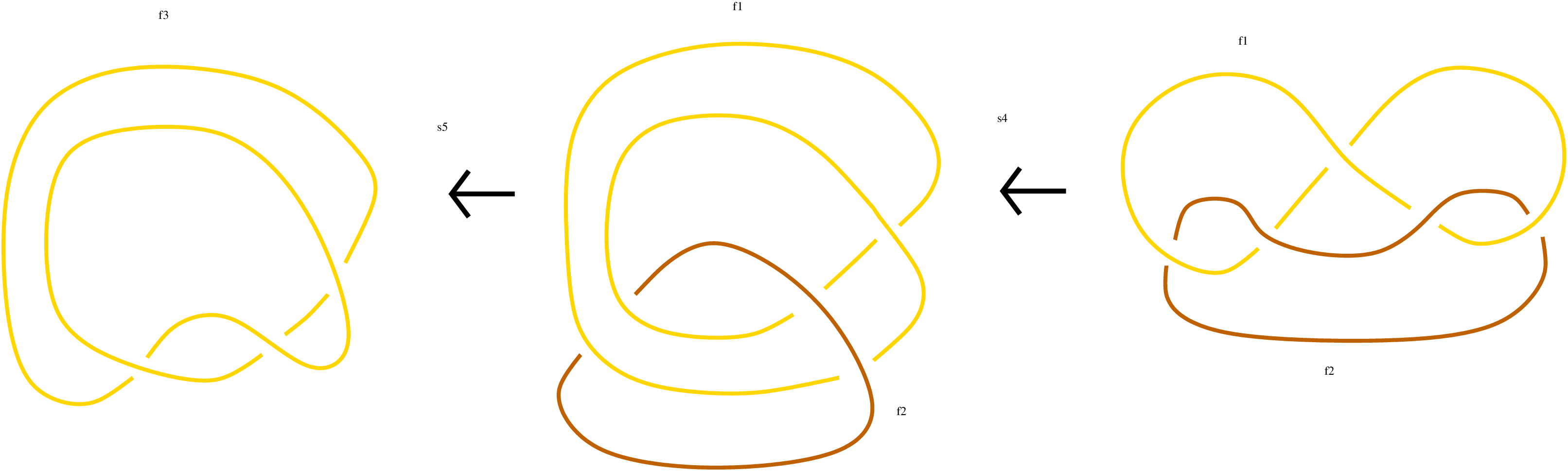}$$
}
The above figure starts off with $(-5,-5)$-surgery on the Whitehead Link, call this
manifold $M$.  Think of {\bf Step 1} as representing a handle attachment to $M \times [0,1]$
on the side of $M \times \{1\}$. {\bf Step 2} represents the Kirby `blow down' move. 
{\bf Step 3} is an isotopy. {\bf Step 4} a further `fold' isotopy.  {\bf Step 5} is a further
`blow down' equivalence of handle presentations. This leaves us with the manifold $S^1 \times S^2$
on the boundary, which we attach a $3$-handle and then a $4$-handle.  In summary, we have
attached a $2$-handle, then a $3$-handle and $4$-handle to $M \times [0,1]$ to construct a manifold
$W_1$ bounding $M \times \{1\}$.  By design $\pi_1 W_1 = \Zed_5$ and the inclusion 
$H_1 (M\times \{1\}) \to H_1 W_1$ has
kernel one of the summands of the hyperbolic splitting $H_1 M = \Zed_5 \oplus \Zed_5$.  By symmetry
of the Whitehead link which switches components, we can repeat the argument on the 
$M \times \{0\}$ side of $M \times [0,1]$, building
a manifold $W_0$ such that the inclusion $M \times \{0\} \to W_0$ kills the
complementary summand of the hyperbolic splitting.  
The union of these two bounding manifolds $W_0 \cup W_1$ is then a homotopy $4$-sphere containing
$M$.

\vskip 5mm
\centerline{$\star$ A symmetric embedding for $S^3 / Q_8$ $\star$}
\vskip 5mm
We describe a particularly symmetric embedding of $S^3 / Q_8$ into $S^4$ that we learned from Rob Kusner. 

Let $\mathcal A$ be the traceless, symmetric $3 \times 3$ real matrices, 
$$\mathcal A = \{ A \in M_{3 \times 3}(\Real) : tr(A) = 0, A^t = A \}.$$

The set $\mathcal A$ is a $5$-dimensional inner product space with the inner product
defined by trace, transpose and the matrix product 
$$\langle A_1, A_2 \rangle = tr( A_1^t A_2 ).$$  

The group $SO_3$ acts
on $\mathcal A$ by conjugation, moreover this action is by isometries.  Thus the
action restricts to the unit sphere of $\mathcal A$, $S \mathcal A \equiv S^4$.  Symmetric matrices
are diagonalizable by orthogonal matrices so the conjugacy class of such a matrix is determined
by its three real eigenvalues $-1 \leq \lambda_1 \leq \lambda_2 \leq \lambda_3 \leq 1$. 
The trace condition $tr(A) = 0$ is equivalent to the condition $\lambda_1 + \lambda_2 + \lambda_3 = 0$ and the
unit sphere condition is similarly equivalent to the condition $\lambda_1^2 + \lambda_2^2 + \lambda_3^2 = 1$.
  
The orbit decomposition of $SO_3$ acting on $S\mathcal A$ consists of two orbits equivalent
to $\mathbb RP^2$ (the $\lambda_1 = \lambda_2$ subspace and the $\lambda_2 = \lambda_3$ subspace) and a 
$1$-parameter family of orbits equivalent to $S^3 / Q_8$, these are the matrices with distinct 
eigenvalues.  The orbits corresponding to distinct eigenvalues isomorphic to $SO_3 / G$ where $G$ is the group of rotations by $\pi$ in the faces of a cube, i.e. the rotations preserving the eigenspaces.  
This group lifts to $Q_8 = \{\pm 1, \pm i, \pm j, \pm k\}$ in $S^3$, thus these orbits are all 
embedded manifolds diffeomorphic to $S^3 / Q_8$. 

\section{Observations and questions from the data}\label{qsec}

A striking feature about the data is that some $3$-manifolds from the census are more susceptible 
to our embedding constructions than others.  For example, if the manifold fibres over $S^1$, we have 
deform-spun embeddings and surgical embeddings at our disposal.  Seifert-fibred spaces have a variety 
of embedding techniques, largely due to Crisp and Hillman.  But when dealing with hyperbolic manifolds, 
the only technique used is the surgical embedding construction.

\begin{question}Do there exist $3$-manifolds $M$ which embed smoothly in 
$S^4$ such that no embedding of $M$ in $S^4$ is a surgical embedding in the sense
of Constructions \ref{surgical} and \ref{tenh}?  
\end{question}

For the above question, I know of no relevant obstructions, although presumably the answer is yes.
Similarly, if the rank of the 1st homology group of $M$ is larger than one, it's not clear if there
are any concordance obstructions that one can use. Given $M$ in $S^4$, let $S^4 = V_1 \cup_M V_2$ be the 
decomposition of $S^4$ into two $4$-manifolds along their common boundary $M$. 

\begin{question}\label{homologyRestrictions}
If the rank of $H_1 M$ is larger than $1$, are there any restrictions on the ranks of $H_1 V_1$ and 
$H_1 V_2$? Similarly, an embedding of $M$ in $S^4$ induces a hyperbolic splitting on $\tau H_1 M$.  How many 
hyperbolic splittings can one generate via embeddings?  Are some hyperbolic splittings of $\tau H_1 M$ 
impossible to realize via embeddings? 
\end{question}

There is at least one such restriction.  The inclusion $M \to V_i$ induces a restriction map on 
cohomology $H^k(V_i) \to H^k(M)$.  So if $\alpha,\beta \in H^1(V_i)$ restrict to classes in $H^1(M)$ 
which have a non-zero cup-product, they must also have a non-zero cup product in $H^2(V_i)$.  This 
gives a restriction in some cases -- for example $M=(S^1)^3$.  Since $H^1(M)$ has non-trivial 
cup-products, it is impossible for $H^1(V_1)$ to have rank three, since $H^2(V_1)$ would 
necessarily have rank zero.  Thus for any embedding of $(S^1)^3$ in $S^4$, 
$rank(H_1(V_i)) \geq 1$ for $i=1,2$. 

One startling observation from the data in this paper and from the references is that there are as of 
yet no examples of $3$-manifolds that embed in homology $4$-spheres which do not embed in $S^4$.  
This leads to two questions.

\begin{question}
If a $3$-manifold $M$ admits a smooth embedding into a homotopy $4$-sphere, does it admit a smooth 
embedding $S^4$? Are there $3$-manifolds that embed in homology $4$-spheres which do not embed in $S^4$?
\end{question}

The earlier question is only interesting if the smooth $4$-dimensional Poincar\'e
conjecture is false.  But it is perhaps surprising that it's not immediately clear
whether an embedding of a 3-manifold into an exotic $S^4$ could be pushed into
a standard $4$-ball. Agol and Freedman \cite{Agol-Freedman} have taken a step towards resolving
this question, giving an obstruction to a $3$-manifold embedding smoothly in $S^4$ in 
terms of the handlebody metric on the curve complex. Technically, the Agol-Freedman obstruction 
obstructs a Heegaard splitting induced by the embedding.  It is unclear, at present, if it can
be used to obstruct embeddings. 

One would think there should be $3$-manifolds that embed in homology $4$-spheres that do not
embed in $S^4$. I am unaware of any obstructions at present.  A reasonable place to look for 
answers to this question would be homology $3$-spheres.  Let $M$ be a homology $3$-sphere. As we 
have observed $M \# (-M)$ embeds smoothly in a homology $4$-sphere but it is not clear $M \# (-M)$ 
embeds in $S^4$ unless we could realize $M$ as something like a cyclic branched cover on a knot
in $S^3$ or some Litherland-style variant on that theme (see Theorem \ref{DSthm}).  
This provides a source of $3$-manifolds that embed in homology $4$-spheres but 
for which there is no clear embedding in $S^4$.

A reoccurring problem in this paper is that even if a $3$-manifold embeds in $S^4$, we have no
uniform, standard way of constructing an embedding. 

\begin{question} Is there an efficient procedure to determine whether or not a triangulated
$3$-manifold admits a locally-flat PL-embedding (equivalently, smooth embedding) in $S^4$?
\end{question}

Costantino and Thurston have recently developed an efficient procedure \cite{CT}
to construct a triangulated $4$-manifold that bounds a triangulated $3$-manifold.  They
do this by perturbing a map $M \to \Real^2$ associated to the triangulation, and
`filling in' the level sets in a natural way.  Perhaps one could devise a combinatorial 
search for embeddings $M \to \Real^4$ by considering such an embedding to be a special pair of generic 
maps $M \to \Real^2$?

\begin{question}
Is there a computable function $\beta : \mathbb N \to \mathbb N$ such that for each $3$-manifold 
that embeds in $S^4$ and admits a triangulation with $n$ tetrahedra, $M$ appears as a vertex-normal
solution to the gluing equations for a triangulation of $S^4$ with no more than $\beta(n)$ pentachora in 
the triangulation of $S^4$? 
\end{question}

Provided we had such a $\beta$, the problem of determining whether or not a $3$-manifold embeds in $S^4$
would be an algorithmically-solvable problem, as there would be a finite list of triangulations of $S^4$
on which to do normal surface enumeration. 

On the pessimistic side, Dranishnikov and Repovs \cite{DR} have shown there exists a smooth embedding
of a $3$-manifold $M$ in $S^4$ such that $S^4 = V_1 \cup_M V_2$ with $\pi_1 V_i$ having an unsolvable
word problem, for some $i \in \{1,2\}$.  Thus if one attempts to find obstructions to $M$ embedding in
$S^4$ based on the fundamental group, one could run into computability problems unless the obstruction
is based on a computable invariant of group presentations.  Computable invariants of group presentations
include things like computable invariants of representation varieties, and the lower central series of the group. 

\begin{question}
If $M$ admits a smooth embedding into $S^4$, does it admit an embedding where $S^4 = V_1 \cup_M V_2$ with
both $\pi_1 V_1$ and $\pi_1 V_2$ having solvable word problems? 
\end{question}

\begin{question}(M. Freedman) Given a smooth $3$-manifold $M$, if $M \# (S^1 \times S^2)$ embeds in $S^4$, does $M$? 
More generally, does stabilization via connect-sum with copies of $S^1 \times S^2$ make the embedding 
problem any easier? 
\end{question}

The question highlights a technical issue with the kinds of invariants we use to obstruct embedding. 
All the invariants we use are additive under connect sum.   

\providecommand{\bysame}{\leavevmode\hbox to3em{\hrulefill}\thinspace}

\Addresses

\begin{thebibliography}{Lan71}

\bibitem{AkKi}
S.~Akbulut, R.~Kirby, 
\emph{Mazur manifolds,}
Michigan Math. J., 26 (1979), 259-284.

\bibitem{Agol}
I.~Agol,
\emph{Bounds on exceptional Dehn filling.}
Geometry \& Topology {\bf 4} (2000) 431--449.

\bibitem{Agol-Freedman}
I.~Agol, M.~Freedman, 
\emph{Embedding heegaard decompositions.}
{\tt [arXiv:1906.03244]} preprint.

\bibitem{Bott}
R.~Bott,
\emph{Nondegenerate critical manifolds.}
Ann. Math. (2) {\bf 60}, 248-261 (1954).

\bibitem{Budney}
R.~Budney,
\emph{A family of embedding spaces,}
Geometry and Topology Monographs {\bf 13} (2008) 41--83. 

\bibitem{BBurton}
B.~Burton, R.~Budney, W.~Pettersson, et al.,
\emph{Regina: Software for 3-manifold topology and normal surface theory,}
{\tt [http://regina.sourceforge.net/]}, 1999--2020. 

\bibitem{Weeks}
P.J.~Callahan, J.C.~Hildebrand, J.R.~Weeks, 
\emph{A Census of Cusped Hyperbolic 3-Manifolds}, 
Mathematics of Computation {\bf 68/225}, 1999.

\bibitem{CH}
A.~Casson, J.~Harer.
\emph{Some homology lens spaces which bound rational homology balls.}
Pacific. J. Math. Volume {\bf 96}, Number 1 (1981), 23-36.

\bibitem{Cerf}
J.~Cerf,
\emph{Sur les diff\'eomorphismes de la sph\`ere de dimension 
trois $(\Gamma \sb{4}=0)$}, 
Lecture Notes in Mathematics, No. 53. 
Springer-Verlag, Berlin-New York 1968.

\bibitem{ChaLiv}
J.C.~Cha, C.~Livingston,
\emph{Unknown values in the table of knots,}
{\tt arXiv [math.GT/0503125]}.

\bibitem{CT}
F.~Costantino, D.~Thurston, 
\emph{3-manifolds efficiently bound 4-manifolds,}
Journal of Topology {\bf 1} (3) 703--745 (2008). 
{\tt arXiv [math.GT/0506577]}.

\bibitem{CGLS}
M.~Culler, C.~Gordon, J.~Luecke, P.~Shalen, 
\emph{Dehn surgery on knots.} 
Ann. of Math. (2) {\bf 125} (1987), no. 2, 237--300.

\bibitem{CrispH}
J.S.~Crisp, J.A.~Hillman,
\emph{Embedding Seifert fibred $3$-manifolds and ${\rm Sol}\sp 3$-manifolds
in $4$-space,}
Proc. London Math Soc. (3) (1998), no. {\bf 3} 685--710.

\bibitem{ADon}
A.~Donald,
\emph{Embedding Seifert manifolds in $S^4$,}
Trans. Amer. Math. Soc. {\bf 367} (2015), 559-595. 
{\tt [arXiv:1203.6008]}. 

\bibitem{DR}
A.~Dranisnikov, D.~Repovs, 
\emph{Embeddings up to homotopy type in Euclidean Space,}
Bull. Austral. Math. Soc (1993). 

\bibitem{Ep}
D.B.A.~Epstein, 
\emph{Embedding punctured manifolds,}
Proc. Amer. Math. Soc, Vol {\bf 16} No. 2 (Apr. 1965), pp. 175--176.

\bibitem{Erle}
D.~Erle,
\emph{Die quadratische form eines knotens, und ein Satz \"uber Knotenmannigfaltigkeiten,}
J. Reine Angew. Math, 236 (1969), 174--218.

\bibitem{SternR}
R.~Fintushel, R.~Stern,
\emph{An exotic free involution on $S^4$,}
Ann. of Math. (2) {\bf 113} (1981) no. 2, 357--365.

\bibitem{FinSter}
R.~Fintushel, R.~Stern,
\emph{Rational homology cobordisms of spherical space forms,}
Topology, {\bf 26} no. 3 pp. 385--393, (1987).

\bibitem{Fickle}
H.~Fickle, 
\emph{Knots, Z-Homology 3-Spheres and Contractible 4-Manifolds,}
pp. 467-493.  Houston Journal of Mathematics Vol. 10, No. 4 (1984).

\bibitem{FQ}
M.~Freedman, F.~Quinn,
\emph{Topology of 4-manifolds,}
Princeton University Press 1990.

\bibitem{GilLiv}
P.M.~Gilmer, C.~Livingston,
\emph{On embedding 3-manifolds in 4-space,}
Topology, {\bf 22}, no. 3, pp. 241--252 (1983).

\bibitem{GS}
R.~Gompf, A.~Stipsicz, 
\emph{4-Manifolds and Kirby Calculus,}
Graduate Studies in Mathematics Vol {\bf 20}, AMS (1999).

\bibitem{Han}
W.~Hantzsche,
\emph{Einlagerung von Mannigflitigkeiten in euklidishe Raume,} Math. Zeit. {\bf 43} (1938),
38--58.

\bibitem{Hillman}
J.A.~Hillman,
\emph{Embedding homology equivalent $3$-manifolds in $4$-space,}
Math. Z. {\bf 223} (1996), no. 3. 473--481.

\bibitem{LinkInv}
J.A.~Hillman,
\emph{Algebraic Invariants of Links,}
Series on Knots and Everything -- Vol. 32.
World Scientific.

\bibitem{Hill2}
J.A.~Hillman, 
\emph{Embedding $3$-manifolds with circle actions in $4$-space,}
University of Sydney preprint,
{\tt [http://www.maths.usyd.edu.au/u/pubs/publist/preprints/2008/hillman-19.html]}
to appear in Proc. AMS.

\bibitem{Hill3}
J.A.~Hillman, 
\emph{Complements of connected hypersurfaces in $S^4$,} 
JKTR (2017). [{\it arXiv:1502.04385}]

\bibitem{Hill4}
J.A.~Hillman, 
\emph{3-manifolds with abelian embeddings in $S^4$,} 
JKTR (2020). [{\tt arXiv:1707.00376}]

\bibitem{Hill5}
J.A.~Hillman, 
\emph{3-manifolds with nilpotent embeddings in $S^4$}, 
[{\tt arXix:1912.03486}] preprint.

\bibitem{Issa}
A.~Issa, D.~McCoy, 
\emph{Smoothly embedding Seifert fibered spaces in $S^4$}, 
[{\tt arXiv:1810.04770}] preprint.

\bibitem{JO}
K.~J\"anich, E.~Ossa, 
\emph{On the signature of an involution,}
Topology {\bf 8} (1969) 27--30.

\bibitem{Kawauchi}
A.~Kawauchi,
\emph{On quadratic forms of $3$-manifolds,}
Invent. Math. {\bf 43} (1977), no. 2 177--198.

\bibitem{KK}
A.~Kawauchi, S.~Kojima,
\emph{Algebraic classification of linking pairings on $3$-manifolds,}
Math. Ann. {\bf 253} (1980), no. 1, 29--42. 

\bibitem{K}
C.~Kearton,
\emph{Signatures of knots and the free differential calculus,}
Quart. J. Math. Oxford (2), {\bf 30} (1979), 157--182.

\bibitem{Kirby}
R.~Kirby, 
\emph{The topology of 4-manifolds,}
Springer LNM 1374.

\bibitem{KA}
\emph{The knot atlas,}
{\tt [http://katlas.org]}

\bibitem{Kosinski}
A.~Kosinski,
\emph{Differential Manifolds,}
Academic Press. Vol {\bf 138} Pure and Applied Mathematics. (1993)
Dover Publications (October 19, 2007).

\bibitem{KM}
P.~Kronheimer, T.~Mrowka, 
\emph{Monopoles and Three-Manifolds,}
Cambridge University Press (2007).

\bibitem{Lisca}
P.~Lisca,
\emph{Sums of lens spaces bounding rational balls,}
Algebraic \& Geometric Topology {\bf 7} (2007) 2141--2164.

\bibitem{Litherland}
R.A.~Litherland,
\emph{Deforming twist-spun knots,}
Trans. Amer. Math. Soc. {\bf 250} (1979), 311--331.

\bibitem{LivingstonSurv}
C.~Livingston,
\emph{A survey of classical knot concordance,}
Handbook of knot theory. 319--347 (2005), Elsevier.

\bibitem{Massey}
W.S.~Massey, 
\emph{Proof of a conjecture of Whitney,}
Pacific J. Math, Vol {\bf 31}, No. 1, 1969. pp 143--156. 

\bibitem{Matv}
S.V.~Matveev, 
\emph{Recognition and tabulation of three-dimensional manifolds,}
Dokl. Akad. Nauk, {\bf 400} (2005) no. 1 26--28.

\bibitem{MJ}
A.~Mijatovic,
\emph{Simplicial structures of knot complements,}
Math. Res. Lett. {\bf 12} (2005), 843--856.

\bibitem{MS}
J.~Milnor, J.~Stasheff,
\emph{Characteristic Classes,}
Princeton University Press (1974).

\bibitem{Mil}
J.~Milnor, 
\emph{Spin structures on manifolds,}
L'enseignement math\'ematique. Vol 9 (1963). pp. 198--203.

\bibitem{NR}
W.~Neumann, F.~Raymond,
\emph{Seifert manifolds, plumbing, $\mu$-invariant and orientation reversing maps,}
Algebraic and geometric topology (Santa Barbara 1977) 163--196. Lecture Notes in Math. {\bf 664}
Springer, 1978.

\bibitem{Orlik}
P.~Orlik,
\emph{Seifert manifolds,}
Lecture notes in mathematics {\bf 291}, Springer (1972).

\bibitem{OzSz}
P.~Oszv\'ath, Z.~Szab\'o, 
\emph{Absolutely graded Floer homologies and intersection forms for four-manifolds
with boundary,} Adv. Math. 173 (2003) 225--254.

\bibitem{Rub}
D.~Ruberman, 
\emph{Seifert surfaces of knots in $S^4$,}
Pacific J. Math {\bf 145} (1990), no. 1, pp. 97--116.

\bibitem{Rub2}
D.~Ruberman,
\emph{Imbedding punctured lens spaces and connected sums,}
Pacific J. Math. {\bf 113} (1984), no. 2, 481--491.

\bibitem{Sav}
N.~Saveliev, 
\emph{Invariants for homology 3-spheres,}
Encyclopedia of Mathematical Sciences {\bf 140}. Springer-Verlag. (2002)

\bibitem{Schubert}
H.~Schubert,
\emph{Knoten und vollringe,}
 Acta Mat. {\bf 90}, 131--286 (1953). 

\bibitem{ST}
H.~Seifert, W.~Threlfall,
\emph{Seifert and Threlfall's Textbook of Topology.}
Academic Press (1980).

\bibitem{Siebenmann}
L.~Siebenmann,
\emph{On vanishing of the Rohlin invariant and nonfinitely amphicheiral homology 3-spheres,}
Topology Symposium, Siegen 1979, 172--222. Lect. Notes in Math. {\bf 788}. Springer, 1980.

\bibitem{Skopenkov}
A.~Skopenkov, 
\emph{Classification of smooth embeddings of $3$-manifolds
in the $6$-space,} Math. Zeitschrift, {\bf 260:3} (2008) 647-672.

\bibitem{Stern}
R.~Stern,
\emph{Some Brieskorn spheres which bound contractible manifolds,}
Notices Amer. Math. Soc {\bf 25} (1978), A448. pp. 313--317.

\bibitem{Larry}
L.~Taylor, 
\emph{Complex $\spin$ structures on 3-manifolds,}
Fields Institute Communications. Vol {\bf 47} (2005) 313--317.

\bibitem{Thurston}
W.~Thurston,
\emph{Three-dimensional manifolds, Kleinian groups and hyperbolic geometry,}
Bull. Amer. Math. Soc. {\bf 6} (1982), 357--381.

\bibitem{Trotter}
H.~Trotter, 
\emph{On s-equivalence of Seifert matrices,}
Invent. Math. {\bf 20}, 173--207. (1973)

\bibitem{Wall}
C.T.C.~Wall, 
\emph{All 3-manifolds imbed in 5-space,}
Bull. Amer. Math. Soc. {\bf 71} (1965) 564--567. 

\bibitem{Whitney}
H.~Whitney,
\emph{On the topology of differentiable manifolds,}
Lectures in topology, Univ. of Michigan Press, 1941, pp. 101--141.


\end{thebibliography}
\end{document}